\documentclass{amsart}
\usepackage{amsfonts, amsbsy, amsmath, amssymb, latexsym}

\newtheorem{thm}{Theorem}[section]
\newtheorem{lem}[thm]{Lemma}
\newtheorem{cor}[thm]{Corollary}
\newtheorem{prop}[thm]{Proposition}
\newtheorem{rem}[thm]{Remark}

\numberwithin{equation}{section}

\begin{document}

\title [Isomorphism Classes of Extensions of $p$-adic Fields]
{Enumeration of Isomorphism Classes of Extensions of $p$-adic Fields}

\author{Xiang-dong Hou and Kevin Keating}

\address{Department of Mathematics and Statistics\\
Wright State University\\
Dayton, Ohio 45435}

\email{xhou@euler.math.wright.edu}

\address{Department of Mathematics\\
University of Florida\\
Gainesville, Florida 32611}

\email{keating@math.ufl.edu}

\keywords{}

\subjclass{}

\begin{abstract}
Let $\Omega$ be an algebraic closure of ${\mathbb Q}_p$ and
let $F$ be a finite extension of ${\mathbb Q}_p$ contained in
$\Omega$.  Given positive integers $f$ and $e$, the number of
extensions $K/F$ contained in $\Omega$ with residue degree
$f$ and ramification index $e$ was computed by Krasner. This
paper is concerned with the number ${\mathfrak I}(F,f,e)$ of
$F$-isomorphism classes of such extensions.  We determine
${\mathfrak I}(F,f,e)$ completely when $p^2\nmid e$ and get
partial results when $p^2\parallel e$.  When $s$ is large,
${\mathfrak I}({\mathbb Q}_p,f,e)$ is equal to the
number of isomorphism classes of finite commutative chain
rings with residue field ${\mathbb F}_{p^f}$, ramification index
$e$, and length $s$.
\end{abstract}

\maketitle

%%%%%%%%%%%%%%%%%%%%%%%%%%%%%%%%%%%%%%%%%%
%    Section 1
%%%%%%%%%%%%%%%%%%%%%%%%%%%%%%%%%%%%%%%%%%%
\section{Introduction} \label{intro}

Fix an algebraic closure $\Omega$ of ${\mathbb Q}_p$ and let $F/{\mathbb Q}_p$ be
a finite extension contained in $\Omega$. Given positive integers $f$ and $e$,
let ${\mathcal E}(F,f,e)$ denote the set of all extensions $K/F$
contained in
$\Omega$ which have residue degree $f$ and ramification index $e$. Krasner's
formulas in \cite{Kra1} -- \cite{Kra5} allow one to compute the cardinality
${\mathfrak N}(F,f,e)$ of the set ${\mathcal E}(F,f,e)$.
Suppose $e=p^me_0$ with
$p\nmid e_0$.  Krasner's formulas state that
\begin{equation} \label{Krasner}
{\mathfrak N}(F,f,e)=e\sum_{s=0}^mp^s(p^{\epsilon(s)N}-
p^{\epsilon(s-1)N}),
\end{equation}
where $N=fe[F:{\mathbb Q}_p]$ and
\begin{equation}
\epsilon(s)=
\begin{cases}
p^{-1}+p^{-2}+\cdots+p^{-s}&\text{if}\ s>0,\cr
0&\text{if}\ s=0,\cr
-\infty&\text{if}\ s=-1.\cr
\end{cases}
\end{equation}
In this paper we consider a related question:
What is the number ${\mathfrak I}(F,f,e)$ of $F$-isomorphism
classes of elements in ${\mathcal E}(F,f,e)$?
Unfortunately, the formulas for ${\mathfrak I}(F,f,e)$ seem to be
much more complicated than those for ${\mathfrak N}(F,f,e)$.

When $p^2\nmid e$, we are able to determine
${\mathfrak I}(F,f,e)$ completely; when $p^2\parallel e$, we
are able to determine ${\mathfrak I}(F,f,e)$ with some
additional assumptions on $f$ and $e$.  It is well-known and
elementary that ${\mathfrak I}(F,f,e)$ can be computed
as a weighted sum over the elements of ${\mathcal E}(F,f,e)$,
\begin{equation} \label{known}
{\mathfrak I}(F,f,e)=\frac{1}{fe}
\sum_{K\in{\mathcal E}(F,e,f)}|\text{Aut}(K/F)|.
\end{equation}
Our method is to use class field theory to determine the
groups $\text{Aut}(K/F)$ explicitly.

Besides Krasner's formulas, another motivation for our work is the
connections between $p$-adic fields and finite commutative
chain rings. A chain ring is a ring whose ideals form a chain
under inclusion.  Finite commutative chain rings have
applications in finite geometry (\cite{Kling}, \cite{Tor}) and
combinatorics (\cite{Hou1}, \cite{Hou2}, \cite{Leu1}, \cite{Leu2}).
Since finite commutative chain rings are precisely the nontrivial
quotients of rings of integers of $p$-adic fields,
classifying isomorphism classes of finite extensions of
${\mathbb Q}_p$
is essentially equivalent to classifying isomorphism classes of finite
commutative chain rings. In particular, in
Section~\ref{rings} we will show that ${\mathfrak I}({\mathbb Q}_p,f,e)$
is equal to the number of isomorphism classes of finite commutative
chain rings with residue field ${\mathbb F}_{p^f}$, ramification index
$e$, and length $s$, for all sufficiently large $s$.

     The paper is organized as follows. Section~\ref{rings} is a summary of
the connections between $p$-adic fields and finite commutative
chain rings.  Section~\ref{prep} contains some preparatory results about
$p$-adic fields.  In particular, we determine the
${\mathbb F}_p[\langle \boldsymbol{\gamma}\rangle]$-module
structure of $K^\times/(K^\times)^p$, where $K$ is a finite
extension of ${\mathbb Q}_p$ and $\boldsymbol{\gamma}$ is a
${\mathbb Q}_p$-automorphism of $K$. In Section~\ref{p0} we consider the
problem of computing ${\mathfrak I}(F,f,e)$ when $p\nmid e$. Besides calculating
${\mathfrak I}(F,f,e)$, we also collect some facts about tamely ramified
extensions of $F$ which will be used later in the paper. In Section~\ref{p1} we determine
${\mathfrak I}(F,f,e)$ in the case $p\parallel e$. Sections \ref{p2}--\ref{conclude} are devoted
to calculating ${\mathfrak I}(F,f,e)$ in the case $p^2\parallel e$, with
some additional restrictions on $f$ and $e$.  In
Section~\ref{p2} we outline the computational plan and
determine the structures of certain Galois groups.
The key ingredients
in the formula for ${\mathfrak I}(F,f,e)$ are computed in Sections \ref{d2}--\ref{d012}, and
the final formula is assembled in Section~\ref{conclude}.

For $K\subset\Omega$ a finite extension of ${\mathbb Q}_p$, we let
$n_K=[K:{\mathbb Q}_p]$ be the degree of $K/{\mathbb Q}_p$.
We denote the ring of integers of $K$ by
${\mathcal O}_K$, the maximal ideal of ${\mathcal O}_K$ by ${\mathcal M}_K$, and
the residue field of $K$ by $\bar{K}={\mathcal O}_K/{\mathcal M}_K$.
Any generator $\pi_K$ for ${\mathcal M}_K$ is called
a uniformizer for $K$.  We let $\nu_K$ denote the valuation
on $K$ normalized so that $\nu_K(\pi_K)=1$ for any
uniformizer $\pi_K$.  Then $\nu_K$ extends uniquely to a
valuation on $\Omega$ which takes values in ${\mathbb Q}$, and is also
denoted $\nu_K$.  In particular, we let $\nu_p=\nu_{{\mathbb Q}_p}$
denote the valuation on $\Omega$ which satisfies
$\nu_p(p)=1$.  Let $L$ be a finite extension of $K$.
Then the residue degree $[\bar{L}:\bar{K}]$ of $L/K$ is denoted
$f(L/K)$, and the ramification index $\nu_L(\pi_K)$ of $L/K$
is denoted $e(L/K)$.  Finally, let
$\{\zeta_a:a\ge 1\}$ be a compatible system of primitive roots
of unity in $\Omega$, with $\zeta_a$ a primitive $a$th root of
unity and $\zeta_{ab}^b=\zeta_a$ for every $a,b\ge 1$.

%%%%%%%%%%%%%%%%%%%%%%%%%%%%%%%%%%%%%%%%%%%%%%%
%   Section 2
%%%%%%%%%%%%%%%%%%%%%%%%%%%%%%%%%%%%%%%%%%%%%%%
\section{$p$-adic Fields and Finite Commutative Chain Rings}
\label{rings}

In addition to the description in terms of $p$-adic fields
given in Section~\ref{intro}, there is another more explicit
construction of finite commutative chain rings based on
Galois rings; we refer the reader to
\cite{McD} for more details.  Choose a prime $p$, positive
integers $n,f$, and a monic polynomial
$\Phi\in({\mathbb Z}/p^n{\mathbb Z})[X]$ of degree
$f$ whose image in $({\mathbb Z}/p{\mathbb Z})[X]$ is irreducible.
The ring $\text{GR}(p^n,f)=({\mathbb Z}/
p^n{\mathbb Z})[X]/(\Phi)$ is called the Galois ring of characteristic
$p^n$ and rank $f$; it is determined up to isomorphism by $p$, $n$,
and $f$.  Every finite commutative chain ring is isomorphic to a
ring of the form $R[X]/(\Psi,p^{n-1}X^t)$,
where $R=\text{GR}(p^n,f)$ is a Galois ring, $\Psi\in R[X]$ is an
Eisenstein polynomial of degree $e$, and
\begin{equation}
\begin{cases}
t=e&\text{if}\ n=1,\cr
1\le t\le e&\text{if}\ n\ge 2.\cr
\end{cases}
\end{equation}
The integers $p,n,f,e,t$ are called the invariants of the finite
commutative chain ring \cite{Clark}.

The following proposition summarizes the connections between
finite commutative chain rings and $p$-adic fields.
%%%%%%%%%%%%%%%%%%%%%%%%%%%%%%%%%%%%
%  Proposition 2.1
%%%%%%%%%%%%%%%%%%%%%%%%%%%%%%%%%%%%
\begin{prop} \label{chain}
Let $K/{\mathbb Q}_p$ be a finite extension, with
residue degree $f$ and ramification index $e$, and let
$k/{\mathbb Q}_p$
be the maximal unramified subextension of $K/{\mathbb Q}_p$.
Let $s,t,n$ be positive integers such that
$s=(n-1)e+t$, with $1\le t\le e$, and let $\tilde{\Psi}$
denote the image of $\Psi\in{\mathcal O}_k[X]$ in
$({\mathcal O}_k/p^n{\mathcal O}_k)[X]$.  Then we have the
following.
\smallskip

(i) Let $a\in{\mathcal O}_k$ be such that $\bar k=({\mathbb Z}/p{\mathbb Z})[\bar a]$,
where $\bar a$ is the image of $a$ in $\bar k$, and let
$\Phi\in{\mathbb Z}_p[X]$ be the minimal polynomial of $a$ over
${\mathbb Q}_p$. Then the image $\tilde \Phi$ of $\Phi$ in
$({\mathbb Z}/p^n{\mathbb Z})[X]$ is monic of degree $f$ and the
image $\bar \Phi$ of $\tilde \Phi$ in $({\mathbb Z}/p{\mathbb Z})[X]$ is
irreducible.  Therefore
\begin{equation}
{\mathcal O}_k/p^n{\mathcal O}_k\cong({\mathbb Z}/p^n{\mathbb Z})[X]/(\tilde \Phi)\cong
\text{\rm GR}(p^n,f).
\end{equation}

(ii) The minimal polynomial
of $\pi_K$ over $k$ is an Eisenstein polynomial
$\Psi\in{\mathcal O}_k[X]$ of degree $e$ such that
\begin{equation}
{\mathcal O}_K/\pi_K^s{\mathcal O}_K\cong({\mathcal O}_k/p^n{\mathcal O}_k)[X]/(\tilde \Psi, p^{n-1}
X^t)\cong\text{\rm GR}(p^n,f)[X]/(\tilde \Psi, p^{n-1}X^t),
\end{equation}
where $\tilde{\Psi}\in({\mathcal O}_k/p^n{\mathcal O}_k)[X]
\cong\text{\rm GR}(p^n,f)[X]$ is an Eisenstein polynomial over
$\text{\rm GR}(p^n,f)$. Thus ${\mathcal O}_K/\pi_K^s{\mathcal O}_K$
is a finite commutative chain ring with invariants
\begin{equation}
\begin{cases}
(p,1,f,t,t)&\text{if}\ n=1,\cr
(p,n,f,e,t)&\text{if}\ n>1.\cr
\end{cases}
\end{equation}
Moreover, every finite commutative chain ring is isomorphic to
${\mathcal O}_K/\pi_K^s{\mathcal O}_K$ for some finite extension $K/{\mathbb Q}_p$ and some
$s\ge1$. \smallskip

(iii) Let $L/{\mathbb Q}_p$ be another finite extension and
let $s\geq\bigl(\frac{p}{p-1}+\nu_p(e)\bigr)e$.  Then
${\mathcal O}_K/\pi_K^s{\mathcal O}_K\cong {\mathcal O}_L/\pi_L^s{\mathcal O}_L$ if and only if $K\cong L$.
\end{prop}

\begin{proof}
(i) and (ii) are well-known.

(iii) We want to prove that if ${\mathcal O}_K/\pi_K^s{\mathcal
O}_K\cong{\mathcal O}_L/\pi_L^s{\mathcal O}_L$
with $s\geq\bigl(\frac{p}{p-1}+\nu_p(e)\bigr)e$ then $K\cong L$.
Note that $s\geq\bigl(\frac{p}{p-1}+\nu_p(e)\bigr)e$ implies
$n>1$.  Thus the residue degree and ramification index of $L/
{\mathbb Q}_p$ are determined by
${\mathcal O}_L/ \pi_L^s {\mathcal O}_L \cong {\mathcal O}_K/\pi_K^s{\mathcal O}_K$,
and so $L/{\mathbb Q}_p$ also has residue degree $f$ and
ramification index $e$.  We may assume that $K$ and $L$ are
both contained in the algebraic
closure $\Omega$ of ${\mathbb Q}_p$. Then $K/{\mathbb Q}_p$ and
$L/{\mathbb Q}_p$ have the same maximal unramified subextension
$k/{\mathbb Q}_p$, and $K/k$ and $L/k$ are both
totally ramified extensions of degree $e$.
We may assume that $e>1$.
Let $\Psi\in{\mathcal O}_k[X]$ be the
minimal polynomial of $\pi_K$ over $k$. The assumption
${\mathcal O}_K/\pi_K^s{\mathcal O}_K\cong
{\mathcal O}_L/\pi_L^s{\mathcal O}_L$ implies that
\begin{equation}
({\mathcal O}_k/p^n{\mathcal O}_k)[X]/(\tilde{\Psi},p^{n-1}X^{t-1})\cong{\mathcal O}_L/
\pi_L^s{\mathcal O}_L.
\end{equation}
By Lemma~XVII.8 in McDonald \cite{McD} there exists $\boldsymbol{\sigma}\in
\text{Aut}({\mathcal O}_k/p^n{\mathcal O}_k)$ such that $\boldsymbol{\sigma}\tilde \Psi$ has a root
$\beta\in {\mathcal O}_L/
\pi_L^s{\mathcal O}_L$. Let $b\in{\mathcal O}_L$ be a lifting of $\beta$.
Since $\boldsymbol{\sigma}\tilde \Psi\in ({\mathcal O}_k/p^n{\mathcal O}_k)[X]$ is an Eisenstein
polynomial of degree $e<s$,
it follows that $\nu_L(b)=1$.  Since $k/{\mathbb Q}_p$ is unramified, the
natural homomorphism
$\text{Gal}(k/{\mathbb Q}_p)\rightarrow\text{Aut}({\mathcal O}_k/p^n{\mathcal O}_k)$
is an isomorphism.  Let
$\boldsymbol{\Sigma}$ be the element of $\text{Gal}(k/{\mathbb Q}_p)$ whose
image in $\text{Aut}({\mathcal O}_k/p^n{\mathcal O}_k)$ is $\boldsymbol{\sigma}$.
Then $\boldsymbol{\sigma}\tilde \Psi=\widetilde{\boldsymbol{\Sigma} \Psi}$,
so the image of $(\boldsymbol{\Sigma} \Psi)(b)$ in
${\mathcal O}_L/\pi_L^s{\mathcal O}_L$ is $(\widetilde{\boldsymbol{\Sigma} \Psi})(\beta)=
(\boldsymbol{\sigma}\tilde \Psi)(\beta)=0$. Therefore
$\nu_L((\boldsymbol{\Sigma} \Psi)(b))\geq s$.

Let $r_1,r_2,\dots,r_e\in\Omega$ be the roots of
$\boldsymbol{\Sigma} \Psi$; then
$(\boldsymbol{\Sigma} \Psi)(X)=\prod_{i=1}^e(X-r_i)$.
We may order the $r_i$ so that
$m=\nu_p(b-r_1)$ is as large as possible.  Then for $2\leq i\leq e$
we have $\nu_p(b-r_i)\geq\min\{m,\nu_p(r_1-r_i)\}$.
If $m>\nu_p(r_1-r_i)$ then $\nu_p(b-r_i)=\nu_p(r_1-r_i)$, while
if $m\leq\nu_p(r_1-r_i)$ then by the maximality of $m$ we get
$\nu_p(b-r_i)\leq\nu_p(r_1-r_i)$.  It follows that for
$2\leq i\leq e$ we have $\nu_p(b-r_i)\leq\nu_p(r_1-r_i)$.  Since
$(\boldsymbol{\Sigma} \Psi)(b)=(b-r_1)(b-r_2)\dots(b-r_e)$, this
implies
\begin{equation}
\frac{s}{e}\leq\nu_p((\boldsymbol{\Sigma} \Psi)(b))\leq m+\sum_{i=2}^e\,\nu_p(r_1-r_i)=
m+\nu_p(\delta_{k(r_1)/k}),
\end{equation}
where $\delta_{k(r_1)/k}=(\boldsymbol{\Sigma}\Psi)'(r_1)$ is the different of the extension
$k(r_1)/k$.  By Remark~1 on p.\,58 of \cite{Ser} we have
$\nu_p(\delta_{K/k})\leq1-e^{-1}+\nu_p(e)$.  Since we are
assuming ${s\geq\bigl(\frac{p}{p-1}+\nu_p(e)\bigr)e}$, this
implies $m>\frac{1}{p-1}$.  Using Lemma~\ref{lower} below we get
${\frac{1}{p-1}\geq\nu_p(r_1-r_i)}$ and hence
$m>\nu_p(r_1-r_i)$ for all $2\leq i\leq e$.  It
follows by Krasner's lemma (see \cite{Kra2}, p.\,224) that
$k(b)\supset k(r_1)$.  Since $[k(b):k]=[k(r_1):k]=e$, we get
$L=k(b)=k(r_1)\cong K$.
\end{proof}

\begin{lem} \label{lower}
Let $k$ be a finite extension of ${\mathbb Q}_p$, let $\Psi\in k[X]$ be
an Eisenstein polynomial of degree $e$, and let
$r_1,r_2,\dots,r_e$ be the roots of $\Psi$.  Then for every
$2\leq i\leq n$ we have $\nu_p(r_1-r_i)\leq\frac{1}{p-1}$.
\end{lem}

\begin{proof}
Let $E=k(r_1)$.  The lemma may be rephrased as a statement
about the higher
ramification theory of the extension $E/k$, which need not be
Galois; for the ramification theory of non-Galois
extensions, see for instance III\,\S3 of \cite{FV}, or the
appendix to \cite{Del}.  In fact the integers $\nu_E(r_1-r_i)$
are the lower ramification breaks for the extension $E/k$.  The
lemma is equivalent to the statement that these breaks are bounded
above by $\frac{1}{p-1}\cdot\nu_E(p)$.  Our method is to
reduce to the case of a Galois extension, where the lemma is
well-known (see for instance Exercise 3(c) on p.\,72 of
\cite{Ser}).

Let $F\subset\Omega$ be the splitting field of $\Psi$, and
set $G=\text{Gal}(F/k)$,
$H=\text{Gal}(F/E)$.  Let $D=F^{G_1}$ be the fixed field of the
wild ramification subgroup of $G$, let $e_1$ be the ramification
index of $D/k$, and let $e_2$ be the ramification index of $ED/E$.
Then $p\nmid e_1$, and hence $p\nmid e_2$.
It follows that the Hasse-Herbrand functions for the
extensions $D/k$ and $ED/E$ are given by
$\phi_{D/k}(x)=x/e_1$ and $\phi_{ED/E}(x)=x/e_2$
for $x\geq0$.  Using the composition rule for towers of
extensions we get $\phi_{E/k}(x)=\frac{1}{e_1}\phi_{ED/D}(e_2x)$.
Since the largest lower ramification break of $E/k$ is
$\text{inf}\{x:\phi_{E/k}'(x)=1/e\}$,
it suffices to prove the lemma for the extension $ED/D$.
Since $G_1=\text{Gal}(F/D)$ is a $p$-group,
there is a refinement $G_1=G_1^{(0)}\geq G_1^{(1)}\geq G_1^{(2)}\geq\ldots\geq
G_1^{(n-1)}\geq G_1^{(n)}=\{1\}$ of the ramification filtration of $G_1$
such that $G_1\trianglerighteq G_1^{(i)}$ and $|G_1^{(i)}/G_1^{(i+1)}|=p$
for all $0\leq i\leq n-1$.  Let $j$ be the largest integer such
that $G_1^{(j)}$ is not contained in $H$, and let $D'$ be the
subfield
of $F$ fixed by $G_1^{(j)}(H\cap G_1)$.  Then the largest lower
ramification break of $ED/D$ is the same as the largest
lower ramification break of $ED/D'$.  Since $|G_1^{(j)}(H\cap G_1)/
(H\cap G_1)|=p$, we have $H\cap G_1\trianglelefteq
G_1^{(j)}(H\cap G_1)$.  Therefore the extension $ED/D'$ is
Galois.  It follows that the lemma holds for $ED/D'$, and
hence also for $E/k$.
\end{proof}

Let ${\mathfrak C}(p,n,f,e,t)$ be the number of isomorphism
classes of finite commutative chain rings with invariants
$(p,n,f,e,t)$. Then Proposition~\ref{chain}(iii) implies that
\begin{equation} \label{classes}
{\mathfrak C}(p,n,f,e,t)={\mathfrak I}({\mathbb Q}_p,f,e)
\text{ when }(n-1)e+t\geq\Bigl(\frac{p}{p-1}+\nu_p(e)\Bigr)e.
\end{equation}
When $p\nmid e$, the number ${\mathfrak C}(p,n,f,e,t)$ was first
determined by Clark and Liang \cite{Clark}. A different
formula for this quantity was given in \cite{Hou3}.
%%%%%%%%%%%%%%%%%%%%%%%%%
%  Theorem 2.2
%%%%%%%%%%%%%%%%%%%%%%%%%
\begin{thm} (Clark and Liang \cite{Clark} and Hou \cite{Hou3})
\label{clh}
Let $p$ be a prime and let $n,f,e,t$ be positive integers
such that $n\ge 2$, $1\le t\le e$, and $p\nmid e$. Then
\begin{equation}
{\mathfrak C}(p,n,f,e,t)=
\sum_{c\mid(e,p^f-1)}\frac{\phi(c)}{\tau(c)}=
\frac 1 f\sum_{i=0}^{f-1}\,\bigl(p^{(i,f)}-1,e\bigr),
\end{equation}
where $\phi$ is the Euler function, $(a,b)$ is the greatest
common divisor of $a$ and $b$, and
$\tau(c)$ is the smallest positive integer $m$ such that
$p^m\equiv1\pmod c$.
\end{thm}

 From Proposition~\ref{chain}(iii) and Theorem~\ref{clh} it
follows that when $p\nmid e$,
\begin{equation}
{\mathfrak I}({\mathbb Q}_p,f,e)=
\sum_{c\mid(e,p^f-1)}\frac{\phi(c)}{\tau(c)}=
\frac 1f\sum_{i=0}^{f-1}\,\bigl(p^{(i,f)}-1,e\bigr).
\end{equation}
In Section~\ref{p0}, we derive a third formula for
${\mathfrak I}({\mathbb Q}_p,f,e)$ in the case $p\nmid e$.
In the other direction,
our formulas for ${\mathfrak I}(F,f,e)$ allow us to compute
${\mathfrak C}(p,n,f,e,t)$ in the following two cases (cf.\
(\ref{classes}), Theorem~\ref{main1}, and Theorem~\ref{main2}):
\smallskip \\
(i) $p\parallel e$ and $n\ge 3+\frac{1}{p-1}-\frac{t}e$;
\smallskip \\
(ii) $p>2$, $p^2\parallel e$,
$n\ge 4+\frac{1}{p-1}-\frac{t}e$, and $(p^f-1,e)=1$.

%%%%%%%%%%%%%%%%%%%%%%%%%%%%%%%%%%%%%%%%%%%%%%%%%%
%    Section 3
%%%%%%%%%%%%%%%%%%%%%%%%%%%%%%%%%%%%%%%%%%%%%%%%%%
\section{Preparatory Results about $p$-adic Fields} \label{prep}

%%%%%%%%%%%%%%%%%%%%%%%%%%%%
%  Proposition 3.1
%%%%%%%%%%%%%%%%%%%%%%%%%%%%
\begin{prop} \label{unique1}
Let $F\subset K$ be finite extensions of ${\mathbb Q}_p$
such that $K/F$ is totally ramified
of degree $p^is$, with $p\nmid s$. Then for each positive integer $d\mid s$,
there is a unique field $K_d$ such that $F\subset K_d\subset K$ and
$[K_d:F]=d$. Furthermore, for $d_1\mid s$ and $d_2\mid s$, we
have $K_{d_1}\subset K_{d_2}$ if and only if $d_1\mid d_2$.
\end{prop}

\begin{proof}
Let $L/F$ be the Galois closure of $K/F$, and set $G=\text{Gal}(L/F)$ and
$H=\text{Gal}(L/K)$.  Let $G_0$ be the inertia subgroup and
$G_1$ the wild inertia subgroup of $G$ (so $G_1$ is
the unique Sylow $p$-subgroup of $G_0$).  Then $G_0/G_1$ is a cyclic
group whose order is prime to $p$ and divisible by $s$.
Since $|G/H|=p^is$ factors as the product of
$|G_1H/H|=|G_1/(G_1\cap H)|$, which is a power of $p$, and
$|G/G_1H|=|G_0/(G_0\cap G_1H)|$, which is prime to $p$, we
have $|G_0/(G_0\cap G_1H)|=s$.
Let $N=G_0\cap G_1H$.  Then $N$ is the unique subgroup of $G_0$ of index $s$ which
contains $G_1$, so $N\trianglelefteq G$.
Since $K/F$ is a totally ramified extension we have $G_0H=G$.
Therefore $G_1H/N$ maps isomorphically onto $G/G_0$,
and hence $G/N$ is a semidirect product of $G_1H/N\cong G/G_0$
acting on $G_0/N$.  This implies that for each $d\mid s$
there is a unique subgroup $S_d\leq G$ of index $d$ such that
$S_d\geq G_1H$.  The fixed field of $S_d$ acting on $L$ is $K_d$.
\end{proof}
%%%%%%%%%%%%%%%%%%%%%%%%%%%%%%%%%%%%%
%   Proposition 3.2
%%%%%%%%%%%%%%%%%%%%%%%%%%%%%%%%%%%%%%
\begin{prop} \label{unique2}
Let $F$ be a finite extension of ${\mathbb Q}_p$ and set
$f(F/{\mathbb Q}_p) =f_0$.  Let $e=p^is$ with $p\nmid s$, and
let $K\in{\mathcal E}(F,f,e)$. \smallskip

(i) There is a unique field $L_K$ such that $F\subset L_K\subset K$ and $K/L_K$ is
totally ramified of degree $p^i$. \smallskip

(ii) If $\bigl(p^{f_0f}-1,s\bigr)=1$, there is a unique field $E_K$
such that $F\subset E_K\subset K$ and $E_K/F$ is totally
ramified of degree $s$. Moreover, we have $E_K\subset L_K$ and
$\text{\rm Aut}(K/F)=\text{\rm Aut}(K/E_K)$.
\end{prop}

\begin{proof}
(i) This is a special case of Proposition~\ref{unique1}.

(ii) Let $k/F$ be the maximal unramified subextension of $K/F$.
Since $p\nmid s$, there are uniformizers $\pi_{L_K}$ for $L_K$ and
$\pi_F$ for $F$ such that $\pi_{L_K}^s/\pi_F\in{\mathcal O}_k^{\times}$
(see II, Prop.\,3.5 in \cite{FV}).
Since $\bigl((p^{f_0f}-1)\cdot p,s\bigr)=1$, we have
$\pi_{L_K}^s/\pi_F=\beta^s$ for some $\beta\in{\mathcal O}_k^{\times}$.  Then
$(\pi_{L_K}/\beta)^s=\pi_F$, and hence $E_K=F(\pi_{L_K}/\beta)$
is a totally ramified extension of $F$ of degree $s$ which is
contained in $L_K$.  To prove the uniqueness of $E_K$, assume that
we have $F\subset E\subset K$ with $E/F$ totally ramified of
degree $s$.
Then there is a uniformizer $\pi_E$
for $E$ such that $\pi_E^s/\pi_F\in{\mathcal O}_k^\times$.
As above we get $\delta\in{\mathcal O}_k^{\times}$
with $\delta^s=\pi_E^s/\pi_F$ and $(\pi_E/\delta)^s=\pi_F=(\pi_{L_K}/\beta)^s$.
Thus $\pi_E/\delta=\zeta\pi_{L_K}/\beta$ for some $\zeta\in L_K$ with
$\zeta^s=1$.  Since $\bigl((p^{f_0f}-1)\cdot p,s\bigr)=1$ we must
have $\zeta=1$, and hence $E=k(\pi_E/\delta)=k(\pi/\beta)=E_K$.

To prove the last statement, we note that for any
$\boldsymbol{\sigma}\in\text{Aut}(K/F)$ we have
$\boldsymbol{\sigma}(E_K)=E_K$ by the uniqueness of $E_K$.
Thus $\boldsymbol{\sigma}|_{E_K}\in\text{Aut}(E_K/F)$. We have
already seen that $E_K=F(\pi_{E_K})$ for some $\pi_{E_K}\in E_K$
such that $\pi_{E_K}^s$ is a uniformizer for $F$. Since $F$ does
not contain any nontrivial $s$th root of unity, it follows that
$\text{Aut}(E_K/F)=\{\text{id}\}$. Thus $\boldsymbol{\sigma}\in
\text{Aut}(K/E_K)$.
\end{proof}

Let $K/{\mathbb Q}_p$ be a finite extension with $f(K/{\mathbb Q}_p)=f$, $e(K/{\mathbb Q}_p)=e$, and let
$\boldsymbol{\gamma}\in\text{Aut}(K/{\mathbb Q}_p)$. Let $K^{\boldsymbol{\gamma}}$
denote the subfield of $K$ fixed by $\langle\boldsymbol{\gamma}\rangle$ and put
$e(\boldsymbol{\gamma})=e(K/K^{\boldsymbol{\gamma}})$ and $f(\boldsymbol{\gamma})=f(K/K^{\boldsymbol{\gamma}})$.
The elementary abelian $p$-group $V_K=K^\times/(K^\times)^p$
can be viewed as a module over the group ring ${\mathbb F}_p[\langle \boldsymbol{\gamma}\rangle]$, where
${\mathbb F}_p={\mathbb Z}/p{\mathbb Z}$. Our approach in Section~\ref{p1} depends on knowing the
${\mathbb F}_p[\langle \boldsymbol{\gamma}\rangle]$-module structure of $V_K$ in the case
$p\nmid e(\boldsymbol{\gamma})$. The rest of this section is devoted to the determination of this
${\mathbb F}_p[\langle \boldsymbol{\gamma}\rangle]$-module structure. The first step is to factor $K^\times$ as
\begin{equation} \label{factorK}
K^\times\cong\langle \pi_K\rangle\times\langle\zeta_{p^f-1}\rangle\times U_1,
\end{equation}
where $U_1=1+\pi_K{\mathcal O}_K$.  We have then
\begin{equation} \label{factorVK}
V_K\cong(\langle\pi_K\rangle/\langle\pi_K^p\rangle)\times(U_1/U_1^p).
\end{equation}

Let $F/K^{\boldsymbol{\gamma}}$ be the maximal unramified subextension
of $K/K^{\boldsymbol{\gamma}}$; then $F$ is the fixed field of
$\langle\boldsymbol{\gamma}^{f(\boldsymbol{\gamma})}\rangle$,
and $K/F$ is a totally ramified cyclic extension of
degree $e(\boldsymbol{\gamma})$.  Since $p\nmid
e(\boldsymbol{\gamma})$, the field $F$ contains an
$e(\boldsymbol{\gamma})$th root of unity, and there is a
uniformizer $\pi_K$ for $K$ such that
$\pi_K^{e(\boldsymbol{\gamma})}\in F$.  In fact we have
$\pi_K^{e(\boldsymbol{\gamma})}=\eta\pi$ for some
$\eta\in{\mathcal O}_F^{\times}$ and
$\pi\in K^{\boldsymbol{\gamma}}$.  Since ${\mathcal O}_F^{\times}/
({\mathcal O}_F^{\times})^{e(\boldsymbol{\gamma})}$ is generated by roots of unity,
we may assume that $\eta$ is a root of unity, of order prime
to $p$.  It follows that
$\boldsymbol{\gamma}(\pi_K^{e(\boldsymbol{\gamma})})=\eta^{q-1}\pi_K^{e(\boldsymbol{\gamma})}$, where $q$
is the cardinality of the residue field of
$K^{\boldsymbol{\gamma}}$.
Since $\eta^{q-1}\in(K^{\times})^p$, this implies that $\boldsymbol{\gamma}$ acts
trivially on the image of $\pi_K^{e(\boldsymbol{\gamma})}$ in
$V_K$, and hence also on the image of $\pi_K$ in $V_K$.  Since $U_1$ is
clearly stabilized by $\boldsymbol{\gamma}$, this implies
that the factors in (\ref{factorVK}) are
${\mathbb F}_p[\langle\boldsymbol{\gamma}\rangle]$-submodules of $V_K$.
Thus it remains only to determine the
${\mathbb F}_p[\langle\boldsymbol{\gamma}\rangle]$-module
structure of $U_1/U_1^p$.

%%%%%%%%%%%%%%%%%%%%%%%%%%%%%%%%%%%%
%   Proposition 3.3
%%%%%%%%%%%%%%%%%%%%%%%%%%%%%%%%%%%%
\begin{prop} \label{U1}
(i) If $\zeta_p\notin K$, there is an isomorphism of ${\mathbb Z}_p[\langle\boldsymbol{\gamma}\rangle]$-modules
\begin{equation}
U_1\cong{\mathbb Z}_p[\langle\boldsymbol{\gamma}\rangle]^{n_K/o(\boldsymbol{\gamma})}.
\end{equation}
Hence there is an isomorphism of ${\mathbb F}_p[\langle\boldsymbol{\gamma}\rangle]$-modules
\begin{equation}
U_1/U_1^p\cong{\mathbb F}_p[\langle\boldsymbol{\gamma}\rangle]^{n_K/o(\boldsymbol{\gamma})}.
\end{equation}
(ii) If $\zeta_p\in K$ and $p\nmid e(\boldsymbol{\gamma})$, there is an isomorphism of ${\mathbb F}_p[\langle\boldsymbol{\gamma}\rangle]$-modules
\begin{equation}
U_1/U_1^p\cong\langle\zeta_p\rangle\times{\mathbb F}_p[\langle\boldsymbol{\gamma}\rangle]^{n_K/o(\boldsymbol{\gamma})}.
\end{equation}
\end{prop}

\begin{proof}[Proof of Proposition~\ref{U1}(i)]
We use Theorem~4(b) in Weiss \cite{Wei}. This theorem implies
that as long as $U_1$ contains no roots of unity, the ${\mathbb Z}_p[\langle\boldsymbol{\gamma}\rangle]$-isomorphism
class of $U_1$ is determined by its ${\mathbb Z}_p$-rank. Since $\langle \boldsymbol{\gamma}\rangle$ is cyclic, this
means that we can replace the extension $K/K^{\boldsymbol{\gamma}}$ with an extension
$K'/k'$ such that $K'/k'$ and $k'/{\mathbb Q}_p$ are unramified,
$n_{k'}=n_{K^{\boldsymbol{\gamma}}}$, and
$\text{Gal}(K'/k')\cong\text{Gal}(K/K^{\boldsymbol{\gamma}})=\langle\boldsymbol{\gamma}\rangle$.
The logarithm gives an isomorphism between the
${\mathbb Z}_p[\text{Gal}(K'/k')]$-modules
$U_1'=1+p{\mathcal O}_{K'}$ and $p{\mathcal O}_{K'}$.  Since $K'/k'$ is unramified,
$p{\mathcal O}_{K'}$ is free over ${\mathbb Z}_p[\text{Gal}(K'/k')]$ of rank
$n_{K'}/o(\boldsymbol{\gamma})$.  It follows that $U_1'$
is free over ${\mathbb Z}_p[\text{Gal}(K'/k')]$ of rank
$n_{K'}/o(\boldsymbol{\gamma})$, and hence that $U_1$ is
free over ${\mathbb Z}_p[\langle\boldsymbol{\gamma}\rangle]$
of rank $n_K/o(\boldsymbol{\gamma})$.
\end{proof}

The ${\mathbb Z}_p[\langle\boldsymbol{\gamma}\rangle]$-module structure of $U_1$
cannot be described as simply when $U_1$ contains roots of unity.
In fact, when $p\nmid e(\boldsymbol{\gamma})$ it follows from another theorem
of Gruenberg and Weiss (Theorem~6.1(a) of \cite{Gru}) that $U_1$
is cohomologically trivial as a
${\mathbb Z}_p[\langle\boldsymbol{\gamma}\rangle]$-module.
Since the ${\mathbb Z}_p$-torsion subgroup
$\langle\zeta_{p^s}\rangle$ of $U_1$ is not in general
cohomologically trivial, $\langle\zeta_{p^s}\rangle$ need not be a
direct ${\mathbb Z}_p[\langle\boldsymbol{\gamma}\rangle]$-summand
of $U_1$.  Therefore to prove Proposition~\ref{U1}(ii),
we work directly with $U_1/U_1^p$.

\begin{proof}[Proof of Proposition~\ref{U1}(ii)]
The group $U_1$ has a filtration $U_1\supset U_2\supset\dots$,
where $U_i=1+\pi_K^i{\mathcal O}_K$.  This induces
a filtration on $U_1/U_1^p$ whose $i$th filtrant is
$\bar U_i=U_iU_1^p/U_1^p$. Put $r=pe/(p-1)$. Since we are assuming
$\zeta_p\in K$, we have $(p-1)\mid e$ and hence $r\in{\mathbb Z}$. Define
\begin{equation}
I=\{i\in {\mathbb Z}:1\le i\le r\ \text{and}\ p\nmid i\}.
\end{equation}
Our strategy is to first determine the
${\mathbb F}_p[\langle\boldsymbol{\gamma}\rangle]$-module structure of the quotients
$\bar U_i/\bar U_{i+1}$ for $i\in I$, and then use this information to
reconstruct $U_1/U_1^p$.
There are isomorphisms
of ${\mathbb F}_p[\langle\boldsymbol{\gamma}\rangle]$-modules
\begin{equation}
\begin{split}
\bar U_i/\bar U_{i+1}&\cong U_iU_1^p/U_{i+1}U_1^p \\
&\cong U_i/(U_i\cap(U_{i+1}U_1^p)).
\end{split}
\end{equation}
For $i\in I$ we have
$U_i\cap(U_{i+1}U_1^p)=U_{i+1}$, while if $i\not\in I$ and $i\not=r$
we have $U_{i+1}U_1^p\supset U_i$.  Therefore if $i\in I$ we have
\begin{equation} \label{layer}
\bar U_i/\bar U_{i+1}\cong U_i/U_{i+1}\cong\pi_K^i{\mathcal
O}_K/\pi_K^{i+1}{\mathcal O}_K,
\end{equation}
while if $i\not\in I$ and $i\not=r$ we have $\bar U_i/\bar
U_{i+1}=\{1\}$.
Finally, we have $|\bar{U}_r/\bar{U}_{r+1}|=p$
since $|U_1/U_1^p|=p^{fe+1}$, $|I|=e$, and $U_{r+1}\subset
U_1^p$.

%%%%%%%%%%%%%%%%%%%%%%%%%%%%%%%%%%%%%%%%%%%%%%%%%
%  Lemma 3.4
%%%%%%%%%%%%%%%%%%%%%%%%%%%%%%%%%%%%%%%%%%%%%%%%%
\begin{lem} \label{indep}
Let $i\geq0$ and let $\hat{\boldsymbol{\gamma}}$ denote the
automorphism of $\pi_K^i{\mathcal O}_K/\pi_K^{i+1}{\mathcal O}_K$
induced
by $\boldsymbol{\gamma}$. Then $\hat{\boldsymbol{\gamma}}^0,\hat{\boldsymbol{\gamma}}^1,\dots,\hat{\boldsymbol{\gamma}}^{f(\boldsymbol{\gamma})-1}$ are linearly independent over $\bar{K}$.
\end{lem}

\begin{proof}
If not then there is a monic polynomial
\begin{equation}
P(X)=X^m+a_{m-1}X^{m-1}+\cdots+a_1X+a_0
\end{equation}
of degree $m<f(\boldsymbol{\gamma})$ in $\bar{K}[X]$ such that $P(\hat{\boldsymbol{\gamma}})=0$.
We assume that $m$ is as small as
possible; then $a_0\ne0$, since $\hat{\boldsymbol{\gamma}}$ is invertible.
Since $1\le m<f(\boldsymbol{\gamma})$, there exists $\alpha\in \bar{K}$ such that
$\boldsymbol{\gamma}^m(\alpha)\ne\alpha$.  Since $P(\hat{\boldsymbol{\gamma}})\cdot v=0$
for all $v\in\pi_K^i{\mathcal O}_K/\pi_K^{i+1}{\mathcal O}_K$, we have
$P(\hat{\boldsymbol{\gamma}})\cdot\alpha v=0$ for all $v$ as well.  This implies that
\begin{equation}
Q(X)=\boldsymbol{\gamma}^m(\alpha)X^m+a_{m-1}\boldsymbol{\gamma}^{m-1}(\alpha)X^{m-1}+
\cdots+a_1\boldsymbol{\gamma}(\alpha)X+a_0\alpha
\end{equation}
satisfies $Q(\hat{\boldsymbol{\gamma}})=0$.  Therefore
$R(X)=Q(X)-\boldsymbol{\gamma}^m(\alpha)P(X)\in \bar{K}[X]$ is
a polynomial of degree $<m$ with nonzero constant term such
that $R(\hat{\boldsymbol{\gamma}})=0$.  This violates the
minimality of $m$, and therefore proves the lemma.
\end{proof}

     It follows from Lemma~\ref{indep} that
$\hat{\boldsymbol{\gamma}}^0,\hat{\boldsymbol{\gamma}}^1,\dots,\hat{\boldsymbol{\gamma}}^{f(\boldsymbol{\gamma})-1}$
are linearly independent over $\bar{K}^{\boldsymbol{\gamma}}$, and hence that the
degree of the minimal polynomial of $\hat{\boldsymbol{\gamma}}$
over $\bar{K}^{\boldsymbol{\gamma}}$
is $\geq f(\boldsymbol{\gamma})$.  Since
$\pi_K^i{\mathcal O}_K/\pi_K^{i+1}{\mathcal O}_K$ has
dimension $f(\boldsymbol{\gamma})$ over
$\bar{K}^{\boldsymbol{\gamma}}$, this implies that
$\pi_K^i{\mathcal O}_K/\pi_K^{i+1}{\mathcal O}_K$ is a cyclic
$\bar{K}^{\boldsymbol{\gamma}}[\langle\boldsymbol{\gamma}\rangle]$-module generated by some
$v_i$, i.e.,
\begin{equation} \label{cyclic}
\pi_K^i{\mathcal O}_K/\pi_K^{i+1}{\mathcal O}_K=\bar{K}^{\boldsymbol{\gamma}}[\langle\boldsymbol{\gamma}\rangle]\cdot v_i.
\end{equation}
Since $\boldsymbol{\gamma}^{f(\boldsymbol{\gamma})}$ generates the
inertia group of the tamely ramified extension
$K/K^{\boldsymbol{\gamma}}$, the image $\xi$ of
$\boldsymbol{\gamma}^{f(\boldsymbol{\gamma})}(\pi_K)/\pi_K$
in $\bar{K}$ is a primitive $e(\boldsymbol{\gamma})$th root
of unity.  On the other hand,
%%%%%%%%%%%%%%%%%%%%%%%%%%%%%%%%%%%%%%%%%%%%%%%%%%%%%%%%%%%%%%%%%%
% since $K/F$ is Galois and totally ramified of
% degree $e(\boldsymbol{\gamma})$, there is a primitive
% $e(\boldsymbol{\gamma})$th root of unity in $F$, and hence
% in $K^{\boldsymbol{\gamma}}$.
%%%%%%%%%%%%%%%%%%%%%%%%%%%%%%%%%%%%%%%%%%%%%%%%%%%%%%%%%%%%%%%%%%
class field theory gives an onto homomorphism $\rho:(K^\gamma)^\times\rightarrow
\langle \gamma\rangle$ such that the inertia subgroup $\langle \gamma^{f(\gamma)}\rangle$
of $\langle\gamma\rangle$ is the image of the unit group of $K^\gamma$.
Hence $(K^\gamma)^\times$ contains an element of order $e(\gamma)$.
Therefore $\bar{K}^{\boldsymbol{\gamma}}$ contains a primitive
$e(\boldsymbol{\gamma})$th root of
unity, so we have $\xi\in \bar{K}^{\boldsymbol{\gamma}}$.
In particular, $\bar{K}^{\boldsymbol{\gamma}}\cdot v_i$ is a
$\bar{K}^{\boldsymbol{\gamma}}[\langle \boldsymbol{\gamma}^{f(\boldsymbol{\gamma})}\rangle]$-submodule of
$\pi_K^i{\mathcal O}_K/\pi_K^{i+1}{\mathcal O}_K$.

Let $\bar\pi_K^i$ denote the image of $\pi_K^i$ in
$\pi_K^i{\mathcal O}_K/\pi_K^{i+1}{\mathcal O}_K$.  Then
$\hat{\boldsymbol{\gamma}}^{f(\boldsymbol{\gamma})}(\bar\pi_K^i)=\xi^i\bar \pi_K^i$.
Since $\langle\boldsymbol{\gamma}^{f(\boldsymbol{\gamma})}\rangle$
acts trivially on $\bar{K}^{\boldsymbol{\gamma}}$,
it follows that there is a
$\bar{K}^{\boldsymbol{\gamma}}[\langle\boldsymbol{\gamma}^{f(\boldsymbol{\gamma})}\rangle]$-module isomorphism
\begin{equation} \label{sum}
\bar{K}^{\boldsymbol{\gamma}}[\langle\boldsymbol{\gamma}^{f(\boldsymbol{\gamma})}\rangle]\cong\bigoplus_{i=0}^{e(\boldsymbol{\gamma})-1}\bar{K}^{\boldsymbol{\gamma}}\cdot v_i.
\end{equation}
Therefore $\bar{K}^{\boldsymbol{\gamma}}\cdot v_i$ is a projective
$\bar{K}^{\boldsymbol{\gamma}}[\langle\boldsymbol{\gamma}^{f(\boldsymbol{\gamma})}\rangle]$-module.
Using (\ref{cyclic}) we get
\begin{equation} \label{induce}
\pi_K^i{\mathcal O}_K/\pi_K^{i+1}{\mathcal O}_K
=\bar{K}^{\boldsymbol{\gamma}}[\langle\boldsymbol{\gamma}\rangle]\otimes_{\bar{K}^{\boldsymbol{\gamma}}[\langle\boldsymbol{\gamma}^{f(\boldsymbol{\gamma})}\rangle]}
\bar{K}^{\boldsymbol{\gamma}}\cdot v_i.
\end{equation}
It follows that $\pi_K^i
{\mathcal O}_K/\pi_K^{i+1}{\mathcal O}_K$ is projective over $\bar{K}^{\boldsymbol{\gamma}}[\langle\boldsymbol{\gamma}\rangle]$,
and hence also over ${\mathbb F}_p[\langle\boldsymbol{\gamma}\rangle]$.
Using (\ref{layer}) we get an
${\mathbb F}_p[\langle\boldsymbol{\gamma}\rangle]$-module
isomorphism
\begin{equation} \label{sumU1}
U_1/U_1^p\cong\biggl(\bigoplus_{i\in I}\pi_K^i{\mathcal O}/\pi_K^{i+1}{\mathcal O}_K\biggr)\oplus
\bar U_r.
\end{equation}

Write $e=dp^t$ with $p\nmid d$; then $e(\boldsymbol{\gamma})\mid d$ and $(p-1)\mid d$. Partition $I$
into $I\cap[0]$, $I\cap[1],\dots,I\cap[d-1]$, where $[i]$ is the congruence class of $i$
modulo $d$. Then each subset $I\cap [i]$ contains either $\frac ed$, $\frac ed-1$, or
$\frac ed +1$ elements. However, using the fact that $(p-1)\mid d$, one can show that
$|I\cap[i]|=\frac ed\pm 1$ if and only if $|I\cap[pi]|=\frac ed\mp 1$.
If $i\equiv j\pmod{d}$ then there is an ${\mathbb
F}_p[\langle\boldsymbol{\gamma}^{f(\boldsymbol{\gamma})}\rangle]$-module isomorphism
$\bar{K}^{\boldsymbol{\gamma}}\cdot v_i\cong
\bar{K}^{\boldsymbol{\gamma}}\cdot v_j$.  On the other hand,
if $v_{pi}$ is chosen suitably,
$x\mapsto x^p$ gives an ${\mathbb F}_p[\langle\boldsymbol{\gamma}^{f(\boldsymbol{\gamma})}\rangle]$-module isomorphism
between $\bar{K}^{\boldsymbol{\gamma}}\cdot v_i$ and $\bar{K}^{\boldsymbol{\gamma}}\cdot v_{pi}$. Therefore we have an
${\mathbb F}_p[\langle\boldsymbol{\gamma}^{f(\boldsymbol{\gamma})}\rangle]$-module
isomorphism
\begin{equation}
\bigoplus_{i\in I}\bar{K}^{\boldsymbol{\gamma}}\cdot v_i\cong\Biggl(\bigoplus_{i=0}^{d-1}\bar{K}^{\boldsymbol{\gamma}}\cdot v_i\Biggr)^{e/d}.
\end{equation}
Using (\ref{sum}) we get
\begin{equation} \label{big1}
\begin{split}
\bigoplus_{i\in I}\bar{K}^{\boldsymbol{\gamma}}\cdot v_i
&\cong\Biggl(\bigoplus_{i=0}^{e(\boldsymbol{\gamma})-1}\bar{K}^{\boldsymbol{\gamma}}\cdot v_i\Biggr)^{e/e(\boldsymbol{\gamma})}\cr
&\cong \bar{K}^{\boldsymbol{\gamma}}[\langle\boldsymbol{\gamma}^{f(\boldsymbol{\gamma})}\rangle]^{e/e(\boldsymbol{\gamma})}\quad \cr
&\cong{\mathbb F}_p[\langle\boldsymbol{\gamma}^{f(\boldsymbol{\gamma})}\rangle]^{n_K/o(\boldsymbol{\gamma})}.
\end{split}
\end{equation}
Therefore by (\ref{induce}) and (\ref{big1}) there are
${\mathbb F}_p[\langle \boldsymbol{\gamma}\rangle]$-module
isomorphisms
\begin{equation} \label{big2}
\begin{split}
\bigoplus_{i\in I}\pi_K^i{\mathcal O}_K/\pi_K^{i+1}{\mathcal O}_K&\cong\bigoplus_{i\in I}
\left(\bar{K}^{\boldsymbol{\gamma}}[\langle\boldsymbol{\gamma}\rangle]\otimes_{\bar{K}^{\boldsymbol{\gamma}}[\langle\boldsymbol{\gamma}^{f(\boldsymbol{\gamma})}\rangle]}\bar{K}^{\boldsymbol{\gamma}}\cdot v_i\right)\cr
&\cong\bigoplus_{i\in I}\left({\mathbb F}_p[\langle\boldsymbol{\gamma}\rangle]\otimes_{{\mathbb F}_p[\langle\boldsymbol{\gamma}^{f(\boldsymbol{\gamma})}\rangle]}
\bar{K}^{\boldsymbol{\gamma}}\cdot v_i\right)\cr
&\cong{\mathbb F}_p[\langle\boldsymbol{\gamma}\rangle]\otimes_{{\mathbb F}_p[\langle\boldsymbol{\gamma}^{f(\boldsymbol{\gamma})}\rangle]}
\biggl(\bigoplus_{i\in I}\bar{K}^{\boldsymbol{\gamma}}\cdot
v_i\biggr)\cr
&\cong{\mathbb F}_p[\langle\boldsymbol{\gamma}\rangle]\otimes_{{\mathbb F}_p[\langle\boldsymbol{\gamma}^{f(\boldsymbol{\gamma})}\rangle]}
{\mathbb F}_p[\langle\boldsymbol{\gamma}^{f(\boldsymbol{\gamma})}\rangle]^{n_K/o(\boldsymbol{\gamma})}\cr
&\cong{\mathbb F}_p[\langle\boldsymbol{\gamma}\rangle]^{n_K/o(\boldsymbol{\gamma})}.\cr
\end{split}
\end{equation}
To complete the proof of Proposition~\ref{U1}(ii), we need only determine the ${\mathbb F}_p[\langle\boldsymbol{\gamma}\rangle]$-module
structure of $\bar U_r$.

\begin{lem} \label{Ur}
There is an ${\mathbb F}_p[\langle\boldsymbol{\gamma}\rangle]$-module isomorphism $\bar U_r\cong\langle\zeta_p\rangle$.
\end{lem}

\begin{proof}
Define a group homomorphism $\psi:U_r\rightarrow{\mathbb F}_p$ by setting
\begin{equation}
\psi(1+p(\zeta_p-1)x)=\text{Tr}_{\bar{K}/{\mathbb F}_p}(\bar x),
\end{equation}
where $x\in{\mathcal O}_K$ and $\bar x$ is the image of $x$ in $\bar{K}$. We claim that $\psi(U_r\cap(U_{r+1}U_1^p))=0$.  If $y\in U_r\cap(U_{r+1}U_1^p)$, then
$y\equiv(1+(\zeta_p-1)x)^p\pmod{\pi_K^{r+1}}$ for some $z\in{\mathcal O}_K$. We have
\begin{equation}
(1+(\zeta_p-1)z)^p\equiv 1+p(\zeta_p-1)\left(z+\frac{(\zeta_p-1)^{p-1}}p z^p\right)
\pmod{\pi_K^{r+1}}.
\end{equation}
Since $(\zeta_p-1)^{p-1}\equiv -p\pmod {p\pi_K}$, we have
\begin{equation}
\text{Tr}_{\bar{K}/{\mathbb F}_p}\biggl(\bar z+\overline{\frac{(\zeta_p-1)^{p-1}}p z^p}\biggr)=
\text{Tr}_{\bar{K}/{\mathbb F}_p}(\bar z-\bar z^p)=0,
\end{equation}
so $\psi(y)=0$. Since $\psi$ is nontrivial, it
follows that $\psi$ induces a group
isomorphism between $\bar U_r\cong U_r/(U_r\cap(U_{r+1}U_1^p))$ and ${\mathbb F}_p$.
Suppose that $\boldsymbol{\gamma}(\zeta_p)=\zeta_p^m$. Then
\begin{equation}
\begin{split}
\boldsymbol{\gamma}(1+p(\zeta_p-1)z)&=1+p(\zeta_p^m-1)\boldsymbol{\gamma}(z)\cr
&=1+p(\zeta_p-1)(1+\cdots+\zeta_p^{m-1})\boldsymbol{\gamma}(z)\cr
&\equiv 1+p(\zeta_p-1)m\boldsymbol{\gamma}(z)\pmod{\pi_K^{r+1}}.
\end{split}
\end{equation}
Since $\text{Tr}_{\bar{K}/{\mathbb F}_p}(\overline{m\boldsymbol{\gamma}(z)})=m\text{Tr}_{\bar{K}/{\mathbb F}_p}
(\bar z)$, we see that $\boldsymbol{\gamma}$ acts on $\bar U_r\cong{\mathbb F}_p$
by raising to the power $m$.
Therefore the ${\mathbb F}_p[\langle\boldsymbol{\gamma}\rangle]$-modules
$\bar U_r$ and $\langle\zeta_p\rangle$ are isomorphic.
\end{proof}

Using (\ref{sumU1}), (\ref{big2}), and Lemma~\ref{Ur}, we get an ${\mathbb F}_p[\langle\boldsymbol{\gamma}\rangle]$-module
isomorphism
\begin{equation}
U_1/U_1^p\cong\langle\zeta_p\rangle\times{\mathbb F}_p[\langle\boldsymbol{\gamma}\rangle]^{n_K/o(\boldsymbol{\gamma})},
\end{equation}
and the proof of Proposition~\ref{U1}(ii) is complete.
\end{proof}
%%%%%%%%%%%%%%%%%%%%%%%%%%%%%%%%%%%%%%%%%%%%
%  Corollary 3.6
%%%%%%%%%%%%%%%%%%%%%%%%%%%%%%%%%%%%%%%%%%

\begin{cor} \label{VK}
(i) If $\zeta_p\notin K$, there is an isomorphism of ${\mathbb F}_p[\langle \boldsymbol{\gamma}\rangle]$-modules
\begin{equation}
K^\times/(K^\times)^p\cong{\mathbb F}_p\times{\mathbb F}_p[\langle\boldsymbol{\gamma}\rangle]^{n_K/o(\boldsymbol{\gamma})}.
\end{equation}

(ii) If $\zeta_p\in K$ and $p\nmid e(\boldsymbol{\gamma})$, there is an isomorphism of ${\mathbb F}_p[\langle\boldsymbol{\gamma}\rangle]$-modules
\begin{equation}
K^\times/(K^\times)^p\cong{\mathbb F}_p\times\langle\zeta_p\rangle\times
{\mathbb F}_p[\langle\boldsymbol{\gamma}\rangle]^{n_K/o(\boldsymbol{\gamma})}.
\end{equation}
\end{cor}

%%%%%%%%%%%%%%%%%%%%%%%%%%%%%%%%%%%%%%%%%%%%%%%%%%
%      Section 4
%%%%%%%%%%%%%%%%%%%%%%%%%%%%%%%%%%%%%%%%%%%%%%%%%%
\section{The Case $p\nmid e$} \label{p0}

In this section we determine ${\mathfrak I}(F,f,e)$ in
the cases where $p\nmid e$.  Thus we are
restricting our attention to tamely ramified extensions of the $p$-adic
field $F$, which are in general well-understood.  Therefore we only
give outlines for the arguments in this section, most of
which are not new.  Besides calculating ${\mathfrak I}(F,f,e)$, we
also
collect some facts about tamely ramified extensions of $F$ which will
be useful in the next section.

We start by listing the elements of ${\mathcal E}(F,f,e)$.
Let $q$ denote the cardinality of the residue field of $F$
and set $g=(q^f-1)e$.  Let $\pi_F$ be a uniformizer of $F$, let
$\pi_E\in \Omega$ be an $e$th root of $\pi_F$, and define
$E=F(\zeta_g,\pi_E)$.  The extension $E/F$ is Galois, with
\begin{equation} \label{GalEF}
\text{Gal}(E/F)=\langle \boldsymbol{\sigma},\boldsymbol{\tau} : \boldsymbol{\sigma}^e=
\boldsymbol{\tau}^{df}=1,\;\boldsymbol{\tau}\boldsymbol{\sigma}\boldsymbol{\tau}^{-1}=
\boldsymbol{\sigma}^q\rangle,
\end{equation}
where $df=[F(\zeta_g):F]$.
The actions of $\boldsymbol{\sigma}$ and $\boldsymbol{\tau}$ on $E$ are given by
\begin{equation}
\boldsymbol{\sigma}(\pi_E)=\zeta_e\pi_E,\ \boldsymbol{\sigma}(\zeta_g)=\zeta_g,\
\boldsymbol{\tau}(\pi_E)=\pi_E,\ \boldsymbol{\tau}(\zeta_g)=\zeta_g^q.
\end{equation}
For $0\le h<e$, we define $K_h=F(\zeta_{q^f-1},\pi_h)$, where
$\pi_h=\zeta_g^h\pi_E$; equivalently, $K_h$ is the subfield of
$E$ fixed by $\langle\boldsymbol{\sigma}^{-h}\boldsymbol{\tau}^f\rangle$.
Then we have
\begin{equation} \label{EO}
{\mathcal E}(F,f,e)=\{K_h:0\le h<e\}.
\end{equation}

We now describe the $F$-automorphisms of $K_h$. The inertia subgroup
$\text{Aut}(K_h/F(\zeta_{q^f-1}))$ of $\text{Aut}(K_h/F)$ is
a cyclic group of order $b=(e,q^f-1)$ generated by an element
${\boldsymbol\mu}$ such that ${\boldsymbol\mu}(\zeta_{q^f-1})=\zeta_{q^f-1}$
and ${\boldsymbol\mu}(\pi_h)=\zeta_b\pi_h$. We need to determine which elements of $\text{Gal}(F(\zeta_{q^f-1})/F)$ can be extended to
automorphisms of $K_h$. Let $\boldsymbol{\rho}$ be the Frobenius automorphism
of $F(\zeta_{q^f-1})/F$. For $c\ge 0$ we attempt to extend
$\boldsymbol{\rho}^c$ to an element $\boldsymbol{\nu}_c$ of $\text{Aut}(K_h/F)$.
Since $\pi_h^g=\pi_E^g=\pi_F^{q^f-1} \in F$,
we have $\boldsymbol{\nu}_c(\pi_h)=\epsilon\pi_h$ for some
$\epsilon\in K_h$ such that $\epsilon^g=1$. But since $p\nmid g$, we must have
$\epsilon=\zeta_{q^f-1}^x$ for some $x\in {\mathbb Z}$.  Therefore
\begin{equation}
\zeta_{q^f-1}^{ex+h}\pi_F=(\zeta_{q^f-1}^x\pi_h)^e=\boldsymbol{\nu}_c(\pi_h)^e= \boldsymbol{\nu}_c(\pi_h^e)=
\boldsymbol{\nu}_c(\zeta_{q^f-1}^h\pi_F)=\zeta_{q^f-1}^{hq^c}\pi_F,
\end{equation}
which implies that
\begin{equation} \label{extend}
ex\equiv(q^c-1)h\pmod{q^f-1}.
\end{equation}
Conversely, if $x$ satisfies (\ref{extend}), then $\boldsymbol{\rho}^c$ can be extended
to $\boldsymbol{\nu}_c\in\text{Aut}(K_h/F)$ such that $\boldsymbol{\nu}_c
(\pi_h)=\zeta_{q^f-1}^x\pi_h$.

The congruence (\ref{extend}) can be solved for $x$ if and
only if $b=(e,q^f-1)$ divides $(q^c-1)h$. Let
$c_h$ be the smallest positive integer $c$ satisfying this
condition.  Then $c_h$ is the order of $q$
in $\bigl({\mathbb Z}/\frac b{(b,h)}{\mathbb Z}\bigr)^\times$.
Let $u\in{\mathbb Z}$ satisfy $eu\equiv b\pmod{q^f-1}$
and set $m_h=\frac{1}{b}(q^{c_h}-1)hu$.  Then
$x=m_h$ is a solution to (\ref{extend}) with $c=c_h$.

The elements in $\text{Gal}(F(\zeta_{q^f-1})/F)$ which can be
extended
to automorphisms of $K_h/F$ are $(\boldsymbol{\rho}^{c_h})^i$ for $0\le i<f/c_h$.
Let ${\boldsymbol\nu}\in\text{Aut}(K_h/F)$ be the extension of $\boldsymbol{\rho}^{c_h}$
defined by ${\boldsymbol\nu}(\pi_h)=\zeta_{q^f-1}^{m_h}\pi_h$. Then ${\boldsymbol\nu}$ generates the group
\begin{equation}
\text{Aut}(K_h/F)/\text{Gal}
(F(\zeta_{q^f-1})/F)=\text{Aut}(K_h/F)/\langle{\boldsymbol\mu}\rangle
\end{equation}
and satisfies ${\boldsymbol\nu}^{f/c_h}\in\langle{\boldsymbol\mu}\rangle$.
The actions of ${\boldsymbol\mu}$ and ${\boldsymbol\nu}$ on $K_h=
F(\zeta_{q^f-1},\pi_h)$ are given by
\begin{equation} \label{action}
{\boldsymbol\mu}(\zeta_{q^f-1})=\zeta_{q^f-1},\
{\boldsymbol\mu}(\pi_h)=\zeta_b\pi_h,\
{\boldsymbol\nu}(\zeta_{q^f-1})=\zeta_{q^f-1}^{q^{c_h}},\
{\boldsymbol\nu}(\pi_h)=\zeta_{q^f-1}^{m_h}\pi_h.
\end{equation}
Using (\ref{action}) we find that ${\boldsymbol\nu}^{f/c_h}={\boldsymbol\mu}^{uh}$
and ${\boldsymbol\nu}{\boldsymbol\mu}{\boldsymbol\nu}^{-1}={\boldsymbol\mu}^{q^{c_h}}$.
Therefore
\begin{equation} \label{Aut}
\text{Aut}(K_h/F)=\langle{\boldsymbol\mu},{\boldsymbol\nu}:{\boldsymbol\mu}^b=1,\;
{\boldsymbol\nu}^{f/c_h}={\boldsymbol\mu}^{uh},\;
{\boldsymbol\nu}{\boldsymbol\mu}{\boldsymbol\nu}^{-1}={\boldsymbol\mu}^{q^{c_h}}
\rangle.
\end{equation}
In particular, $|\text{Aut}(K_h/F)|=bf/c_h$. Combining this fact with
(\ref{known}), we get the following proposition.

%%%%%%%%%%%%%%%%%%%%%%%%%%%%%%%%%%%%%%%
%  Proposition 4.1
%%%%%%%%%%%%%%%%%%%%%%%%%%%%%%%%%%%%%%
\begin{prop} \label{I0}
Assume that $p\nmid e$. Then
\begin{equation}
{\mathfrak I}(F,f,e)=\frac be\sum_{h=0}^{e-1}\frac 1{c_h},
\end{equation}
where $b=(e,q^f-1)$, $q$ is the cardinality of the residue
field of $F$, and $c_h$ is the smallest positive integer such
that $b$ divides $(q^{c_h}-1)h$.

\end{prop}

\begin{rem}
\rm
The method of Corollary~4.3 in \cite{Hou3} gives the alternative
formula
\begin{equation}
{\mathfrak I}(F,f,e)=\frac{1}{f}\sum_{i=0}^{f-1}\,
\bigl(q^{(i,f)}-1,e\bigr).
\end{equation}
On the other hand, Proposition~\ref{I0} with $F={\mathbb Q}_p$
gives a third
formula (in addition to those in Theorem~\ref{clh}) for ${\mathfrak C}(p,n,f,e,t)$ when $p\nmid e$ and $n\ge 2$,
\begin{equation}
{\mathfrak C}(p,n,f,e,t)=\frac be\sum_{h=0}^{e-1}\frac 1{c_h}.
\end{equation}
\end{rem}

\begin{rem} \label{future}
\rm
For future use we note that ${\boldsymbol\mu}$ is the restriction to $K_h$
of $\boldsymbol{\sigma}^{e/b}$ and ${\boldsymbol\nu}$ is the restriction to
$K_h$ of $\boldsymbol{\sigma}^{v_h}\boldsymbol{\tau}^{c_h}$, where
\begin{equation}
v_h=\frac{em_h+h(1-q^{c_h})}{q^f-1}.
\end{equation}
\end{rem}

For each $\boldsymbol{\gamma}\in\text{Aut}(K_h/F)$, let $K_h^{\boldsymbol{\gamma}}$ denote
the subfield of $K_h$ fixed by $\langle\boldsymbol{\gamma}\rangle$.  By
(\ref{Aut}), we can write $\boldsymbol{\gamma}$ uniquely in the form
$\boldsymbol{\gamma}={\boldsymbol\mu}^i{\boldsymbol\nu}^j$ with $0\le i<b$ and $0\le j<f/c_h$.
The smallest power of $\boldsymbol{\gamma}$ which lies in the inertia subgroup
$\langle {\boldsymbol\mu}\rangle$ of $\text{Aut}(K_h/F)$ is
\begin{equation}
\boldsymbol{\gamma}^{\frac{f}{(f,c_hj)}}=({\boldsymbol\mu}^i
{\boldsymbol\nu}^j)^{\frac{f}{(f,c_hj)}}
={\boldsymbol\mu}^{t(\boldsymbol{\gamma})},
\end{equation}
where $t(\boldsymbol{\gamma})$ is computed using (\ref{Aut}) to be
\begin{equation}
t(\boldsymbol{\gamma})=\frac{q^{\text{lcm}(f,c_hj)}-1}{q^{c_hj}-1}\cdot i+\frac{uhc_hj}{(f,c_hj)}.
\end{equation}
This implies that the extension $K_h/K_h^{\boldsymbol{\gamma}}$ has residue degree
$f(\boldsymbol{\gamma})=f/(f,c_hj)$ and
ramification index
$e(\boldsymbol{\gamma})=b/(b,t(\boldsymbol{\gamma}))$. The order $o(\boldsymbol{\gamma})=[K_h:K_h^{\boldsymbol{\gamma}}]$ of $\boldsymbol{\gamma}$ is given by
\begin{equation}
o(\boldsymbol{\gamma})=e(\boldsymbol{\gamma})\cdot f(\boldsymbol{\gamma})
=\frac{b}{(b,t(\boldsymbol{\gamma}))}\cdot\frac{f}{(f,c_hj)}.
\end{equation}
%
%%%%%%%%%%%%%%%%%%%%%%%%%%%%%%%%%%%%%%%%%%%%%%%%%%
%    Section 5
%%%%%%%%%%%%%%%%%%%%%%%%%%%%%%%%%%%%%%%%%%%%%%%%%%
\section{The Case $p\parallel e$} \label{p1}

In this section, we assume $e=pe_0$ with $p\nmid e_0$.
We use the notation
of Section~\ref{p0} with $e_0$ in place of $e$. In particular,
\begin{equation}
\left\{ \begin{split}
%\begin{cases}
q&=\text{the cardinality of the residue class field of}\ F,\cr
g&=(q^f-1)e_0,\cr
b&=(e_0,q^f-1),\cr
c_h&=\text{the smallest positive integer such that $b\mid(q^{c_h}-1)h$},\cr
u&\in{\mathbb Z}\ \text{satisfies}\ e_0u\equiv b\pmod{q^f-1},\cr
\displaystyle m_h&=\frac{(q^{c_h}-1)h}b u,\cr
\displaystyle v_h&=\frac{e_0m_h+h(1-q^{c_h})}{q^f-1}.\cr
\end{split} \right.
%\end{cases}
\end{equation}
Let $L\in {\mathcal E}(F,f,e)$.  Then by
Proposition~\ref{unique1}, there is a unique
$K\in {\mathcal E}(F,f,e_0)$ which is contained in $L$.
It follows from (\ref{known}) and (\ref{EO}) that
\begin{equation} \label{sumsum}
\begin{split}
{\mathfrak I}(F,f,e)&=\frac 1{fe}\sum_{K\in{\mathcal E}(F,f,e_0)}\sum_{L\in{\mathcal E}
(K,1,p)}|\text{Aut}(L/F)|\cr
&=\frac 1{fe}\sum_{h=0}^{e_0-1}\sum_{L\in{\mathcal E}(K_h,1,p)}|\text{Aut}
(L/F)|.\cr
\end{split}
\end{equation}
For the time being, we fix $K=K_h$ and concentrate on evaluating the inner sum of (\ref{sumsum}).

Let $L\in{\mathcal E}(K,1,p)$.  Then since $K\in{\mathcal E}(F,f,e_0)$
is uniquely determined by $L$, restriction induces a homomorphism
$\text{Aut}(L/F)\rightarrow\text{Aut}(K/F)$.
Let $H_L\subset\text{Aut}(K/F)$ be the image of this homomorphism. Then
$\text{Aut}(L/F)$ is an extension of $H_L$ by $\text{Aut}(L/K)$.
Thus if $L/K$
is not Galois then $\text{Aut}(K/F)\cong H_L$, while if $L/K$ is Galois then
$\text{Aut}(K/F)$ is an extension of $H_L$ by a cyclic group of order $p$.
For each $\boldsymbol{\gamma}\in\text{Aut}(K/F)$, let
\begin{align}
S_1(\boldsymbol{\gamma})&=\{L\in{\mathcal E}(K,1,p):\boldsymbol{\gamma}\in H_L\ \text{and}\
L/K\ \text{is Galois}\}, \\
S_2(\boldsymbol{\gamma})&=\{L\in{\mathcal E}(K,1,p):\boldsymbol{\gamma}\in H_L\ \text{and}\
L/K\ \text{is not Galois}\},
\end{align}
and define $m_i(\boldsymbol{\gamma})=|S_i(\boldsymbol{\gamma})|$.
Then by counting the elements in the set
\begin{equation}
\{(L,\delta):L\in{\mathcal E}(K,1,p),\;\delta\in\text{Aut}(L/F)\}
\end{equation}
in two different ways we find that
\begin{equation} \label{different}
\sum_{L\in{\mathcal E}(K,1,p)}|\text{Aut}(L/F)|=\sum_{\boldsymbol{\gamma}\in\text{Aut}(K/F)}
\bigl(p\cdot m_1(\boldsymbol{\gamma})+m_2(\boldsymbol{\gamma})\bigr).
\end{equation}
%

%%%%%%%%%%%%%%%%%%%%%%%%%%%%%%%%%%%%%
%  Lemma 5.1
%%%%%%%%%%%%%%%%%%%%%%%%%%%%%%%%%%%%%
\begin{lem} \label{n1lemma}
Let $\boldsymbol{\gamma}\in\text{\rm Aut}(K/F)$, let
$K^{\boldsymbol{\gamma}}$ be the subfield of $K/F$
fixed by $\langle \boldsymbol{\gamma}\rangle$, and let
$d(\boldsymbol{\gamma})=(o(\boldsymbol{\gamma}),p-1)$.  Then
\begin{equation}
m_1(\boldsymbol{\gamma})=
\begin{cases}
\displaystyle(d(\boldsymbol{\gamma})-1)\frac{p^{\frac{n_Fe_0f}{o(\boldsymbol{\gamma})}}-1}{p-1}+
\frac{p^{\frac{n_Fe_0f}{o(\boldsymbol{\gamma})}+1}-1}{p-1}-1&\text{if}\ \zeta_p\notin K,\\[.4cm]
\displaystyle(d(\boldsymbol{\gamma})-2)\frac{p^{\frac{n_Fe_0f}{o(\boldsymbol{\gamma})}}-1}{p-1}+2\cdot
\frac{p^{\frac{n_Fe_0f}{o(\boldsymbol{\gamma})}+1}-1}{p-1}-1&\text{if}\ \zeta_p\in K\setminus K^{\boldsymbol{\gamma}}, \\[.4cm]
\displaystyle(d(\boldsymbol{\gamma})-1)\frac{p^{\frac{n_Fe_0f}{o(\boldsymbol{\gamma})}}-1}{p-1}+
\frac{p^{\frac{n_Fe_0f}{o(\boldsymbol{\gamma})}+2}-1}{p-1}-1&\text{if}\ \zeta_p\in K^{\boldsymbol{\gamma}}.
\end{cases}
\end{equation}
\end{lem}

\begin{proof}
By class field theory, cyclic extensions $L/K$ of degree $p$ such that
$\boldsymbol{\gamma}\in H_L$ correspond to $\langle\boldsymbol{\gamma}\rangle$-invariant subgroups
of $V_K=K^\times/(K^\times)^p$ of index $p$. The unramified degree-$p$
extension of $K$ is excluded from $S_1(\boldsymbol{\gamma})$, and so $m_1(\boldsymbol{\gamma})$ is
one less than the number of $\langle\boldsymbol{\gamma}\rangle$-invariant subgroups of
$V_K$ of index $p$.

For each of the $d(\boldsymbol{\gamma})$ homomorphisms
$\psi: \langle\boldsymbol{\gamma}\rangle\rightarrow
{\mathbb F}_p^\times$, the largest quotient on which
$\langle\boldsymbol{\gamma}\rangle$ acts through $\psi$ is
$V_K(\psi)=V_K/(\boldsymbol{\gamma}-\psi(\boldsymbol{\gamma}))V_K$.
On the other hand, if $H$ is a
$\langle\boldsymbol{\gamma}\rangle$-invariant subgroup
of $V_K$ of index $p$ then
$\langle\boldsymbol{\gamma}\rangle$ acts on $V_K/H$ through
some homomorphism $\psi:\langle\boldsymbol{\gamma}\rangle\rightarrow{\mathbb F}_p^\times$.
Thus $H\ge(\boldsymbol{\gamma}-\psi(\boldsymbol{\gamma}))V_K$ and $H/(\boldsymbol{\gamma}-\psi(\boldsymbol{\gamma}))V_K$
is a subgroup of $V_K(\psi)$ of index $p$.  In fact, $H\leftrightarrow
H/(\boldsymbol{\gamma}-\psi(\boldsymbol{\gamma}))V_K$ gives a one-to-one correspondence between
the set of $\langle\boldsymbol{\gamma}\rangle$-invariant subgroups of $V_K$ of index
$p$ and the set
\begin{equation}
\bigcup_{\psi:\langle\boldsymbol{\gamma}\rangle\rightarrow{\mathbb
F}_p^\times}\{U\subset V_K(\psi):[V_K(\psi):U]=p\}.
\end{equation}
If $\zeta_p\notin K$, it follows from Corollary~\ref{VK}(i) that
there is an isomorphism of ${\mathbb
F}_p[\langle\boldsymbol{\gamma}\rangle]$-modules
$V_K\cong{\mathbb F}_p\times{\mathbb F}_p[\langle\boldsymbol{\gamma}\rangle]^{\frac{n_Fe_0f}{o(\boldsymbol{\gamma})}}$.
Thus $\dim_{{\mathbb F}_p}V_K(\psi)
=\frac{1}{o(\boldsymbol{\gamma})}n_Fe_0f$ for $\psi\ne 1$, and $\dim_{{\mathbb F}_p}V_K(1)=\frac
{1}{o(\boldsymbol{\gamma})}n_Fe_0f+1$. The formula for
$m_1(\boldsymbol{\gamma})$ in the case $\zeta_p\not\in K$
follows from this.

Now assume $\zeta_p\in K$.  By Corollary~\ref{VK}(ii) we have
$V_K\cong{\mathbb F}_p\times\langle\zeta_p\rangle\times
{\mathbb F}_p[\langle\boldsymbol{\gamma}\rangle]^{\frac{n_Fe_0f}{o(\boldsymbol{\gamma})}}$.
If $\zeta_p\notin K^{\boldsymbol{\gamma}}$ then $\boldsymbol{\gamma}(\zeta_p)=\zeta_p^m$
for some $m\not\equiv 1\pmod p$.  Therefore $\dim_{{\mathbb F}_p}V_K(\psi)=
\frac{1}{o(\boldsymbol{\gamma})}n_Fe_0f$ when $\psi(\boldsymbol{\gamma})\not\in\{1,m\}$,
and $\dim_{{\mathbb F}_p}V_K(\psi)=\frac{1}{o(\boldsymbol{\gamma})}n_Fe_0f+1$
when $\psi(\boldsymbol{\gamma})\in\{1,m\}$.
In the case where $\zeta_p\in K^{\boldsymbol{\gamma}}$ we have
$\dim_{{\mathbb F}_p}V_K(\psi)=\frac{1}{o(\boldsymbol{\gamma})}n_Fe_0f$
when $\psi(\boldsymbol{\gamma})\ne 1$ and
$\dim_{{\mathbb F}_p}V_K(\psi)=
\frac{1}{o(\boldsymbol{\gamma})}n_Fe_0f+2$ when
$\psi(\boldsymbol{\gamma})=1$.  The remaining formulas for
$m_1(\boldsymbol{\gamma})$ follow from these observations.
\end{proof}

%%%%%%%%%%%%%%%%%%%%%%%%%%%%%%%%%%%%
%   Lemma 5.2
%%%%%%%%%%%%%%%%%%%%%%%%%%%%%%%%%%%%
\begin{lem} \label{n2lemma}
We have
\begin{equation} \label{n2}
m_2(\boldsymbol{\gamma})=
\begin{cases}
\displaystyle p^{\frac{n_Fe_0f}{o(\boldsymbol{\gamma})}+2}-p^2+p-pd(\boldsymbol{\gamma})\cdot
\frac{p^{\frac{n_Fe_0f}{o(\boldsymbol{\gamma})}}-1}{p-1}&\text{if}\ \zeta_p\notin K,\\[.4cm]
\displaystyle (p^2-p)p^{\frac{n_Fe_0f}{o(\boldsymbol{\gamma})}}-p^2+p-pd(\boldsymbol{\gamma})\cdot
\frac{p^{\frac{n_Fe_0f}{o(\boldsymbol{\gamma})}}-1}{p-1}&\text{if}\ \zeta_p\in K.\cr
\end{cases}
\end{equation}
\end{lem}

\begin{proof}
Suppose $L\in S_2(\boldsymbol{\gamma})$. Since $\text{Aut}(L/F)\cong H_L$,
there is a unique $\tilde{\boldsymbol{\gamma}}\in\text{Aut}(L/F)$ which
extends $\boldsymbol{\gamma}$. Let $L^{\tilde{\boldsymbol{\gamma}}}$
be the subfield of $L$ fixed by $\langle\tilde{\boldsymbol{\gamma}}\rangle$.
Then $L^{\tilde{\boldsymbol{\gamma}}}/K^{\boldsymbol{\gamma}}$ is a ramified extension
of degree $p$ such that $L^{\tilde{\boldsymbol{\gamma}}}K=L$.
Conversely, let $M/K^{\boldsymbol{\gamma}}$ be a ramified extension of degree $p$ such that
the compositum $MK$ is not Galois over $K$. Then $M$ and $K$ are linearly
disjoint over $K^{\boldsymbol{\gamma}}$ and so $\boldsymbol{\gamma}\in\text{Gal}(K/K^{\boldsymbol{\gamma}})$ can be
uniquely extended to an element $\tilde{\boldsymbol{\gamma}}\in\text{Aut}(MK/M)$. Therefore
$MK\in S_2(\boldsymbol{\gamma})$. Thus $M\leftrightarrow MK$ gives a bijection
between the set $\{M\in{\mathcal E}(K^{\boldsymbol{\gamma}},1,p):MK/K\ \text{not Galois}\}$ and
$S_2(\boldsymbol{\gamma})$.  Let
\begin{align}
{\mathcal Y}(\boldsymbol{\gamma})&=\{M\in{\mathcal E}(K^{\boldsymbol{\gamma}},1,p): M/K^{\boldsymbol{\gamma}}\ \text{is Galois}\} \\
{\mathcal Z}(\boldsymbol{\gamma})&=\{M\in{\mathcal E}(K^{\boldsymbol{\gamma}},1,p): M/K^{\boldsymbol{\gamma}}\ \text{is not Galois
but $MK/K$ is Galois}\}.
\end{align}
Then we have
\begin{equation} \label{n2formula}
m_2(\boldsymbol{\gamma})=|{\mathcal E}(K^{\boldsymbol{\gamma}},1,p)|-|{\mathcal Y}(\boldsymbol{\gamma})|-|{\mathcal Z}(\boldsymbol{\gamma})|.
\end{equation}
By Krasner's formula (\ref{Krasner}) we have
\begin{equation}
|{\mathcal E}(K^{\boldsymbol{\gamma}},1,p)|=p^{\frac{n_Fe_0f}{o(\boldsymbol{\gamma})}+2}-p^2+p,
\end{equation}
and by class field theory we have
\begin{equation}
|{\mathcal Y}(\boldsymbol{\gamma})|=
\begin{cases}
\displaystyle\frac{p^{\frac{n_Fe_0f}{o(\boldsymbol{\gamma})}+1}-1}{p-1}-1&\text{if}\ \zeta_p\notin K^{\boldsymbol{\gamma}},\\[.4cm]
\displaystyle\frac{p^{\frac{n_Fe_0f}{o(\boldsymbol{\gamma})}+2}-1}{p-1}-1&\text{if}\ \zeta_p\in K^{\boldsymbol{\gamma}}.\cr
\end{cases}
\end{equation}

It remains to determine $|{\mathcal Z}(\boldsymbol{\gamma})|$. Let $M\in
{\mathcal Z}(\boldsymbol{\gamma})$.
Then $\boldsymbol{\gamma}\in\text{Gal}(K/K^{\boldsymbol{\gamma}})$ lifts to
$\tilde{\boldsymbol{\gamma}}\in\text{Gal}(MK/M)$, so $MK/K^{\boldsymbol{\gamma}}$ is
Galois and $\text{Gal}(MK/K^{\boldsymbol{\gamma}})$ is the semidirect product of
$\text{Gal}(MK/M)=\langle\tilde{\boldsymbol{\gamma}}\rangle$
acting on $\text{Gal}(MK/K)\cong{\mathbb Z}/p{\mathbb Z}$.  This action is nontrivial
since $M/K^{\boldsymbol{\gamma}}$ is not Galois.  On the other hand, suppose that $L/K$
is a cyclic extension of degree $p$ such that $L/K^{\boldsymbol{\gamma}}$ is Galois and
$\text{Gal}(L/K^{\boldsymbol{\gamma}})$ is nonabelian.  Then $\boldsymbol{\gamma}$
can be lifted to an automorphism $\tilde{\boldsymbol{\gamma}}$ of $L$, which
must satisfy $o(\tilde{\boldsymbol{\gamma}})=o(\boldsymbol{\gamma})$, since otherwise
$\text{Gal}(L/K^{\boldsymbol{\gamma}})=\langle\tilde{\boldsymbol{\gamma}}\rangle$ is abelian.
Therefore $\text{Gal}(L/K^{\boldsymbol{\gamma}})$ is a semidirect product of
$\langle\tilde{\boldsymbol{\gamma}}\rangle$ acting nontrivially on $\text{Gal}(L/K)$.
For such an $L$, the group $\text{Gal}(L/K^{\boldsymbol{\gamma}})$
contains $p$ different
subgroups which map isomorphically onto $\text{Gal}(K/K^{\boldsymbol{\gamma}})$,
so there are $p$ elements $M\in{\mathcal Z}(\boldsymbol{\gamma})$ such that $MK=L$.
Therefore we have $|{\mathcal Z}(\boldsymbol{\gamma})|=p|{\mathcal W}(\boldsymbol{\gamma})|$,
where
\begin{equation}
{\mathcal W}(\boldsymbol{\gamma})=\{L:K\subset L\subset\Omega,\,[L:K]=p,\,
L/K^{\boldsymbol{\gamma}}\ \text{is Galois and nonabelian}\}.
\end{equation}
By class field theory, elements in ${\mathcal W}(\boldsymbol{\gamma})$ correspond
to subgroups $H\le V_K$ of index $p$ such that
$H$ is invariant under the action of $\langle \boldsymbol{\gamma}\rangle$ and such
that
$\langle \boldsymbol{\gamma}\rangle$ acts nontrivially on $V_K/H$. The number of such
subgroups $H$ is equal to the number of 1-dimensional subspaces
of $V_K$ with nontrivial action by $\langle\boldsymbol{\gamma}\rangle$.  Using
Corollary~\ref{VK}, we find that
\begin{equation} \label{W}
|{\mathcal W}(\boldsymbol{\gamma})|=
\begin{cases}
\displaystyle(d(\boldsymbol{\gamma})-1)\frac{p^{\frac{n_Fe_0f}{o(\boldsymbol{\gamma})}}-1}{p-1}
&\text{if}\ \zeta_p\notin K\ \text{or}\ \zeta_p\in K^{\boldsymbol{\gamma}},\\[.4cm]
\displaystyle(d(\boldsymbol{\gamma})-2)\frac{p^{\frac{n_Fe_0f}{o(\boldsymbol{\gamma})}}-1}{p-1}
+\frac{p^{\frac{n_Fe_0f}{o(\boldsymbol{\gamma})}+1}-1}{p-1}&\text{if}\ \zeta_p\in K\
\text{but}\ \zeta_p\notin K^{\boldsymbol{\gamma}}.\cr
\end{cases}
\end{equation}
Equation (\ref{n2}) now follows from (\ref{n2formula})--(\ref{W}).
\end{proof}

It follows from Lemmas \ref{n1lemma} and \ref{n2lemma} that
\begin{equation} \label{n1n2}
p\cdot m_1(\boldsymbol{\gamma})+m_2(\boldsymbol{\gamma})=
\begin{cases}
(p^2+p)p^{\frac{n_Fe_0f}{o(\boldsymbol{\gamma})}}-p^2&\text{if}\ \zeta_p\notin K^{\boldsymbol{\gamma}},
\\[.2cm]
2p^{\frac{n_Fe_0f}{o(\boldsymbol{\gamma})}+2}-p^2&\text{if}\ \zeta_p\in K^{\boldsymbol{\gamma}}.
\cr
\end{cases}
\end{equation}
In order to write down explicit formulas for $\sum_{L\in{\mathcal E}(K,1,p)}|\text{Aut}(L/F)|$,
we change notation slightly: We restore the subscript $h$
to $K$, and instead of $o(\boldsymbol{\gamma})$, $t(\boldsymbol{\gamma})$, $K^{\boldsymbol{\gamma}}$, we write
$o(h,i,j)$, $t(h,i,j)$, $K_h^{ij}$, where $\boldsymbol{\gamma}={\boldsymbol\mu}^i{\boldsymbol\nu}^j$.
If $\zeta_p\notin K_h$, then by combining (\ref{different}) with
(\ref{n1n2})
we get the following result.

%%%%%%%%%%%%%%%%%%%%%%%%%%%%%%%%
%   Proposition 5.3
%%%%%%%%%%%%%%%%%%%%%%%%%%%%%%%%
\begin{prop} \label{nozeta}
If $\zeta_p\notin K_h$ then
\begin{equation}
\begin{split}
\sum_{L\in{\mathcal E}(K_h,1,p)}|\text{\rm Aut}(L/F)|&=\sum_{i=0}^{b-1}
\sum_{j=0}^{\frac f{c_h}-1}\left((p^2+p)p^{\frac{n_Fe_0f}{o(h,i,j)}}-p^2\right) \\
&=-\frac{p^2bf}{c_h}+(p^2+p)\sum_{i=0}^{b-1}\sum_{j=0}^{\frac f{c_h}-1}
p^{\frac{n_Fe_0f}{o(h,i,j)}}, \\
\end{split}
\end{equation}
where
\begin{align}
o(h,i,j)&=\frac{f}{(f,c_hj)}\cdot\frac{b}{(b,t(h,i,j))}, \\[.2cm]
t(h,i,j)&=\frac{q^{\text{\rm lcm}(f,c_hj)}-1}{q^{c_hj}-1} i+\frac{uhc_hj}{(f,c_hj)}.
\end{align}
\end{prop}

In order to evaluate (\ref{sumsum}), we need to be able
to tell when $\zeta_p\in K_h$, and when $\zeta_p\in K_h^{ij}$.
Let $f_p$ be the residue degree and $e_p$ the ramification index of the extension
$F(\zeta_p)/F$. A necessary condition for $K_h$ to contain $\zeta_p$ is that $f_p\mid
f$ and $e_p\mid e_0$.  Therefore, in what follows, we will assume $f_p\mid f$ and $e_p\mid e_0$.
Then $E=F(\zeta_g,\pi_E)$ contains all the fields in
${\mathcal E}(F,f_p,e_p)$, including $F(\zeta_p)$.  Therefore
$F(\zeta_p)$
is the fixed field of a normal subgroup $H$ of $\text{Gal}(E/F)$.
Since the residue degree of $E/F(\zeta_p)$ is $f/f_p$, and the
ramification index is $e/e_p$, we easily see that
$H=\langle \boldsymbol{\sigma}^{e_p},\boldsymbol{\sigma}^l\boldsymbol{\tau}^{f_p}
\rangle$ for some $l$.
The following lemma shows that we can assume $l=1$.

%%%%%%%%%%%%%%%%%%%%%%%%%%%%%%
%  Lemma 5.4
%%%%%%%%%%%%%%%%%%%%%%%%%%%%%%
\begin{lem} \label{auto}
There is an automorphism $\Psi$ of $\text{\rm Gal}(E/F)=\langle\boldsymbol{\sigma},
\boldsymbol{\tau}\rangle$ such that

(i) $\Psi$ maps the inertia group $\langle\boldsymbol{\sigma}\rangle$ onto
itself.

(ii) $\Psi$ acts trivially on the quotient
$\text{\rm Gal}(E/F)/\langle\boldsymbol{\sigma}\rangle$.

(iii) $\Psi(H)=\langle\boldsymbol{\sigma}^{e_p},\boldsymbol{\sigma}
\boldsymbol{\tau}^{f_p}\rangle$.
\end{lem}

\begin{proof}
Since $\text{Gal}(E/F)/H\cong\text{Gal}(F(\zeta_p)/F)$ is cyclic,
it is generated by $\boldsymbol{\sigma}^a\boldsymbol{\tau}$ for some $a\in
{\mathbb Z}$. Define an automorphism
$\Psi_1$ of $\text{Gal}(E/F)$ by setting $\Psi_1(\boldsymbol{\sigma})=\boldsymbol{\sigma}$ and
$\Psi_1(\boldsymbol{\tau})=\boldsymbol{\sigma}^{-a}\boldsymbol{\tau}$.  Then
$\Psi_1(H)=\langle\boldsymbol{\sigma}^{e_p},\boldsymbol{\sigma}^d
\boldsymbol{\tau}^{f_p}\rangle$, where
\begin{equation}
d=l-\frac{q^{f_p}-1}{q-1} a.
\end{equation}
Furthermore, $\boldsymbol{\tau}=\Psi_1(\boldsymbol{\sigma}^a
\boldsymbol{\tau})$ generates the quotient
$\text{Gal}(E/F)/\Psi_1(H)$.  This implies that
$(e_p,d)=1$, so there is $k\in{\mathbb Z}$ such
that $dk\equiv1\pmod{e_p}$.  In addition, since the homomorphism
$({\mathbb Z}/e_0{\mathbb Z})^\times\rightarrow({\mathbb Z}/e_p{\mathbb Z})^\times$
is onto, we may choose $k$ so that $(e_0,k)=1$.
Define an automorphism $\Psi_2$ of
$\text{Gal}(E/F)$ by setting $\Psi_2(\boldsymbol{\sigma})=
\boldsymbol{\sigma}^k$ and
$\Psi_2(\boldsymbol{\tau})=\boldsymbol{\tau}$.  Then $\Psi=\Psi_2\circ\Psi_1$ satisfies the given conditions.
\end{proof}

By Lemma~\ref{auto} we have $H=\langle (\boldsymbol{\sigma}^s)^{e_p},\boldsymbol{\sigma}^s
(\boldsymbol{\sigma}^t\boldsymbol{\tau})^{f_p}\rangle$ for some $s,t\in{\mathbb Z}$ such that $(s,e_0)=1$.
Let $\tilde{\boldsymbol{\sigma}}=\boldsymbol{\sigma}^s$, $\tilde{\boldsymbol{\tau}}
=\boldsymbol{\sigma}^t\boldsymbol{\tau}$, $\tilde\pi_E=\zeta_{e_0(q-1)}^{-t}
\pi_E$, and $\tilde\pi_F=\zeta_{q-1}^{-t}\pi_F$. Then $\tilde\pi_F$
is a uniformizer of $F$ and $\tilde\pi_E^{e_0}=\tilde\pi_F$.
Furthermore,
\begin{equation}
\tilde{\boldsymbol{\sigma}}(\tilde\pi_E)=\zeta_e^s\tilde\pi_E,\
\tilde{\boldsymbol{\sigma}}(\zeta_g)=\zeta_g,\
\tilde{\boldsymbol{\tau}}(\tilde\pi_E)=\tilde\pi_E,\
\tilde{\boldsymbol{\tau}}(\zeta_g)=\zeta_g^q.
\end{equation}
It follows that by replacing $\pi_F$, $\pi_E$, $\boldsymbol{\sigma}$,
$\boldsymbol{\tau}$ with $\tilde\pi_F$, $\tilde\pi_E$,
$\tilde{\boldsymbol{\sigma}}$, $\tilde{\boldsymbol{\tau}}$
we may assume that
$H=\langle \boldsymbol{\sigma}^{e_p},\boldsymbol{\sigma}
\boldsymbol{\tau}^{f_p}\rangle$.  Under this assumption, for
every $x,y\in{\mathbb Z}$ the element
$\boldsymbol{\sigma}^x\boldsymbol{\tau}^y$ fixes $\zeta_p$ if and only if
$\boldsymbol{\sigma}^x\boldsymbol{\tau}^y\in\langle
\boldsymbol{\sigma}^{e_p},\boldsymbol{\sigma}\boldsymbol{\tau}^{f_p}\rangle$.
By (\ref{GalEF}) this is equivalent to $f_p\mid y$ and
\begin{equation}
x\equiv\frac{q^y-1}{q^{f_p}-1}\pmod{e_p}.
\end{equation}
Since $K_h$ is the subfield of $E$ fixed by
$\langle\boldsymbol{\sigma}^{-h}\boldsymbol{\tau}^f\rangle$, we see that $\zeta_p\in K_h$ if and only if
\begin{equation}
-h\equiv \frac{q^f-1}{q^{f_p}-1}\pmod {e_p}.
\end{equation}

     To determine whether $\zeta_p$ is in $K_h^{ij}$ we
recall that by Remark~\ref{future}, ${\boldsymbol\mu}$ is the restriction
of $\boldsymbol{\sigma}^{e_0/b}$ to $K_h$ and ${\boldsymbol\nu}$
is the restriction of
$\boldsymbol{\sigma}^{v_h}\boldsymbol{\tau}^{c_h}$ to $K_h$. It follows that $\boldsymbol{\gamma}={\boldsymbol\mu}^i{\boldsymbol\nu}^j$ is the restriction to $K_h$ of
\begin{equation}
(\boldsymbol{\sigma}^{e_0/b})^i(\boldsymbol{\sigma}^{v_h}
\boldsymbol{\tau}^{c_h})^j=
\boldsymbol{\sigma}^{r}\boldsymbol{\tau}^{c_hj},
\end{equation}
where
\begin{equation}
r=\frac{e_0i}{b}+\frac{q^{c_hj}-1}{q^{c_h}-1}v_h.
\end{equation}
Thus $\zeta_p\in K_h^{ij}$ if and only if $f_p\mid c_hj$ and
\begin{equation} \label{modep}
\frac{e_0i}{b}+\frac{q^{{c_h}j}-1}{q^{c_h}-1}v_h \equiv \frac{q^{{c_h}j}-1}{q^{f_p}-1} \pmod{e_p}.
\end{equation}

Using (\ref{different}), (\ref{n1n2}), and (\ref{modep}) we get a
formula for $\sum_{L\in{\mathcal E}(K_h,1,p)}
|\text{Aut}(L/F)|$ in the case $\zeta_p\in K_h$.

%%%%%%%%%%%%%%%%%%%%%%%%%%%%
%  Proposition 5.5
%%%%%%%%%%%%%%%%%%%%%%%%%%%%%
\begin{prop} \label{zeta}
If $\zeta_p\in K_h$ then
\begin{equation}
\begin{split}
\sum_{L\in{\mathcal E}(K_h,1,p)}|\text{\rm Aut}(L/F)|
=&-\frac{p^2bf}{c_h}+(p^2+p)\sum_{i=0}^{b-1}\sum_{j=0}^{\frac{f}{c_h}-1}
p^{\frac{n_Fe_0f}{o(h,i,j)}}\cr
&\hspace{.5cm}+(p^2-p)\sum_{\substack{0\le j<f/c_h\\[.05cm] f_p\mid c_hj}}
\,\sum_{i\in R_{h,j}}p^{\frac{n_Fe_0f}{o(h,i,j)}},
\end{split}
\end{equation}
where $R_{h,j}$ denotes the set of integers $0\le i<b$
satisfying the congruence (\ref{modep}).
\end{prop}

By combining Propositions \ref{nozeta} and \ref{zeta} with
equation (\ref{sumsum}), we get the main result of this
section.

%%%%%%%%%%%%%%%%%%%%%%%%%%%%%%%%%%%%
%  Theorem 5.6
%%%%%%%%%%%%%%%%%%%%%%%%%%%%%%%%%%%%
\begin{thm} \label{main1}
Let $F$ be a finite extension of ${\mathbb Q}_p$, let $f$ and $e$ be positive
integers such that $p\parallel e$, and set $e_0=e/p$.  Then the number of
$F$-isomorphism classes of extensions of $F$ with residue class degree $f$ and  ramification index $e$ is
\begin{equation}
{\mathfrak I}(F,f,e)=\frac{1}{fe}\sum_{h=0}^{e_0-1}\Biggl(
-\frac{p^2bf}{c_h}+\sum_{i=0}^{b-1}\sum_{j=0}^{\frac f{c_h}-1}(p^2+
\omega_{hij})p^{\frac{n_Fe_0f}{o(h,i,j)}}\Biggr),
\end{equation}
where \smallskip

(i) $q=$ the cardinality of the residue field of $F$,

(ii) $b=(e_0,q^f-1)$,

(iii) $u\in{\mathbb Z}$ satisfies $e_0u\equiv b\pmod{q^f-1}$,

(iv)  $c_h$ is the smallest positive integer such that $b\mid
(q^{c_h}-1)h$, \vspace{.1cm}

(v) $\displaystyle t(h,i,j)=
\frac{q^{\text{\rm lcm}(f,c_hj)}-1}{q^{c_hj}-1} i+
\frac{uhc_hj}{(f,c_hj)}$, \vspace{.1cm}

(vi) $\displaystyle o(h,i,j)=
\frac{f}{(f,c_hj)}\cdot
\frac{b}{(b,t(h,i,j))}$, \vspace{.2cm}

(vii) $\omega_{hij}=
\begin{cases}
p^2&\text{if}\ e(F(\zeta_p)/F)\mid e,\;f(F(\zeta_p)/F)\mid f,\;f(F(\zeta_p)/F)\mid c_hj,\ \\
&\hspace{.5cm}\text{and}\ h, i, j\ \text{satisfy (\ref{modep})},
\\
p&\text{otherwise}.
\end{cases}$
\end{thm}

%%%%%%%%%%%%%%%%%%%%%%%%%%%%%%%%%%%%%%%%%%%%%%%%%%
%  Section 6
%%%%%%%%%%%%%%%%%%%%%%%%%%%%%%%%%%%%%%%%%%%%%%%%%%
\section{The Case $p^2\parallel e$} \label{p2}

For the remainder of the paper we consider the case
$p^2\parallel e$. Let $f_1=f(F/{\mathbb Q}_p)$ and
$e_1=e(F/{\mathbb Q}_p)$, and write $e=p^2e_0$ with $p\nmid e_0$.
We make the following simplifying assumptions:
\begin{equation} \label{assumptions}
\begin{cases}
(p^{f_1f}-1,e_0)=1,\cr
f(F(\zeta_p)/F)\nmid f\ \text{or}\ e(F(\zeta_p)/F)>1.\cr
\end{cases}
\end{equation}
Some consequences of these assumptions are given in the
following proposition.

%%%%%%%%%%%%%%%%%%%%%%%%%%%%%%%%
%  Proposition 6.1
%%%%%%%%%%%%%%%%%%%%%%%%%%%%%%%%
\begin{prop} \label{consequences}
Assume the conditions in (\ref{assumptions}), and let
$K\in{\mathcal E}(F,f,e)$.  Then: \smallskip

(i) There is a unique field $L_K$ such that $F\subset L_K\subset K$ and
$K/L_K$ is totally ramified of degree $p^2$.

(ii) There is a unique field $E_K$ such that $F\subset E_K\subset K$ and
$E_K/F$ is totally ramified of degree $e_0$.

(iii) $E_K$ is a subfield of $L_K$.

(iv) $\text{\rm Aut}(K/F)$ acts trivially on $E_K$.

(v) $\zeta_p\notin K$.
\end{prop}

\begin{proof}
(i) -- (iv) follow from Proposition~\ref{unique2}. To prove
(v), note that ${e(F(\zeta_p)/F)\mid(p-1)}$
and $(p-1,e)=1$. Thus if $e(F(\zeta_p)/F)>1$ then
$e(F(\zeta_p)/F)\nmid e$.  Since $\zeta_p\in K$ implies
$e(F(\zeta_p))\mid e$ and $f(F(\zeta_p))\mid f$, we must
have $\zeta_p\notin K$.
\end{proof}

The approach we take here is somewhat different from that of
Sections \ref{p0} and \ref{p1}.  Even with the conditions in
(\ref{assumptions}), the computations that we will face are
quite lengthy.  To control the overall length of the paper
and maintain readability, we will describe the reasoning
behind our computations but omit the details.  All the
computations in these sections have been checked using
{\it Mathematica}.

For each positive integer $d$, let
\begin{equation}
{\mathcal B}_d=\{K\in{\mathcal E}(F,f,e):d\bigm||\text{Aut}(K/F)|\}.
\end{equation}
Using (\ref{known}) and the fact that $\sum_{d\mid n}\phi(d)=n$
we get
\begin{equation} \label{Mobius}
{\mathfrak I}(F,f,e)=\frac1{fe}\sum_{d>0}\,\phi(d)|{\mathcal B}_d|,
\end{equation}
where $\phi$ is the Euler function.
For $d>0$ and $i=0,1,2$, put
\begin{equation}
\begin{split}
{\mathcal C}_d^i=&\{K\in{\mathcal B}_d:\exists F\subset N\subset K
\text{ such that } K/N\ \text{is Galois}, \\
&\hspace{2cm}e(K/N)=p^i,\text{ and }d\mid[K:N]\}.
\end{split}
\end{equation}
By Proposition~\ref{consequences}(iv) we have
${\mathcal B}_d={\mathcal C}_d^2\cup{\mathcal C}_d^1\cup{\mathcal C}_d^0$.
Therefore
\begin{equation} \label{Bd}
|{\mathcal B}_d|=|{\mathcal C}_d^2|+|{\mathcal C}_d^1\setminus{\mathcal C}_d^2|+
|{\mathcal C}_d^0\setminus({\mathcal C}_d^1\cup{\mathcal C}_d^2)|.
\end{equation}
More precisely, ${\mathcal C}_d^2$ consists of the fields
$K\in{\mathcal B}_d$ such
that the ramification index of $K$ over the fixed field of $\text{Aut}(K/F)$ is $p^2$;
${\mathcal C}_d^1\setminus{\mathcal C}_d^2$ consists of the fields
$K\in{\mathcal B}_d$ such
that the ramification index of $K$ over the fixed field of $\text{Aut}(K/F)$ is $p$;
and ${\mathcal C}_d^0\setminus({\mathcal C}_d^1\cup{\mathcal C}_d^2)$ consists of the fields $K\in{\mathcal B}_d$ such
that the ramification index of $K$ over the fixed field of $\text{Aut}(K/F)$ is 1.

We will determine $|{\mathcal C}_d^2|$, $|{\mathcal C}_d^1\setminus{\mathcal C}_d^2|$, and
$|{\mathcal C}_d^0\setminus({\mathcal C}_d^1\cup{\mathcal C}_d^2)|$
separately in Sections {\ref{d2} -- \ref{d012}}. In our computations, we will frequently encounter a tower
of finite extensions ${\mathbb Q}_p\subset T\subset L\subset K$,
where $K/T$ is Galois, $L/T$ is unramified
of degree $d$, $K/L$ is abelian of degree $p^i$,
and $\zeta_p\notin L$.  We need to give an explicit description
of $\text{Gal}(K/T)$ in terms of $\text{Gal}(K/L)$.

Let $L/T$ be an unramified extension of degree $d$ and define
\begin{equation}
{\mathfrak K}(T,L;p^i)=\{K:L\subset K\subset\Omega,\;
K/L\text{ abelian of degree }p^i,\;
K/T\text{ Galois}\}.
\end{equation}
For positive integers $m$ and $n$, put
\begin{align}
C(m)&=\left[
\begin{matrix}
0&1\\
&0\\
&&\ddots\\
&&&0&1\\
1&&&&0
\end{matrix}
\right]_{m\times m} \\[.1cm]
D(m,n)&=\left[\left.
\begin{matrix}
C(m)\\
&\ddots\\
&&C(m)\\
\end{matrix}
\right]\right\}{\scriptstyle n} \\[.1cm]
E(m,n)&=\left[
\begin{matrix}
1\cr
&D(m,n)\cr
\end{matrix}\right].
\end{align}
We have an isomorphism of groups
\begin{equation}
L^\times/(L^\times)^{p^i}\cong\frac{\langle\pi_L\rangle}{\langle\pi_L^{p^i}\rangle}
\times\frac{1+\pi_L{\mathcal O}_L}{(1+\pi_L{\mathcal O}_L)^{p^i}}\cong
({\mathbb Z}/p^i{\mathbb Z})^{1+n_L}.
\end{equation}
Let $\boldsymbol{\sigma}$ be the Frobenius map of $L/T$. Then by
Proposition~\ref{U1}(i), there is a
$({\mathbb Z}/p^i{\mathbb Z})$-basis $\mathcal S$ for
$L^\times/(L^\times)^{p^i}\cong
({\mathbb Z}/p^i{\mathbb Z})^{1+n_L}$ such that
the matrix of $\boldsymbol{\sigma}$ with respect to $\mathcal S$ is
$E(d,n_T)$.

Define ${\mathcal H}(m,n;p^i)$ to be the set of all
$E(m,n)$-invariant subgroups of
$({\mathbb Z}/p^i{\mathbb Z})^{1+mn}$ of index $p^i$.
By class field theory, there is a bijection between the set of all
$\boldsymbol{\sigma}$-invariant subgroups of $L^\times/(L^\times)^{p^i}$ of
index $p^i$ and ${\mathfrak K}(T,L;p^i)$. This bijection induces a bijection
between ${\mathcal H}(d,n_T;p^i)$ and ${\mathfrak K}(T,L;p^i)$, which
is denoted by $H\mapsto K_H$. Furthermore, $\text{Gal}(K_H/L)\cong
({\mathbb Z}/p^i{\mathbb Z})^{1+n_L}/H$ for each $H\in
{\mathcal H}(d,n_T;p^i)$.  For each
$u\in({\mathbb Z}/p^i{\mathbb Z})^{1+n_L}/H$ let
$\boldsymbol{\omega}(u)$ be the element of $\text{Gal}(K_H/L)$
which corresponds to $u$ under this isomorphism.

Let $M/T$ be the maximal unramified subextension of $K_H/T$.
Since $\boldsymbol{\sigma}$
is the Frobenius of $L/T$, it can be extended to the Frobenius
$\boldsymbol{\sigma}'$ of $M/T$.  Let $\boldsymbol{\theta}$ be
an arbitrary extension of $\boldsymbol{\sigma}'$ to an element
of $\text{Gal}(K_H/T)$.  Then $\text{Gal}(K_H/T)$ is generated
by $\text{Gal}(K_H/L)\cong({\mathbb Z}/p^i{\mathbb Z})^{1+n_L}/H$
and $\boldsymbol{\theta}$, where
\begin{equation}
\boldsymbol{\theta} \boldsymbol{\omega}(u)\boldsymbol{\theta}^{-1}
=\boldsymbol{\omega}(u^{\boldsymbol{\sigma}})=
\boldsymbol{\omega}(E(d,n_T)u)
\end{equation}
for all $u\in({\mathbb Z}/p^i{\mathbb Z})^{1+n_L}/H$, and
$\boldsymbol{\theta}^d=\boldsymbol{\omega}
\left(\genfrac{[}{]}{0pt}{1}a\alpha\right)$
for some $a\in{\mathbb Z}/p^i{\mathbb Z}\cong
\langle\pi_L\rangle/\langle\pi_L^{p^i}\rangle$ and
$\alpha\in({\mathbb Z}/p^i{\mathbb Z})^{n_L}\cong
(1+\pi_L{\mathcal O}_L)/(1+\pi_L{\mathcal O}_L)^{p^i}$.  It
follows from the definition of $\boldsymbol{\theta}$ that the
restriction of $\boldsymbol{\theta}^d$ to $M/L$ is the Frobenius.
By class field theory this implies $a\equiv1\pmod{p^i}$.

To summarize, we have the following description of the structure of
$\text{Gal}(K_H/T)$.

%%%%%%%%%%%%%%%%%%%%%%%%%%%%%%%%
%  Proposition 6.2
%%%%%%%%%%%%%%%%%%%%%%%%%%%%%%%%
\begin{prop} \label{GalKHT}
Let ${\mathbb Q}_p\subset T\subset L$ be finite extensions such
that $L/T$ is unramified of degree $d$ and $\zeta_p\notin L$.
Let $H\mapsto K_H$ be the bijection between
${\mathcal H}(d,n_T;p^i)$ and ${\mathfrak K}(T,L;p^i)$ induced
by class field theory.  Then
for each $H\in{\mathcal H}(d,n_T;p^i)$ we have an isomorphism
$\boldsymbol{\omega}:({\mathbb Z}/p^i{\mathbb Z})^{1+n_L}/H
\rightarrow\text{Gal}(K_H/L)$ such that $\text{Gal}(K_H/T)$ is
generated by $\text{Gal}(K_H/L)$ and an element
$\boldsymbol{\theta}$ satisfying
\begin{equation} \label{relations}
\begin{cases}
\boldsymbol{\theta}^d=\boldsymbol{\omega}
\left(\genfrac{[}{]}{0pt}{1}1\alpha\right), \\
\boldsymbol{\theta}\boldsymbol{\omega}(u)\boldsymbol{\theta}^{-1}=
\boldsymbol{\omega}(E(d,n_T)u)\ \text{for all}\
u\in({\mathbb Z}/p^i{\mathbb Z})^{1+n_L}/H,
\end{cases}
\end{equation}
where $\genfrac{[}{]}{0pt}{1}1\alpha\in({\mathbb Z}/p^i{\mathbb Z})^{n_L+1}/H$
depends only on $T$ and $L$.
\end{prop}

Next, we list the elements $H$ of
${\mathcal H}(d,n_T;p^i)$ for $i=1,2$.
For each such $H$ we will give a more explicit description
of the  structure of $\text{Gal}(K_H/T)$ than
that given in Proposition~\ref{GalKHT}.  These explicit descriptions
will be essential in our later calculations.

%%%%%%%%%%%%%%%%%%%%%%%%%%%%%%%%%
%   Corollary 6.3
%%%%%%%%%%%%%%%%%%%%%%%%%%%%%%%%%
\begin{cor} \label{explicit}
Let ${\mathbb Q}_p\subset T\subset L\subset\Omega$ be as in Proposition~\ref{GalKHT}.
Then the elements of ${\mathcal H}(d,n_T;p)$ are the groups of the form
\begin{equation}
H(\lambda,a)=a^\bot=\{x\in({\mathbb Z}/p{\mathbb Z})^{1+n_L}:x^ta=0\},
\end{equation}
where $0\ne a\in({\mathbb Z}/p{\mathbb Z})^{1+n_L}$ satisfies
$E(d,n_T)^t a=\lambda a$ for some
$\lambda\in{\mathbb Z}/p{\mathbb Z}$ such that $\lambda^d=1$.
Furthermore, $\text{\rm Gal}(K_{H(\lambda,a)}/T)$ is generated by
\begin{equation}
\text{\rm Gal}(K_{H(\lambda,a)}/L)=
\langle\boldsymbol{\kappa}\rangle\cong{\mathbb Z}/p{\mathbb Z}
\end{equation}
and an element $\boldsymbol{\theta}$ such that
$\boldsymbol{\theta}^d=\boldsymbol{\kappa}^{c(a)}$ and
$\boldsymbol{\theta}\boldsymbol{\kappa}\boldsymbol{\theta}^{-1}
=\boldsymbol{\kappa}^{\lambda}$,
where $c(a)=\genfrac{[}{]}{0pt}{1}1\alpha^t a\in{\mathbb Z}/p{\mathbb Z}$
for some fixed $\alpha\in({\mathbb Z}/p{\mathbb Z})^{n_L}$.
\end{cor}

We omit the proof of Corollary~\ref{explicit} since it is a simpler version of the proof of the
next corollary.

%%%%%%%%%%%%%%%%%%%%%%%%%%%%%%%%%
%   Corollary 6.4
%%%%%%%%%%%%%%%%%%%%%%%%%%%%%%%%%
\begin{cor} \label{long}
Let ${\mathbb Q}_p\subset T\subset L\subset\Omega$ be as in
Proposition~\ref{GalKHT}, and let $H$ be a subgroup of
$({\mathbb Z}/p^2{\mathbb Z})^{1+n_L}$ of index $p^2$.
\smallskip

(i) If $({\mathbb Z}/p^2{\mathbb Z})^{1+n_L}/H\cong
{\mathbb Z}/p^2{\mathbb Z}$, then $H\in {\mathcal H}(d,n_T;p^2)$
if and only if $H$ is of the form
\begin{equation}
H(\lambda, a)=a^\bot=
\{x\in({\mathbb Z}/p^2{\mathbb Z})^{1+n_L}:x^ta=0\},
\end{equation}
where $\lambda\in{\mathbb Z}/p^2{\mathbb Z}$ satisfies
$\lambda^d=1$ and $a\in({\mathbb Z}/p^2{\mathbb Z})^{1+
n_L}\setminus (p{\mathbb Z}/p^2{\mathbb Z})^{1+
n_L}$ satisfies $E(d,n_T)^ta=\lambda a$.
Furthermore, $\text{\rm Gal}(K_{H(\lambda,a)}/T)$
is generated by
\begin{eqnarray}
\text{\rm Gal}(K_{H(\lambda,a)}/L)
=\langle\boldsymbol{\kappa}\rangle
\cong{\mathbb Z}/p^2{\mathbb Z}
\end{eqnarray}
and an element $\boldsymbol{\theta}$ such that
$\boldsymbol{\theta}^d=\boldsymbol{\kappa}^{c(a)}$ and
$\boldsymbol{\theta}\boldsymbol{\kappa}\boldsymbol{\theta}^{-1}
=\boldsymbol{\kappa}^{\lambda}$,
where ${c(a)=\genfrac{[}{]}{0pt}{1}1\alpha^t a\in{\mathbb Z}/p^2{\mathbb Z}}$
for some fixed $\alpha\in({\mathbb Z}/p^2{\mathbb Z})^{n_L}$.
\smallskip

(ii) If $({\mathbb Z}/p^2{\mathbb Z})^{1+
n_L}/H\cong ({\mathbb Z}/p{\mathbb Z})^2$, then $H\supset (p{\mathbb Z}/p^2{\mathbb Z})^{1+
n_L}$.
In this case $H\in{\mathcal H}(d,n_T;p^2)$
if and only if $(H/(p{\mathbb Z}/p^2{\mathbb Z})^{1+
n_L})^\bot$ is the column space of some
${(1+n_L)\times 2}$ matrix $A$ over
${\mathbb Z}/p{\mathbb Z}$ such that $\text{\rm rank}(A)=2$ and $E(d,n_T)^tA=A\Lambda$ for some
$\Lambda\in \text{\rm GL}(2,p)$ satisfying $\Lambda^d=I_2$. Equivalently,
$H\in{\mathcal H}(d,n_T;p^2)$ if and only of $H$ is of the form
\begin{equation}
H(\Lambda,A)=\{x\in({\mathbb Z}/p^2{\mathbb Z})^{1+n_L}:x^tA\equiv 0\pmod p\},
\end{equation}
where $A$ and $\Lambda$ are as above. Furthermore, $\text{\rm Gal}(K_{H(\Lambda,A)}/T)$ is generated
by
\begin{equation}
\text{\rm Gal}(K_{H(\Lambda,A)}/L)=
\langle\boldsymbol{\kappa}_1,\boldsymbol{\kappa}_2\rangle
\cong({\mathbb Z}/p{\mathbb Z})^2
\end{equation}
and an element
$\boldsymbol{\theta}$ such that
$\boldsymbol{\theta}^d=\boldsymbol{\kappa}_1^{c_1(A)}\boldsymbol{\kappa}_2^{c_2(A)}$,
$\boldsymbol{\theta}\boldsymbol{\kappa}_1\boldsymbol{\theta}^{-1}
=\boldsymbol{\kappa}_1^{\lambda_{11}}\boldsymbol{\kappa}_2^{\lambda_{12}}$,
and $\boldsymbol{\theta}\boldsymbol{\kappa}_2\boldsymbol{\theta}^{-1}
=\boldsymbol{\kappa}_1^{\lambda_{21}}\boldsymbol{\kappa}_2^{\lambda_{22}}$,
where $(c_1(A),c_2(A))=\genfrac{[}{]}{0pt}{1}1\alpha^t A\in({\mathbb Z}/p{\mathbb Z})^2$
for some fixed $\alpha\in({\mathbb Z}/p{\mathbb Z})^{n_L}$
and $[\lambda_{ij}]=\Lambda$.
\end{cor}

\begin{proof}
In both cases the necessary and sufficient conditions for $H$ to
be an element of ${\mathcal H}(d,n_T;p^2)$ are straightforward
from the definitions.  It remains to show that for
$H\in {\mathcal H}(d,n_T;p^2)$, the structure of
$\text{Gal}(K_H/T)$ is as
described.  We will only give the argument for case (i), as
case (ii) is quite similar.

By Proposition~\ref{GalKHT}, $\text{Gal}(K_H/T)$ is generated
by $\text{Gal}(K_H/L)\cong({\mathbb Z}/p^2{\mathbb Z})^{1+n_L}/H$
and an element $\boldsymbol{\theta}$ satisfying the relations
(\ref{relations}).
We have $H=H(\lambda,a)=a^{\bot}$ for some $\lambda\in{\mathbb
Z}/p^2{\mathbb Z}$ and
$a\in({\mathbb Z}/p^2{\mathbb Z})^{n_L}\setminus
(p{\mathbb Z}/p^2{\mathbb Z})^{n_L}$
such that $E(d,n_T)^ta=\lambda a$ and
$\lambda^d=1$.  There is a canonical isomorphism
\begin{equation}
\begin{matrix}
\psi:&({\mathbb Z}/p^2{\mathbb Z})^{n_L+1}/a^\bot&\longrightarrow&
\text{Hom}_{{\mathbb Z}/p^2{\mathbb Z}}(\langle a\rangle,{\mathbb Z}/p^2{\mathbb Z}),\cr
&x+a^\bot&\longmapsto&\langle\,\cdot\,,x\rangle
\end{matrix}
\end{equation}
where $\langle\,\cdot\,,\cdot\,\rangle$ is the standard inner product
on $({\mathbb Z}/p^2{\mathbb Z})^{n_L}$.  It follows that the
conjugation action of $\boldsymbol{\theta}$ on
$\text{Gal}(K_H/L)$ induces an action of $\boldsymbol{\theta}$ on
$\text{Hom}_{{\mathbb Z}/p^2{\mathbb Z}}(\langle a\rangle,{\mathbb Z}/p^2{\mathbb
Z})$.
Let $\phi$ be the unique element of
$\text{Hom}(\langle a\rangle,{\mathbb Z}/p^2{\mathbb Z})$
such that $\phi(a)=1$. Then
$\psi\bigl(\genfrac{[}{]}{0pt}{1}1\alpha\bigr)
=\bigl(\genfrac{[}{]}{0pt}{1}1\alpha^t a\bigr)\phi$ and
$\boldsymbol{\theta}\cdot\phi=\lambda\phi$ for all $u\in{\mathbb Z}/p^2{\mathbb Z}$.
Therefore, by identifying $\text{Gal}(K_H/L)$ with
$\text{Hom}_{{\mathbb Z}/p^2{\mathbb Z}}(\langle a\rangle,{\mathbb Z}/p^2{\mathbb Z})$ using $\psi$ and
identifying $\text{Hom}_{{\mathbb Z}/p^2{\mathbb Z}}(\langle a\rangle,{\mathbb Z}/p^2{\mathbb Z})$ with
${\mathbb Z}/p^2{\mathbb Z}$ using the basis $\{\phi\}$, we see that $\text{Gal}(K_H/T)$ is generated by
$\text{Gal}(K_H/L)=\langle\boldsymbol{\kappa}\rangle\cong{\mathbb Z}/p^2{\mathbb Z}$ and an element $\boldsymbol{\theta}$ satisfying the relations
$\boldsymbol{\theta}^d=\boldsymbol{\kappa}^{c(a)}$ and
$\boldsymbol{\theta}\boldsymbol{\kappa}\boldsymbol{\theta}^{-1}=
\boldsymbol{\kappa}^{\lambda}$
with $c(a)=\genfrac{[}{]}{0pt}{1}1\alpha^t a$, as claimed.
\end{proof}

%%%%%%%%%%%%%%%%%%%%%%%%%%%%%%%%%%%%%%%%%%%%%%%%%%
%    Section 7
%%%%%%%%%%%%%%%%%%%%%%%%%%%%%%%%%%%%%%%%%%%%%%%%%%
\section{Determination of $|{\mathcal C}_d^2|$} \label{d2}

The goal of this section is to determine $|{\mathcal C}_d^2|$.
We retain the notation of Section~\ref{p2}.
In addition, we set $n=n_F=f_1e_1$. Observe that
$|{\mathcal C}_d^2|=0$ if $d\nmid p^2f$. Also note that
${\mathcal C}_d^2={\mathcal C}_{d'}^2$ where $d'=\text{lcm}(p^2,d)$.
Thus we assume that $d\mid p^2f$ and $p^2\mid d$.

Let ${\mathcal X}$ be the set of all $(T,L,K)$ in the diagram
\setlength{\unitlength}{3mm}
\begin{equation*}
\begin{picture}(2,17)
\put(1,0){\makebox(0,0){$F$}}
\put(1,5){\makebox(0,0){$T$}}
\put(1,9){\makebox(0,0){$L$}}
\put(1,13){\makebox(0,0){$K$}}
\put(1,17){\makebox(0,0){$\Omega$}}
\put(1.5,10.8){$\scriptstyle p^2$}
\put(-1.8,10.8){$\scriptstyle\text{(ram)}$}
\put(1.5,6.8){$\scriptstyle \frac d{p^2}$}
\put(-1.2,6.8){$\scriptstyle \text{(un)}$}
\put(1.5,3){$\scriptstyle f(T/F)=\frac{p^2f}d$}
\put(1.5,1.5){$\scriptstyle e(T/F)=e_0$}
\put(1,1){\line(0,1){3}}
\put(1,6){\line(0,1){2}}
\put(1,10){\line(0,1){2}}
\put(1,14){\line(0,1){2}}
\end{picture}
\end{equation*}
such that $K/T$ is Galois. Using Proposition~\ref{consequences} (i) -- (iv),
we see that
$(T,L,K)\mapsto K$ gives a bijection between ${\mathcal X}$ and
${\mathcal C}_d^2$. Hence
$|{\mathcal C}_d^2|=|{\mathcal X}|$.
Meanwhile, $|{\mathcal X}|$ can be computed by counting the elements $(T,L,K)\in
{\mathcal X}$ in the order $T, L, K$.

Fix $T\in{\mathcal E}\bigl(F,\frac{p^2f}d,e_0\bigr)$, let $L/T$ be
unramified of degree
$d/p^2$, and let $M/L$ be unramified of degree $p$.  Then we have
\begin{equation} \label{difference}
|\{K:(T,L,K)\in{\mathcal X}\}|=|{\mathfrak K}(T,L;p^2)|-|{\mathfrak K}(T,M;p)|.
\end{equation}
Using Corollary~\ref{explicit}, we find that
\begin{equation} \label{kTMp}
\begin{split}
|{\mathfrak K}(T,M;p)|&=\biggl|{\mathcal H}\!\biggl(\frac dp, \frac{p^2e_0fn}{d};p\biggr)\biggr|\cr
&=(p-1,d)\frac{p^{\frac 1d p^2e_0fn}-1}{p-1}+p^{\frac 1d p^2e_0fn}.\cr
\end{split}
\end{equation}
%
%%%%%%%%%%%%%%%%%%%%%%%%%%%%%%
%  Lemma 7.1
%%%%%%%%%%%%%%%%%%%%%%%%%%%%%%
\begin{lem} \label{kTLp2}
We have
\begin{equation} \label{sizeKTLp2}
|{\mathfrak K}(T,L;p^2)|=
\begin{cases}
\left.\begin{split}
&\frac 12 (p-1,d)^2\frac{(p^{\frac 1d p^2e_0fn}-1)^2}{(p-1)^2}\cr
&\;\;+(p-1,d)\frac{p^{\frac 1d p^2e_0fn}(p^{\frac 1d p^2e_0fn}-1)}{p-1}
\cr
&\;\;+p^{\frac 1d 2p^2e_0fn}+\frac 12(p^2-1,d)\frac{p^{\frac 1d 2p^2e_0fn}-1}{p^2-1}
\end{split}
\right\}&\text{ if } p^3\nmid d,\cr
\cr
\left.\begin{split}
&\frac 12 (p-1,d)^2\frac{(p^{\frac 1d p^2e_0fn}-1)^2}{(p-1)^2}\cr
&\;\;+2(p-1,d)\frac{p^{\frac 1d p^2e_0fn}(p^{\frac 1d p^2e_0fn}-1)}{p-1}
\cr
&\;\;+p^{\frac 1d p^2e_0fn}(2p^{\frac 1d p^2e_0fn}-1)\cr
&\;\;+\frac 12(p^2-1,d)\frac{p^{\frac 1d 2p^2e_0fn}-1}{p^2-1}\cr
\end{split}
\right\}&\text{ if }p^3\mid d.\cr
\end{cases}
\end{equation}
\end{lem}

\begin{proof}
First, we have
\begin{equation} \label{H1H2}
|{\mathfrak K}(T,L;p^2)|
=\biggl|{\mathcal H}\!\biggl(\frac d{p^2},\frac{p^2e_0fn}{d};p^2\biggr)\biggr|
=|{\mathcal H}_1|+|{\mathcal H}_2|,
\end{equation}
where ${\mathcal H}_1$ and ${\mathcal H}_2$ are the subsets
of ${\mathcal H}\bigl(\frac d{p^2},\frac 1d p^2e_0fn;p^2\bigr)$
corresponding to the two cases of Corollary~\ref{long}. Put
$E=E\bigl(\frac d{p^2},\frac 1d p^2e_0fn\bigr)$.  Then
\begin{align}
{\mathcal H}_1=&\{H(\lambda,a):a\in({\mathbb Z}/p^2{\mathbb Z})^{1+e_0fn}\setminus
(p{\mathbb Z}/p^2{\mathbb Z})^{1+e_0fn}, \\
\nonumber&\hspace{2cm}\lambda\in{\mathbb Z}/p^2{\mathbb Z},
\ E^ta=\lambda a,\ \lambda^{d/p^2}=1\},
\\[.1cm]
{\mathcal H}_2=&\{H(\Lambda,A):A\in M_{(1+e_0fn)\times 2}({\mathbb Z}/p{\mathbb Z}),\ \text{rank}(A)=2, \\
\nonumber&\hspace{2cm}\Lambda\in \text{GL}(2,p),\ E^tA=A\Lambda,\ \Lambda^{d/p^2}=I_2\}.
\end{align}

Let $\lambda\in {\mathbb Z}/p^2{\mathbb Z}$ satisfy $\lambda^{d/p^2}=1$.  Then
\begin{equation}
|\{H(\lambda,a)\in{\mathcal H}_1\}|=\frac{|{\mathcal A}_\lambda|}{p^2-p}
\end{equation}
where
\begin{equation}
{\mathcal A}_\lambda=\{a\in ({\mathbb Z}/p^2{\mathbb Z})^{1+e_0fn}\setminus
(p{\mathbb Z}/p^2{\mathbb Z})^{1+e_0fn}: E^ta=\lambda a\}.
\end{equation}
Hence
\begin{equation} \label{H1}
|{\mathcal H}_1|=\sum_{
\substack{
\lambda\in{\mathbb Z}/p{\mathbb Z}\\ \lambda^{d/p^2}=1}}
|\{H(\lambda,a)\in{\mathcal H}_1\}|
=\frac 1{p^2-p}\sum_{\substack{\lambda\in{\mathbb Z}/p{\mathbb Z}\\ \lambda^{d/p^2}=1}}|{\mathcal A}_\lambda|.
\end{equation}
The cardinality of ${\mathcal A}_\lambda$ can be easily determined:
\begin{equation} \label{Alambda}
|{\mathcal A}_\lambda|=
\begin{cases}
p^{2(1+\frac 1d p^2e_0fn)}-p^{1+\frac 1d p^2e_0fn}&\text{if}\ \lambda=1,\cr
p^{1+\frac 2d p^2e_0fn}-p^{1+\frac 1d p^2e_0fn}&\text{if}
\ \lambda\equiv 1\pmod p\text{ but }\lambda \ne 1,\cr
p^{\frac 2d p^2e_0fn}-p^{\frac 1d p^2e_0fn}&\text{if}
\ \lambda\not\equiv 1\pmod p.\cr
\end{cases}
\end{equation}
Combining (\ref{H1}) and (\ref{Alambda}), we find that
\begin{equation} \label{H1formula}
\begin{split}
|{\mathcal H}_1|
=\,&(p-1,d)\biggl(p,\frac{d}{p^2}\biggr)\frac{p^{\frac 1d p^2e_0fn-1}(p^{\frac 1d p^2e_0fn}-1)}{p-1}\cr
&\;\;+\biggl(p,\frac{d}{p^2}\biggr)p^{\frac 1d p^2e_0fn-1}(p^{\frac 1d p^2e_0fn}-1)+
p^{\frac 1d 2p^2e_0fn}.\cr
\end{split}
\end{equation}

To compute $|{\mathcal H}_2|$, we first observe that for $H(\Lambda, A)\in{\mathcal H}_2$ and $Q\in
\text{GL}(2,p)$, we have
\begin{equation} \label{HPX}
H(\Lambda,A)=H(Q^{-1}\Lambda Q,AQ).
\end{equation}
Also note that for $\Lambda\in \text{GL}(2,p)$ with
$\Lambda^{d/p^2}=I$, we have $H(\Lambda, A_1)=H(\Lambda,A_2)$
if and only if $A_1=A_2Q$ for some $Q$ in the centralizer $\text{cent}(\Lambda)$ of $\Lambda$ in
$\text{GL}(2,p)$. Thus for each $\Lambda\in GL(2,p)$,
\begin{equation} \label{cardHPX}
|\{H(\Lambda,A)\in{\mathcal H}_2\}|=\frac{|{\mathcal A}_{\Lambda}|}{|\text{cent}(\Lambda)|},
\end{equation}
where
\begin{equation}
{\mathcal A}_{\Lambda}=\{A\in M_{(1+e_0fn)\times 2}({\mathbb Z}/p{\mathbb Z}):\text{rank}(A)=2,\;
E^tA=A\Lambda\}.
\end{equation}
By (\ref{HPX}) and (\ref{cardHPX}), we have
\begin{equation} \label{H2}
|{\mathcal H}_2|=\sum_{\Lambda}\;\frac{|{\mathcal A}_{\Lambda}|}{|\text{cent}(\Lambda)|},
\end{equation}
where $\Lambda$ runs over the set of canonical forms in $M_{2\times 2}(
{\mathbb Z}/p{\mathbb Z})$ with $\Lambda^{d/p^2}=I$.  We can
compute $|\text{cent}(\Lambda)|$ using
the formula in \cite{Hou4}, and $|{\mathcal A}_{\Lambda}|$
using the well-known formula for the dimension of the solution
set of $E^tA=A\Lambda$ (\cite{Horn}, Theorem~4.4.14). We omit
the details of these
computations and record the result for $|{\mathcal H}_2|$ below.
\begin{equation} \label{H2formula}
|{\mathcal H}_2|=
\begin{cases}
\left.\begin{split}
&\frac 12(p-1,d)^2\frac {(p^{\frac 1d p^2e_0fn}-1)^2}{(p-1)^2}
\cr
&\;\;+(p-1,d)p^{-1+\frac 1d p^2e_0fn}
(p^{\frac 1d p^2e_0fn}-1)\cr
&\;\;-(p^{\frac 1d p^2e_0fn}-1)p^{-1+\frac 1d p^2e_0fn}
\cr
&\;\;+\frac 12(p^2-1,d)\frac {p^{\frac 1d 2p^2e_0fn}-1}{p^2-1}
\end{split}
\right\}&\text{ if }p^3\nmid d\cr
\cr
\left.\begin{split}
&\frac 12(p-1,d)^2\frac {(p^{\frac 1d p^2e_0fn}-1)^2}{(p-1)^2}
\cr
&\;\;+(p-1,d)\frac{p^{\frac 1d p^2e_0fn}(p^{\frac 1d
p^2e_0fn}-1)}{p-1}
\cr
&\;\;+\frac 12(p^2-1,d)\frac {p^{\frac 1d 2p^2e_0fn}-1}{p^2-1}
\end{split}
\right\}&\text{ if }p^3\mid d.\cr
\end{cases}
\end{equation}
Finally, (\ref{sizeKTLp2}) follows from (\ref{H1H2}), (\ref{H1formula}), and (\ref{H2formula}).
\end{proof}

We now state the main result of this section.

%%%%%%%%%%%%%%%%%%%%%%%%%%%%
%  Proposition 7.2
%%%%%%%%%%%%%%%%%%%%%%%%%%%%
\begin{prop} \label{Cd2}
Let $d'=\text{\rm lcm}(p^2,d)$. Then we have
\begin{equation*}
\begin{split}
&|{\mathcal C}_d^2|=
\begin{cases}
0&\text{if}\ d\nmid p^2f,\cr
\cr
\left.\begin{split}
&e_0\biggl[\frac 12(p-1,d)^2\frac{(p^{\frac 1{d'}
p^2e_0fn}-1)^2}{(p-1)^2}
\cr
&\;\;\;+(p-1,d)\frac{(p^{\frac 1{d'} p^2e_0fn}-1)^2}{p-1}
\cr
&\;\;\;+p^{\frac 1{d'} p^2e_0fn}(p^{\frac 1{d'} p^2e_0fn}-1)
\cr
&\;\;\;+\frac 12(p^2-1,d)
\frac 1{p^2-1}(p^{\frac 1{d'} 2p^2e_0fn}-1)\biggr]\cr
\end{split}
\right\}&\text{if }d\mid p^2f\text{ and }p^3\nmid d,\cr
\cr
\left.\begin{split}
&e_0\biggl[\frac 12(p-1,d)^2\frac{(p^{\frac 1{d'}
p^2e_0fn}-1)^2}{(p-1)^2}
\cr
&\;\;\;+(p-1,d)\frac{(p^{\frac 1{d'} p^2e_0fn}-1)(2p^{\frac 1{d'} p^2e_0fn}
-1)}{p-1}\cr
&\;\;\;+2p^{\frac 1{d'} p^2e_0fn}(p^{\frac 1{d'} p^2e_0fn}-1)
\cr
&\;\;\;+\frac 12(p^2-1,d)\frac{p^{\frac 1{d'} 2p^2e_0fn}-1}{p^2-1}\biggr]\cr
\end{split}
\right\}&\text{if }d\mid p^2f\text{ and }p^3\mid d.\cr
\end{cases}\cr
\end{split}
\end{equation*}
\end{prop}

\begin{proof}
By the comments at the beginning of the section, it
suffices to prove the proposition under the assumptions
$d\mid p^2f$ and $p^2\mid d$.  Using (\ref{difference}) we
get
\begin{equation}
|{\mathcal C}_d^2|
=|{\mathcal X}|
=\sum_{T\in{\mathcal E}\left(F, \frac{1}{d}p^2f,e_0\right)}\bigl(|{\mathfrak K}(T,L;p^2)|-
|{\mathfrak K}(T,M;p)|\bigr).
\end{equation}
The proposition now follows from (\ref{kTMp}), (\ref{sizeKTLp2}), and
the fact that
$\left|{\mathcal E}\left(F,\frac1d p^2f,e_0\right)\right|=e_0$.
\end{proof}

%%%%%%%%%%%%%%%%%%%%%%%%%%%%%%%%%%%%%%%%%%%%%%%%%%
%    Section 8
%%%%%%%%%%%%%%%%%%%%%%%%%%%%%%%%%%%%%%%%%%%%%%%%%%
\section{Determination of $|{\mathcal C}_d^1\setminus{\mathcal C}_d^2|$}
\label{d12}

In this section we compute the cardinality of
${\mathcal C}_d^1\setminus{\mathcal C}_d^2$.  Recall that
${\mathcal C}_d^1\setminus{\mathcal C}_d^2$ consists of the
fields $K\in{\mathcal E}(F,f,e)$ such that
$d\bigm||\text{Aut}(K/F)|$ and the ramification index of $K$
over the fixed field $N_K$ of $\text{Aut}(K/F)$ is $p$.
Thus $|{\mathcal C}_d^1\setminus{\mathcal C}_d^2|=0$ if
$d\nmid pf$.  Furthermore, setting $d'=\text{lcm}(p,d)$, we
have ${\mathcal C}_d^1\setminus{\mathcal C}_d^2=
{\mathcal C}_{d'}^1\setminus{\mathcal C}_{d'}^2$.
Therefore we may assume that $d\mid pf$ and $p\mid d$.

Let ${\mathcal X}$ be the set of all $(M,E,K)$ in the diagram
\setlength{\unitlength}{3mm}
\begin{equation*}
\begin{picture}(2,17)
\put(1,0){\makebox(0,0){$F$}}
\put(1,5){\makebox(0,0){$M$}}
\put(1,9){\makebox(0,0){$E$}}
\put(1,13){\makebox(0,0){$K$}}
\put(1,17){\makebox(0,0){$\Omega$}}
\put(1.5,10.8){$\scriptstyle p$}
\put(-1.8,10.8){$\scriptstyle\text{(ram)}$}
\put(1.5,6.8){$\scriptstyle \frac{d}{p}$}
\put(-1.2,6.8){$\scriptstyle \text{(un)}$}
\put(1.5,3){$\scriptstyle f(M/F)=\frac{pf}{d}$}
\put(1.5,1.5){$\scriptstyle e(M/F)=pe_0$}
\put(1,1){\line(0,1){3}}
\put(1,6){\line(0,1){2}}
\put(1,10){\line(0,1){2}}
\put(1,14){\line(0,1){2}}
\end{picture}
\end{equation*}
such that $K/M$ is Galois. The cardinality of ${\mathcal X}$ can be calculated
by counting the elements $(M,E,K)\in{\mathcal X}$ in the order of $M,E,K$;
when $M$ and $E$ are fixed, the number of $K$ is $\bigl|{\mathcal H}
\bigl(\frac{d}{p},\frac1d p^2e_0fn;p\bigr)\bigr|-1$, which is
computed in (\ref{kTMp}). It turns out that
\begin{equation} \label{cardX}
|{\mathcal X}|
=pe_0(p^{1+\frac 1d pe_0fn}-p+1)
\Bigl[(p-1,d)\frac{p^{\frac 1d p^2e_0fn}-1}{p-1}+p^{\frac 1d p^2e_0fn}-1\Bigr].
\end{equation}
On the other hand,
${\mathcal X}={\mathcal X}_1\overset{\cdot}\cup{\mathcal X}_2$,
where
\begin{align}
{\mathcal X}_1&=\{(M,E,K)\in{\mathcal X}:K\in{\mathcal C}_d^1\setminus{\mathcal C}_d^2\}, \\
{\mathcal X}_2&=\{(M,E,K)\in{\mathcal X}:K\in{\mathcal C}_d^1\cap{\mathcal C}_d^2\}.
\end{align}
%
%%%%%%%%%%%%%%%%%%%%%%%%%%%%%%%%
%  Lemma 8.1
%%%%%%%%%%%%%%%%%%%%%%%%%%%%%%%%
\begin{lem} \label{bijection1}
For each $K\in{\mathcal C}_d^1\setminus{\mathcal C}_d^2$, there
are unique subfields $M\subset E\subset K$ such that
$(M,E,K)\in{\mathcal X}$.
\end{lem}

\begin{proof}
The existence of such an $(M,E)$ follows from the definition of ${\mathcal C}_d^1$.
To see the uniqueness of $(M,E)$, we assume that there are
$(M_1,E_1)$ and $(M_2,E_2)$ such that $(M_1,E_1,K)\in{\mathcal X}$ and
$(M_2,E_2,K)\in{\mathcal X}$.  Then
$K/E_1\cap E_2$ is totally ramified, and also Galois since both $K/E_1$
and $K/E_2$ are Galois. Therefore by Proposition~\ref{consequences}, $[K:E_1\cap E_2]\mid p^2$.
If $E_1\ne E_2$, we must have $[K:E_1\cap E_2]=p^2$. However, this would imply that
$K\in{\mathcal C}_d^2$, contrary to assumption. Thus $E_1=E_2$.
In the diagram
\setlength{\unitlength}{3mm}
\begin{equation*}
\begin{picture}(8,12)
\put(4,0){\makebox(0,0){$M_1\cap M_2$}}
\put(0,4){\makebox(0,0){$M_1$}}
\put(8,4){\makebox(0,0){$M_2$}}
\put(4,8){\makebox(0,0){$E_1$}}
\put(4,12){\makebox(0,0){$K$}}
\put(5,1){\line(1,1){2}}
\put(1,3){\line(1,-1){2}}
\put(1,5){\line(1,1){2}}
\put(5,7){\line(1,-1){2}}
\put(4,9){\line(0,1){2}}
\put(0,6.2){$\scriptstyle \text{(un)}$}
\put(6.2,6.2){$\scriptstyle \text{(un)}$}
\put(2.2,5.2){$\scriptstyle \frac{d}{p}$}
\put(5.3,5.2){$\scriptstyle \frac{d}{p}$}
\put(3.2,9.8){$\scriptstyle p$}
\put(4.2,9.8){$\scriptstyle \text{(ram)}$}
\end{picture}
\end{equation*}
the extension $K/M_1\cap M_2$ is Galois since both $K/M_1$ and $K/M_2$ are Galois.
If $e(K/M_1\cap M_2)>p$ then $e(K/M_1\cap M_2)=p^2$,
which implies $K\in{\mathcal C}_d^2$, contrary to assumption.
Thus $E_1/M_1\cap M_2$ is unramified, and hence $M_1=M_2$.
\end{proof}

Lemma~\ref{bijection1} implies that $(M,E,K)\mapsto K$ gives a
bijection between ${\mathcal X}_1$
and ${\mathcal C}_d^1\setminus{\mathcal C}_d^2$. Hence
\begin{equation} \label{sizeC1dC2d}
|{\mathcal C}_d^1\setminus{\mathcal C}_d^2|=|{\mathcal X}_1|
=|{\mathcal X}|-|{\mathcal X}_2|.
\end{equation}
In order to calculate $|{\mathcal X}_2|$, we let ${\mathcal Y}$ be the set of all $(T,L,M,E,K)$ in the diagram
\setlength{\unitlength}{3mm}
\begin{equation*}
\begin{picture}(8,21)
\put(4,0){\makebox(0,0){$F$}}
\put(4,5){\makebox(0,0){$T$}}
\put(0,9){\makebox(0,0){$L$}}
\put(8,9){\makebox(0,0){$M$}}
\put(4,13){\makebox(0,0){$E$}}
\put(4,17){\makebox(0,0){$K$}}
\put(4,21){\makebox(0,0){$\Omega$}}
\put(4,1){\line(0,1){3}}
\put(5,6){\line(1,1){2}}
\put(1,8){\line(1,-1){2}}
\put(1,10){\line(1,1){2}}
\put(5,12){\line(1,-1){2}}
\put(4,14){\line(0,1){2}}
\put(4,18){\line(0,1){2}}
\put(4.5,1.4){$\scriptstyle e(T/F)=e_0$}
\put(4.5,3){$\scriptstyle f(T/F)=\frac{pf}d$}
\put(2.2,7.4){$\scriptstyle \frac dp$}
\put(0,6.5){$\scriptstyle\text{(un)}$}
\put(5.6,7.5){$\scriptstyle p$}
\put(2.2,10.4){$\scriptstyle p$}
\put(6.2,11.2){$\scriptstyle\text{(un)}$}
\put(5.5,10){$\scriptstyle\frac dp$}
\put(3,14.8){$\scriptstyle p$}
\end{picture}
\end{equation*}
such that $T=L\cap M$ and $K/T$ is Galois. Put
\begin{align}
{\mathcal Z}_1&=\{(T,L,M,E,K)\in{\mathcal Y}: E/L\ \text{is unramified}\}, \\
{\mathcal Z}_2&=\{(T,L,M,E,K)\in{\mathcal Y}: E/L\ \text{is ramified but}\ K/E\
\text{is unramified}\}.
\end{align}
%
%%%%%%%%%%%%%%%%%%%%%%%%%%%%
%   Lemma 8.2
%%%%%%%%%%%%%%%%%%%%%%%%%%%%
\begin{lem} \label{bijection2}
The map
\begin{equation}
\begin{matrix}
\eta:&{\mathcal Y}\setminus({\mathcal Z}_1\cup{\mathcal Z}_2)&\longrightarrow&{\mathcal X}_2\cr
&(T,L,M,E,K)&\longmapsto&(M,E,K)\cr
\end{matrix}
\end{equation}
is a bijection.
\end{lem}

\begin{proof}
For each $(T,L,M,E,K)\in{\mathcal Y}\setminus
({\mathcal Z}_1\cup{\mathcal Z}_2)$, the extension $K/L$
is totally ramified of degree $p^2$. By
Proposition~\ref{consequences}(i), $L$ is uniquely determined
by $K$. Consequently, $T=L\cap M$ is determined by $K$ and $M$,
so $\eta$ is one-to-one.

On the other hand, for each $(M,E,K)\in{\mathcal X}_2$, we have a diagram
\setlength{\unitlength}{3mm}
\begin{equation*}
\begin{picture}(5,8)
\put(5,0){\makebox(0,0){$M$}}
\put(1,4){\makebox(0,0){$E$}}
\put(1,8){\makebox(0,0){$K$}}
\put(2,3){\line(1,-1){2}}
\put(1,5){\line(0,1){2}}
\put(0.2,5.8){$\scriptstyle p$}
\put(1.2,5.8){$\scriptstyle\text{(ram)}$}
\put(3.2,2.2){$\scriptstyle \text{(un)}$}
\put(2.2,1){$\scriptstyle\frac dp$}
\end{picture}
\end{equation*}
with $K/M$ Galois.  Let $L/F$ be the unique subextension
of $K/F$ such that $K/L$ is totally ramified  of degree $p^2$.
Since $K\in{\mathcal C}_d^2$, the extension $K/L$ is Galois.
Clearly, $L\subset E$, and in the diagram
\setlength{\unitlength}{3mm}
\begin{equation*}
\begin{picture}(8,12)
\put(4,0){\makebox(0,0){$L\cap M$}}
\put(0,4){\makebox(0,0){$L$}}
\put(8,4){\makebox(0,0){$M$}}
\put(4,8){\makebox(0,0){$E$}}
\put(4,12){\makebox(0,0){$K$}}
\put(5,1){\line(1,1){2}}
\put(1,3){\line(1,-1){2}}
\put(1,5){\line(1,1){2}}
\put(5,7){\line(1,-1){2}}
\put(4,9){\line(0,1){2}}
\put(2.2,5.5){$\scriptstyle p$}
\put(5.5,5.2){$\scriptstyle \frac dp$}
\put(-0.5,6.2){$\scriptstyle\text{(ram)}$}
\put(6.2,6.2){$\scriptstyle \text{(un)}$}
\put(3.2,9.8){$\scriptstyle p$}
\put(4.2,9.8){$\scriptstyle \text{(ram)}$}
\put(2.2,2.3){$\scriptstyle \frac dp$}
\put(5.6,2.3){$\scriptstyle p$}
\put(6.2,1.5){$\scriptstyle\text{(ram)}$}
\put(0,1.5){$\scriptstyle \text{(un)}$}
\end{picture}
\end{equation*}
the extension $K/L\cap M$ is Galois. Hence $(L\cap M, L, M, E, K)\in{\mathcal Y}\setminus({\mathcal Z}_1\cup
{\mathcal Z}_2)$. This proves that $\eta$ is onto.
\end{proof}

It follows from Lemma~\ref{bijection2} that
\begin{equation} \label{X2YZ1Z2}
|{\mathcal X}_2|=|{\mathcal Y}|-|{\mathcal Z}_1|-|{\mathcal Z}_2|.
\end{equation}
The cardinality of ${\mathcal Z}_1$ can be calculated by counting
the elements ${(T,L,M,E,K)}$ in the order
$T,E,K,L,M$.  Fix $T$ and $E$.  Then the number of $K$ is
$\left|{\mathcal H}\left(d,\frac1d pe_0fn;p\right)\right|$, and
the number of $(L,M)$ is
0 if $(p,d/p)\ne 1$, and 1 if $(p,d/p)=1$.  It follows that
\begin{equation} \label{Z1cases}
|{\mathcal Z}_1|=
\begin{cases}
\displaystyle
e_0\biggl[(p-1,d)\frac{p^{\frac 1d pe_0fn}-1}{p-1}+p^{\frac 1d pe_0fn}\biggr]&
\text{if}\ p^2\nmid d,\cr
0&\text{if}\ p^2\mid d.\cr
\end{cases}
\end{equation}
%
%%%%%%%%%%%%%%%%%%%%%%%%%%%%%%
%    Lemma 8.3
%%%%%%%%%%%%%%%%%%%%%%%%%%%%%%
\begin{lem} \label{cardYlemma}
We have
\begin{equation} \label{cardY}
\begin{split}
|{\mathcal Y}|
=\begin{cases}
\left.\begin{split}
&e_0\bigg[(p-1,d)^2\cdot\frac {p(p^{\frac 1d pe_0fn}-1)^2}{(p-1)^2}+p^{\frac 1d pe_0fn}\cr
&\;\;\;+(p-1,d)\frac{(p^{\frac 1d pe_0fn}-1)(p^{1+\frac 1d pe_0fn}+1)}{p-1}
\bigg]\cr
\end{split}
\right\}&\text{if } p^2\nmid d,\cr
\cr
\left.\begin{split}
&e_0\bigg[(p-1,d)^2\cdot\frac{p(p^{\frac 1d pe_0fn}-1)^2}{(p-1)^2}
\cr
&\;\;\;+2(p-1,d)\frac{p^{1+\frac 1d pe_0fn}(p^{\frac 1d pe_0fn}-1)}{p-1}\cr
&\;\;\;+p^{1+\frac 1d pe_0fn}(p^{\frac 1d pe_0fn}-1)
\cr
&\;\;\;-\frac{(p^{\frac 1d pe_0fn}-1)(p^{\frac 1d pe_0fn}-p)}{p^2-1}\bigg]\cr
\end{split}
\right\}&\text{if }p^2\mid d.\cr
\end{cases}
\end{split}
\end{equation}
\end{lem}

\begin{proof}
First fix $T\in{\mathcal E}\left(F,\frac{1}{d}pf,e_0\right)$
and let $L/T$ be unramified of degree $d/p$.  Put
${\mathcal H}={\mathcal H}\bigl(\frac dp,\frac 1d pe_0fn;p^2\bigr)$,
and let $H\mapsto K_H$ be the bijection between ${\mathcal H}$
and ${\mathfrak K}(T,L;p^2)$ induced by class field theory.
For each $H\in{\mathcal H}$, the number of $(M,E)$ such that
$(T,L,M,E,K_H)\in{\mathcal Y}$ is equal to
$|{\mathcal J}(H)|$, where
\begin{equation}
{\mathcal J}(H)=\{J\leq\text{Gal}(K_H/T):|J|=d,\;\text{Gal}(K_H/T)=
\text{Gal}(K_H/L)\cdot J\}.
\end{equation}
Alternatively, $|{\mathcal J}(H)|$ is the number of
subextensions $M/T$ of $K_H/T$ such that $[M:T]=p$ and
$L\cap M=T$. For any such $M$, the compositum $E=LM$ is the
only field such that $(T,L,M,E, K_H)\in{\mathcal Y}$. Therefore
\begin{equation} \label{MEK}
|\{(M,E,K):(T,L,M,E,K)\in{\mathcal Y}\}|=\sum_{H\in{\mathcal H}}|{\mathcal J}(H)|.
\end{equation}

For each $H\in{\mathcal H}$, the structure of $\text{Gal}(K_H/T)$
and the position of $\text{Gal}(K_H/L)$ in $\text{Gal}(K_H/T)$
are explicitly described in Corollary~\ref{long}.
In fact, in the notation of Corollary~\ref{long}(i),
\begin{equation}
{\mathcal J}(H(\lambda,a))=\{\langle\boldsymbol{\kappa}^p,\boldsymbol{\kappa}^c\boldsymbol{\theta}\rangle:
c\in{\mathbb Z}/p^2{\mathbb Z},\;
(\boldsymbol{\kappa}^c\boldsymbol{\theta})^{d/p}\in
\langle\boldsymbol{\kappa}^p\rangle\},
\end{equation}
and
$\langle\boldsymbol{\kappa}^p,\boldsymbol{\kappa}^c\boldsymbol{\theta}\rangle=
\langle\boldsymbol{\kappa}^p,\boldsymbol{\kappa}^{c'}\boldsymbol{\theta}\rangle$
if and only if $c'-c\in p{\mathbb Z}/p^2{\mathbb Z}$. Therefore
\begin{equation}
\begin{split}
|{\mathcal J}(H(\lambda,a))|&
=\frac 1p|\{c\in{\mathbb Z}/p^2{\mathbb Z}:(\boldsymbol{\kappa}^c\boldsymbol{\theta})^{d/p}\in\langle\boldsymbol{\kappa}^p\rangle\}|\cr
&=\frac 1p\bigl|\bigl\{c\in{\mathbb Z}/p^2{\mathbb Z}:(1+\lambda+\cdots+\lambda^{\frac dp-1})c+\genfrac{[}{]}{0pt}{1}1\alpha^ta
\in p{\mathbb Z}/p^2{\mathbb Z}\bigr\}\bigr|\cr
&=
\begin{cases}
p&\text{if}\ \lambda\ne 1,\cr
1&\text{if}\ \lambda = 1,\ p^2\nmid d,\cr
p&\text{if}\ \lambda=1,\ p^2\mid d,\text{ and }
\genfrac{[}{]}{0pt}{1}1\alpha^ta
 \in p{\mathbb Z}/p^2{\mathbb Z},\cr
0&\text{if}\ \lambda=1,\ p^2\mid d,\text{ and }
\genfrac{[}{]}{0pt}{1}1\alpha^ta
\notin p{\mathbb Z}/p^2{\mathbb Z}.\cr
\end{cases}
\end{split}
\end{equation}

Suppose instead we are in the situation of Corollary~\ref{long}(ii).
Then $H=H(\Lambda,A)$ for some $\Lambda$ and $A$ satisfying the
conditions of the corollary.  For each
$J\in {\mathcal J}(H(\Lambda,A))$ the group
$B=J\cap\langle\boldsymbol{\kappa}_1,\boldsymbol{\kappa}_2\rangle$
is a normal subgroup of $J$ of order $p$, and
$\boldsymbol{\kappa}_1^{c_1}\boldsymbol{\kappa}_2^{c_2}\boldsymbol{\theta}\in J$ for some
$(c_1,c_2)\in({\mathbb Z}/p{\mathbb Z})^2$.
It follows that $B$ is
invariant under conjugation by $\boldsymbol{\theta}$, so
$B=\langle\boldsymbol{\kappa}_1^{b_1}\boldsymbol{\kappa}_2^{b_2}\rangle$
with $\genfrac{[}{]}{0pt}{1}{b_1}{b_2}$ an eigenvector
of $\Lambda^t$.
Therefore $J=
\langle\boldsymbol{\kappa}_1^{b_1}\boldsymbol{\kappa}_2^{b_2},
\boldsymbol{\kappa}_1^{c_1}\boldsymbol{\kappa}_2^{c_2}\boldsymbol{\theta}\rangle$,
with
$(\boldsymbol{\kappa}_1^{c_1}\boldsymbol{\kappa}_2^{c_2}\boldsymbol{\theta})^{d/p}\in B$.
Moreover, $\langle \boldsymbol{\kappa}_1^{b_1}\boldsymbol{\kappa}_2^{b_2}, \boldsymbol{\kappa}_1^{c_1}\boldsymbol{\kappa}_2^{c_2}\boldsymbol{\theta}\rangle=\langle \boldsymbol{\kappa}_1^{b_1}\boldsymbol{\kappa}_2^{b_2}, \boldsymbol{\kappa}_1^{c_1'}\boldsymbol{\kappa}_2^{c_2'}\boldsymbol{\theta}\rangle$ if
and only if $\boldsymbol{\kappa}_1^{c_1'-c_1}\boldsymbol{\kappa}_2^{c_2'-c_2}\in B$. Hence
\begin{equation}
\begin{split}
|{\mathcal J}(H(\Lambda,A))|
=\frac{1}{p}&\bigl|\bigl\{(W,(c_1,c_2)):W\text{ is a 1-dimensional eigenspace of $\Lambda^t$} \\
&\;\;\;\text{and }(\boldsymbol{\kappa}_1^{c_1}\boldsymbol{\kappa}_2^{c_2}\boldsymbol{\theta})^{d/p}\in \langle\boldsymbol{\kappa}_1^{b_1}\boldsymbol{\kappa}_2^{b_2}\rangle\bigr\}\bigr|\cr
=\frac{1}{p}&\Bigl|\Bigl\{(W,(c_1,c_2)):W\text{ is a 1-dimensional
eigenspace of $\Lambda^t$}\cr
&\;\;\;\text{and }(c_1,c_2)(I_2+\Lambda+\cdots+\Lambda^{\frac dp-1})+\genfrac{[}{]}{0pt}{1}1\alpha^tA\in W \Bigr\}\Bigr|.
\end{split}
\end{equation}
Thus $|{\mathcal J}(H(\Lambda,A))|$ can be determined from the canonical form of $\Lambda$.
This allows us to compute $\sum_{H\in{\mathcal H}}|{\mathcal J}(H)|$; we omit
the details. Since $\sum_{H\in{\mathcal H}}|{\mathcal J}(H)|$
is independent of the choice of $(T,L)$, by (\ref{MEK}) we
have
\begin{equation} \label{Z2easy}
|{\mathcal Y}|=\Bigl|{\mathcal E}\Bigl(F,\frac{pf}{d},e_0\Bigr)\Bigr|
\cdot\sum_{H\in{\mathcal H}}|{\mathcal J}(H)|,
\end{equation}
and the formula (\ref{cardY}) for $|{\mathcal Y}|$ follows.
\end{proof}

%%%%%%%%%%%%%%%%%%%%%%%%%%%%%%%
%   Lemma 8.4
%%%%%%%%%%%%%%%%%%%%%%%%%%%%%%%
\begin{lem}
We have
\begin{equation} \label{sizeZ2}
|{\mathcal Z}_2|=e_0(p-1,d)\frac{p(p^{\frac 1d pe_0fn}-1)}{p-1}.
\end{equation}
\end{lem}

\begin{proof}
Fix $T\in{\mathcal E}\left(F,\frac{1}{d}pf,e_0\right)$, let $L/T$
be unramified of degree $d/p$, and let $L'/L$ be unramified of
degree $p$.  Let $H\mapsto K_H$ be the bijection between
${\mathcal H}={\mathcal H}\left(d,\frac 1d pe_0fn;p\right)$ and
${\mathfrak K}(T,L';p)$ induced by class field theory.
For each $H\in{\mathcal H}$, put
\begin{equation}
\begin{split}
{\mathcal J}(H)=\{&(B,J):B\leq J\leq \text{Gal}(K_H/T),\ |B|=p,\ |J|=d, \\
&\hspace{1.5cm}\text{and Gal}(K_H/T)=\text{Gal}(K_H/L')\cdot J\}.\cr
\end{split}
\end{equation}
When $K_H/L'$ is ramified, the number of $(M,E)$ such that
$(T,L,M,E,K_H)\in{\mathcal Z}_2$ is equal to $|{\mathcal J}(H)|$.
In the following diagram we have $T=L'\cap M$,
%%%%%%%%%%%%%%%%%%%%%%%%
%  (??)
%%%%%%%%%%%%%%%%%%%%%%%%
and hence
${\text{Gal}(K_H/T)=\text{Gal}(K_H/L')\cdot \text{Gal}(K_H/M)}$.
\setlength{\unitlength}{3mm}
\begin{equation*}
\begin{picture}(8,12)
\put(4,0){\makebox(0,0){$T$}}
\put(0,4){\makebox(0,0){$L$}}
\put(0,8){\makebox(0,0){$L'$}}
\put(4,8){\makebox(0,0){$E$}}
\put(8,4){\makebox(0,0){$M$}}
\put(4,12){\makebox(0,0){$K_H$}}
\put(1,3){\line(1,-1){2}}
\put(0,5){\line(0,1){2}}
\put(1,9){\line(1,1){2}}
\put(5,1){\line(1,1){2}}
\put(1,5){\line(1,1){2}}
\put(5,7){\line(1,-1){2}}
\put(4,9){\line(0,1){2}}
\put(2.2,2.4){$\scriptstyle \frac dp$}
\put(0,1.5){$\scriptstyle \text{(un)}$}
\put(-0.5,10){$\scriptstyle \text{(ram)}$}
\put(2.2,9.5){$\scriptstyle p$}
\put(3.2,9.8){$\scriptstyle p$}
\put(-2,5.8){$\scriptstyle \text{(un)}$}
\put(0.2,5.8){$\scriptstyle p$}
\put(2.2,5.5){$\scriptstyle p$}
\put(5.6,2.2){$\scriptstyle p$}
\put(5.5,5){$\scriptstyle \frac dp$}
\put(4.2,9.8){$\scriptstyle \text{(un)}$}
\put(0.4,6.8){$\scriptstyle \text{(ram)}$}
\put(6,6.3){$\scriptstyle \text{(un)}$}
\put(6.2,1.5){$\scriptstyle \text{(ram)}$}
\end{picture}
\end{equation*}
Therefore
\begin{equation} \label{MEKformula}
|\{(M,E,K):(T,L,M,E,K)\in{\mathcal Z}_2\}|\;
=\sum_{\substack{H\in{\mathcal H}\\
K_H/L'\ \text{ramified}}}|{\mathcal J}(H)|.
\end{equation}
Corollary~\ref{explicit} allows us to compute $|{\mathcal J}(H)|$
for each $H\in{\mathcal H}$, and hence to compute
(\ref{MEKformula}).  Omitting the details, we get
\begin{equation} \label{sumIH}
\sum_{\substack{H\in{\mathcal H}\\
K_H/L'\ \text{ramified}}}
|{\mathcal J}(H)|=
(p-1,d)\frac{p(p^{\frac 1de_0fn}-1)}{p-1}.
\end{equation}
Using (\ref{MEKformula}) we see that $|{\mathcal Z}_2|$ is
equal to (\ref{sumIH}) multiplied by the number of pairs
$(T,L)$, which is
$\left|{\mathcal E}\left(F,\frac{1}{d}pf,e_0\right)\right|=e_0$.
\end{proof}

By (\ref{sizeC1dC2d}) and (\ref{X2YZ1Z2}) we have
\begin{equation}
|{\mathcal C}_d^1\setminus{\mathcal C}_d^2|=|{\mathcal X}|-|{\mathcal Y}|
+|{\mathcal Z}_1|+|{\mathcal Z}_2|,
\end{equation}
where $|{\mathcal X}|$, $|{\mathcal Y}|$, $|{\mathcal Z}_1|$, and $|{\mathcal Z}_2|$
are given in (\ref{cardX}), (\ref{cardY}), (\ref{Z1cases}), and
(\ref{sizeZ2}). Hence we have a formula
for $|{\mathcal C}_d^1\setminus{\mathcal C}_d^2|$.

%%%%%%%%%%%%%%%%%%%%%%%%%%%%%%%%%
%   Proposition 8.5
%%%%%%%%%%%%%%%%%%%%%%%%%%%%%%%%%
\begin{prop} \label{Cd1Cd2}
Let $d'=\text{\rm lcm}(p,d)$. Then we have
\begin{equation*}
\begin{split}
&|{\mathcal C}_d^1\setminus{\mathcal C}_d^2|=
\begin{cases}
0&\text{if}\ d\nmid pf,\cr
\cr
\left.\begin{split}
&e_0\biggl[-(p-1,d)^2\frac p{(p-1)^2}(p^{\frac 1{d'}pe_0fn}-1)^2 \\
&\;\;\;+(p-1,d)\frac p{p-1}\Bigl(p^{\frac 1{d'}p^2e_0fn}-1 \\
&\;\;\;+(p^{\frac 1{d'}pe_0fn}-1)
(p^{1+\frac 1{d'}p^2e_0fn}-p^{\frac 1{d'}pe_0fn}-p+1)\Bigr) \\
&\;\;\;+p(p^{1+\frac 1{d'}pe_0fn}-p+1)(p^{\frac 1{d'}p^2e_0fn}-1)\biggr]
\end{split} \right\} \hspace{-.3cm} &
\text{if}\ d\mid pf,\ p^2\nmid d,\cr
\cr
\left.\begin{split}
&e_0\biggl[-(p-1,d)^2\frac p{(p-1)^2}(p^{\frac 1{d'}pe_0fn}-1)^2 \\
&\;\;\;+(p-1,d)\frac p{p-1}\Bigl(-(p^{\frac 1{d'}pe_0fn}-1)
(2p^{\frac 1{d'}pe_0fn}-1)\cr
&\;\;\;+(p^{1+\frac 1{d'}pe_0fn}-p+1)
(p^{\frac 1{d'}p^2e_0fn}-1)\Bigr) \\
&\;\;\;+p(p^{1+\frac 1{d'}pe_0fn}-p+1)(p^{\frac 1{d'}p^2e_0fn}-1) \\
&\;\;\;-p^{1+\frac 1{d'}pe_0fn}(p^{\frac 1{d'}pe_0fn}-1)\cr
&\;\;\;+\frac{(p^{\frac 1{d'}pe_0fn}-1)(p^{\frac 1{d'}pe_0fn}-p)}{p^2-1}\biggr]
\end{split}
\right\} \hspace{-.3cm} &
\text{if}\ d\mid pf,\ p^2\mid d.\cr
\end{cases}
\end{split}
\end{equation*}
\end{prop}
%%%%%%%%%%%%%%%%%%%%%%%%%%%%%%%%%%%%%%%%%%%%%%%%%%
%    Section 9
%%%%%%%%%%%%%%%%%%%%%%%%%%%%%%%%%%%%%%%%%%%%%%%%%%
\section{Determination of $|{\mathcal C}_d^0\setminus({\mathcal C}_d^1\cup{\mathcal C}_d^2)|$}
\label{d012}

Recall that ${\mathcal C}_d^0\setminus({\mathcal C}_d^1\cup{\mathcal C}_d^2)$ consists
of those fields $K\in{\mathcal E}(F,f,e)$ such that ${d\bigm| |\text{Aut}(K/F)|}$ and
the ramification index of $K$ over
the fixed field $N_K$ of $\text{Aut}(K/F)$ is $1$.  Thus
$|{\mathcal C}_d^0\setminus({\mathcal C}_d^1\cup{\mathcal C}_d^2)|=0$ if $d\nmid f$,
so we will assume that $d\mid f$.

Let ${\mathcal X}$ be the set of all $(M,K)$ in the diagram
\setlength{\unitlength}{3mm}
\begin{equation*}
\begin{picture}(2,13)
\put(1,0){\makebox(0,0){$F$}}
\put(1,5){\makebox(0,0){$M$}}
\put(1,9){\makebox(0,0){$K$}}
\put(1,13){\makebox(0,0){$\Omega$}}
\put(1.5,6.8){$\scriptstyle d$}
\put(-1.2,6.8){$\scriptstyle \text{(un)}$}
\put(1.5,3){$\scriptstyle f(M/F)=\frac{f}d$}
\put(1.5,1.5){$\scriptstyle e(M/F)=p^2e_0$}
\put(1,1){\line(0,1){3}}
\put(1,6){\line(0,1){2}}
\put(1,10){\line(0,1){2}}
\end{picture}
\end{equation*}
Then
\begin{equation} \label{XEF}
\begin{split}
|{\mathcal X}|&=\Bigl|{\mathcal E}\Bigl(F,\frac fd, p^2e_0\Bigr)\Bigr| \\
&=e_0p^2(p^{2+\frac 1d(p+1)e_0fn}-p^{2+\frac 1d pe_0fn}+p^{1+\frac 1d pe_0fn}-p+1).
\end{split}
\end{equation}
On the other hand,
${\mathcal X}={\mathcal X}_1\overset{\cdot}{\cup}{\mathcal X}_2\overset{\cdot}{\cup}{\mathcal X}_3$,
where
\begin{align}
{\mathcal X}_1&=\{(M,K)\in{\mathcal X}:K\in{\mathcal C}_d^0\setminus({\mathcal C}_d^1\cup{\mathcal C}_d^2)\}, \\
{\mathcal X}_2&=\{(M,K)\in{\mathcal X}:K\in{\mathcal C}_d^0\cap({\mathcal C}_d^1\setminus{\mathcal C}_d^2)\},
\label{X2} \\
{\mathcal X}_3&=\{(M,K)\in{\mathcal X}:K\in{\mathcal C}_d^0\cap{\mathcal C}_d^2\}.
\end{align}

For each $K\in{\mathcal C}_d^0\setminus({\mathcal C}_d^1\cup{\mathcal C}_d^2)$, we claim that there is a unique $M$ such that $(M,K)\in{\mathcal X}$. The existence of $M$ follows from the definition of
${\mathcal C}_d^i$. To see the uniqueness of
$M$, assume to the contrary that we have two different
subextensions $M_1/F$ and $M_2/F$ of $K/F$ such that $K/M_1$
and $K/M_2$ are both unramified of degree $d$.
Then $K/M_1\cap M_2$ is Galois and not unramified. By Proposition~\ref{consequences}(iv),
we must have $e(K/M_1\cap M_2)=p$ or $p^2$. This means that $K\in{\mathcal C}_d^1\cup
{\mathcal C}_d^2$, which is a contradiction.
It follows that
\begin{equation} \label{XX2X3}
|{\mathcal C}_d^0\setminus({\mathcal C}_d^1\cup{\mathcal C}_d^2)|=|{\mathcal X}_1|=
|{\mathcal X}|-|{\mathcal X}_2|-|{\mathcal X}_3|.
\end{equation}

We now determine $|{\mathcal X}_3|$. Let ${\mathcal Y}$ be the
set of all $(T,L,M,K)$ in the diagram
\setlength{\unitlength}{3mm}
\begin{equation*}
\begin{picture}(8,17)
\put(4,0){\makebox(0,0){$F$}}
\put(4,5){\makebox(0,0){$T$}}
\put(0,9){\makebox(0,0){$L$}}
\put(8,9){\makebox(0,0){$M$}}
\put(4,13){\makebox(0,0){$K$}}
\put(4,17){\makebox(0,0){$\Omega$}}
\put(4,1){\line(0,1){3}}
\put(5,6){\line(1,1){2}}
\put(1,8){\line(1,-1){2}}
\put(1,10){\line(1,1){2}}
\put(5,12){\line(1,-1){2}}
\put(4,14){\line(0,1){2}}
\put(4.5,1.2){$\scriptstyle e(T/F)=e_0$}
\put(4.5,2.8){$\scriptstyle f(T/F)=\frac fd$}
\put(1.2,6.2){$\scriptstyle  d$}
\put(2,7.2){$\scriptstyle\text{(un)}$}
\put(6.2,6.5){$\scriptstyle p^2$}
\put(1.2,11.2){$\scriptstyle p^2$}
\put(4,10.5){$\scriptstyle\text{(un)}$}
\put(6.2,11.2){$\scriptstyle d$}
\end{picture}
\end{equation*}
such that $T=L\cap M$ and $K/T$ is Galois. Write
${\mathcal Y}={\mathcal Y}_1\overset{\cdot}{\cup}{\mathcal Y}_2$,
where
\begin{align}
{\mathcal Y}_1&=\{(T,L,M,K)\in{\mathcal Y}:e(K/L)=p^2\}, \\
{\mathcal Y}_2&=\{(T,L,M,K)\in{\mathcal Y}:e(K/L)=1\ \text{or}\ p\}.
\label{Y2}
\end{align}
Then $(T,L,M,K)\mapsto (M,K)$ is a bijection between ${\mathcal Y}_1$
and ${\mathcal X}_3$ (cf.\ Proposition~\ref{consequences}).
Therefore
\begin{equation} \label{X3Y1}
|{\mathcal X}_3|=|{\mathcal Y}_1|=|{\mathcal Y}|-|{\mathcal Y}_2|.
\end{equation}
%
%%%%%%%%%%%%%%%%%%%%%%%%%%%%%%%%%
%    Lemma 9.1
%%%%%%%%%%%%%%%%%%%%%%%%%%%%%%%%%
\begin{lem} \label{Ylemma}
We have
\begin{equation} \label{Ycases}
|{\mathcal Y}|=
\begin{cases}
\left.\begin{split}
&e_0\bigg[\frac 12 (p-1,d)^2\frac{p^2(p^{\frac 1d e_0fn}-1)^2}{(p-1)^2}
\cr
&\;\;\;+(p-1,d)\frac{p(p^{\frac 1d e_0fn}-1)}{p-1}\cr
&\;\;\;+\frac 12(p^2-1,d)\frac{p^2(p^{\frac 1d 2e_0fn}-1)}{p^2-1}+1\bigg]
\end{split}
\right\}&\text{if }p\nmid d,\cr
\cr
\left.\begin{split}
&e_0\bigg[\frac 12 (p-1,d)^2\frac{p^2(p^{\frac 1d e_0fn}-1)^2}{(p-1)^2}
\cr
&\;\;\;+(p-1,d)\frac{p^{2+\frac 1d e_0fn}(p^{\frac 1d e_0fn}-1)}{p-1}\cr
&\;\;\;+\frac 12(p^2-1,d)\frac{p^2(p^{\frac 1d 2e_0fn}-1)}{p^2-1}\bigg]
\end{split}
\right\}&\text{if }p\mid d.\cr
\end{cases}
\end{equation}
\end{lem}

\begin{proof}
The proof is similar to that of Lemma~\ref{cardYlemma}. Once
again, we only describe the method and omit the computational
details. Fix $T\in{\mathcal E}\bigl(F,\frac fd,e_0\bigr)$ and
let $L/T$ be the unramified extension
of degree $d$. Let $H\mapsto K_H$ be the bijection between
${\mathcal H}={\mathcal H}\left(d,\frac 1d e_0fn;p^2\right)$
and ${\mathfrak K}(T,L;p^2)$ induced by class field theory.
For each $H\in{\mathcal H}$ we have
\begin{equation} \label{cardIH}
|\{M:(T,L,M,K_H)\in{\mathcal Y}\}|=|{\mathcal J}(H)|
\end{equation}
where
\begin{equation}
{\mathcal J}(H)=\{J\leq\text{Gal}(K_H/T):|J|=d,\;\text{Gal}(K_H/T)=\text{Gal}(K_H/L)\cdot
J\}.
\end{equation}
Using (\ref{cardIH}), we see that
\begin{equation}
|\{(M,K):(T,L,M,K)\in{\mathcal Y}\}|=
\sum_{H\in{\mathcal H}}|{\mathcal J}(H)|.
\end{equation}
The sum $\sum_{H\in{\mathcal H}}|{\mathcal J}(H)|$ can be computed
using Corollary~\ref{long}, and the result is independent
of the choice of $(T,L)$ (cf.\ the proof of
Lemma~\ref{cardYlemma}). Therefore we get
$|{\mathcal Y}|=\bigl|{\mathcal E}\bigl(F,\frac fd,e_0\bigr)\bigr|\cdot
\sum_{H\in{\mathcal H}}|{\mathcal J}(H)|$.
\end{proof}

%%%%%%%%%%%%%%%%%%%%%%%%%%%%%%
%  Lemma 9.2
%%%%%%%%%%%%%%%%%%%%%%%%%%%%%%
\begin{lem} \label{Y2lemma}
We have
\begin{equation} \label{Y2cases}
|{\mathcal Y}_2|=
\begin{cases}
\displaystyle e_0\Bigl[(p-1,d)\frac{p(p^{\frac 1d e_0fn}-1)}{p-1}+1\Bigr]&\text{if}\ p\nmid d,\cr
0&\text{if}\ p\mid d.\cr
\end{cases}
\end{equation}
\end{lem}

\begin{proof}
By (\ref{Y2}) we have $|{\mathcal Y}_2|=0$ if $p\mid d$. Thus we
may assume $p\nmid d$.  Fix
$T\in {\mathcal E}\bigl(F,\frac fd,e_0\bigr)$, let $L/T$ be
unramified of degree $d$,
and let $L'/L$ be unramified of degree $p$. Let $H\mapsto K_H$ be the bijection
between ${\mathcal H}={\mathcal H}\left(dp,\frac 1d e_0fn;p\right)$ and
${\mathfrak K}(T,L';p)$ induced by class field theory. Then
for each $H\in{\mathcal H}$ we have
\begin{equation}
|\{M:(T,L,M,K_H)\in{\mathcal Y}_2\}|=|{\mathcal J}(H)|,
\end{equation}
where
\begin{equation}
{\mathcal J}(H)=\{J\leq\text{Gal}(K_H/T):|J|=d\}.
\end{equation}
It follows that
\begin{equation} \label{cardMKY2}
|\{(M,K):(T,L,M,K)\in{\mathcal Y}_2\}|
=\sum_{H\in{\mathcal H}}|{\mathcal J}(H)|,
\end{equation}
which can be computed using Corollary~\ref{explicit}. The
lemma then follows from (\ref{cardMKY2}) and the formula
$|{\mathcal Y}_2|=\bigl|{\mathcal E}\bigl(F,\frac fd,e_0\bigr)\bigr|\cdot
\sum_{H\in{\mathcal H}}|{\mathcal J}(H)|$.
\end{proof}

Since $|{\mathcal X}_3|$ is determined by
$|{\mathcal Y}|$ and $|{\mathcal Y}_2|$, the only part of
(\ref{XX2X3}) that remains to be computed is $|{\mathcal X}_2|$.
Let ${\mathcal Z}$ be the set of all $(T,L,M,K)$ in the diagram
\setlength{\unitlength}{3mm}
\begin{equation*}
\begin{picture}(8,17)
\put(4,0){\makebox(0,0){$F$}}
\put(4,5){\makebox(0,0){$T$}}
\put(0,9){\makebox(0,0){$L$}}
\put(8,9){\makebox(0,0){$M$}}
\put(4,13){\makebox(0,0){$K$}}
\put(4,17){\makebox(0,0){$\Omega$}}
\put(4,1){\line(0,1){3}}
\put(5,6){\line(1,1){2}}
\put(1,8){\line(1,-1){2}}
\put(1,10){\line(1,1){2}}
\put(5,12){\line(1,-1){2}}
\put(4,14){\line(0,1){2}}
\put(4.5,1.2){$\scriptstyle e(T/F)=pe_0$}
\put(4.5,2.8){$\scriptstyle f(T/F)=\frac fd$}
\put(1.2,6.2){$\scriptstyle  d$}
\put(2,7.2){$\scriptstyle\text{(un)}$}
\put(6.2,6.5){$\scriptstyle p$}
\put(1.5,11.2){$\scriptstyle p$}
\put(4,10.5){$\scriptstyle\text{(un)}$}
\put(6.5,11.2){$\scriptstyle d$}
\end{picture}
\end{equation*}
such that $T=L\cap M$ and $K/T$ is Galois. In addition, define
\begin{align}
{\mathcal Z}_1&=\{(T,L,M,K)\in{\mathcal Z}:e(K/L)=p\}, \\
{\mathcal Z}_2&=\{(T,L,M,K)\in{\mathcal Z}_1:K/L_K\ \text{is Galois}\},
\label{Z2}
\end{align}
where $L_K$ is the unique field between $F$ and $K$ such that
$K/L_K$ is totally ramified of degree $p^2$.  We claim that
\begin{equation} \label{Z1Z2}
\begin{matrix}
\psi:&{\mathcal Z}_1\setminus{\mathcal Z}_2&\longrightarrow&{\mathcal X}_2\cr
&(T,L,M,K)&\longmapsto&(M,K)\cr
\end{matrix}
\end{equation}
is a bijection. By the definitions of ${\mathcal X}_2$,
${\mathcal Z}_1$, and
${\mathcal Z}_2$, it is clear that $\psi$ is onto. Suppose that $\psi$ is
not one-to-one. Then there are elements $(T_1,L_1,M,K)$ and
$(T_2,L_2,M,K)$ of ${\mathcal Z}_1\setminus{\mathcal Z}_2$ such that
$(T_1,L_1)\ne (T_2,L_2)$. Since $T_1= L_1\cap M$ and $T_2=L_2\cap M$, we must have
$L_1\ne L_2$. Since $K/L_1$ and $K/L_2$ are Galois and totally ramified of degree $p$,
we see that $K/L_1\cap L_2$ is Galois with $e(K/L_1\cap L_2)>p$.
Thus by Proposition~\ref{consequences}(iv) we get
$e(K/L_1\cap L_2)=p^2$. This implies that
$(T_1,L_1,M,K)\in{\mathcal Z}_2$, which is a contradiction.
It follows that
\begin{equation} \label{X2Z1Z2}
|{\mathcal X}_2|=|{\mathcal Z}_1|-|{\mathcal Z}_2|.
\end{equation}
%
%%%%%%%%%%%%%%%%%%%%%%%%%%%%%%%%
%  Lemma 9.3
%%%%%%%%%%%%%%%%%%%%%%%%%%%%%%%
\begin{lem}
We have
\begin{equation} \label{Z1formula}
|{\mathcal Z}_1|
=e_0(p-1,d)\frac{p^2(p^{1+\frac 1d e_0fn}-p+1)(p^{\frac 1d pe_0fn}-1)}{p-1}.
\end{equation}
\end{lem}

\begin{proof}
Observe that
\begin{equation} \label{Z1Z}
|{\mathcal Z}_1|=|{\mathcal Z}|-|\{(T,L,M,K)\in{\mathcal Z}:
K/L\ \text{unramified}\}|,
\end{equation}
and that
\begin{equation} \label{TLMK}
\begin{split}
|\{(T,L,M,K)\in{\mathcal Z}:
K/L\ \text{unram.}\}|&=
\begin{cases}
\bigl|{\mathcal E}\bigl(F,\frac fd,pe_0\bigr)\bigr|&\text{if}\ p\nmid d,\cr
0&\text{if}\ p\mid d;\cr
\end{cases}\cr
&=
\begin{cases}
pe_0(p^{1+\frac 1d e_0fn}-p+1)&\text{if}\ p\nmid d,\cr
0&\text{if}\ p\mid d.\cr
\end{cases}\cr
\end{split}
\end{equation}
Meanwhile $|{\mathcal Z}|$ can be computed as before: Fix
$T\in{\mathcal E}\bigl(F,\frac fd,pe_0\bigr)$ and
let $L/T$ be unramified of degree $d$. Let $H\mapsto K_H$ be
the bijection between ${\mathcal H}=
{\mathcal H}\left(d,\frac 1d pe_0fn;p\right)$ and
${\mathfrak K}(T,L;p)$ induced by class field theory.
Then for each $H\in{\mathcal H}$ we have
\begin{equation}
|\{M:(T,L,M,K_H)\in{\mathcal Z}\}|=|{\mathcal J}(H)|,
\end{equation}
where
\begin{equation}
{\mathcal J}(H)=\{J\leq\text{Gal}(K_H/T):|J|=d,\;\text{Gal}(K_H/T)=\text{Gal}(K_H/L)\cdot J\}.
\end{equation}
Consequently,
\begin{align} \label{Zsum}
|\{(M,K):(T,L,M,K)\in{\mathcal Z}\}|&=
\sum_{H\in{\mathcal H}}|{\mathcal J}(H)| \\[.2cm]
|{\mathcal Z}|&=
\Bigl|{\mathcal E}\Bigl(F,\frac fd,pe_0\Bigr)\Bigr|\cdot
\sum_{H\in{\mathcal H}}|{\mathcal J}(H)|.
\end{align}
Using (\ref{Zsum}) and Corollary~\ref{explicit} we get
\begin{equation} \label{Zcases}
|{\mathcal Z}|=
\begin{cases}
\displaystyle e_0p(p^{1+\frac 1d e_0fn}-p+1)\Bigl[(p-1,d)\frac{p(p^{\frac 1d pe_0fn}-1)}{p-1}+1\Bigr]
&\text{if}\ p\nmid d,\\[.2cm]
\displaystyle e_0(p-1,d)\frac{p^2(p^{1+\frac 1d e_0fn}-p+1)(p^{\frac 1d pe_0fn}-1)}{p-1}&\text{if}\
p\mid d.\cr
\end{cases}
\end{equation}
Equation (\ref{Z1formula}) now follows from (\ref{Z1Z}),
(\ref{TLMK}), and (\ref{Zcases}).
\end{proof}

%%%%%%%%%%%%%%%%%%%%%%%%%%%%%%%
%  Lemma 9.4
%%%%%%%%%%%%%%%%%%%%%%%%%%%%%%%
\begin{lem}
We have
\begin{equation} \label{Z2formula}
|{\mathcal Z}_2|
=
\begin{cases}
\left.\begin{split}
&e_0\biggl[\frac 12(p-1,d)^2\frac{p^2(p+1)(p^{\frac 1d e_0fn}-1)^2}{(p-1)^2}\cr
&\;\;\;-(p-1,d)
\frac{p^{2+\frac 1d e_0fn}(p^{\frac 1d e_0fn}-1)}{p-1}\cr
&\;\;\;+\frac 12(p^2-1,d)\frac{p^2(p^{\frac 1d 2e_0fn}-1)}{p-1}\biggr]
\end{split}
\right\}&\text{if }p\nmid d,\cr
\cr
\left.\begin{split}
&e_0\biggl[\frac 12(p-1,d)^2\frac{p^2(p+1)(p^{\frac 1d e_0fn}-1)^2}{(p-1)^2}\cr
&\;\;\;+(p-1,d)\frac{p^{2+\frac 1d e_0fn}(p^{\frac 1d e_0fn}-1)}{p-1}\cr
&\;\;\;+\frac 12(p^2-1,d)\frac{p^2(p^{\frac 1d 2e_0fn}-1)}{p-1}\biggr]\quad
\end{split}
\right\}&\text{if }p\mid d.\cr
\end{cases}
\end{equation}
\end{lem}

\begin{proof}
For each $(T,L,M,K)\in{\mathcal  Z}_2$, we have a diagram
\setlength{\unitlength}{3mm}
\begin{equation*}
\begin{picture}(12,17)
\put(4,0){\makebox(0,0){$F$}}
\put(4,5){\makebox(0,0){$S$}}
\put(0,9){\makebox(0,0){$L'$}}
\put(4,13){\makebox(0,0){$L$}}
\put(8,9){\makebox(0,0){$T$}}
\put(12,13){\makebox(0,0){$M$}}
\put(8,17){\makebox(0,0){$K$}}
\put(4,1){\line(0,1){3}}
\put(5,6){\line(1,1){2}}
\put(9,10){\line(1,1){2}}
\put(1,10){\line(1,1){2}}
\put(5,14){\line(1,1){2}}
\put(1,8){\line(1,-1){2}}
\put(5,12){\line(1,-1){2}}
\put(9,16){\line(1,-1){2}}
\put(4.2,1.2){$\scriptstyle e(S/F)=e_0$}
\put(4.2,2.8){$\scriptstyle f(S/F)=\frac fd$}
\put(2.2,7){$\scriptstyle d$}
\put(0,6.2){$\scriptstyle\text{(un)}$}
\put(6.2,11){$\scriptstyle d$}
\put(4,10.2){$\scriptstyle\text{(un)}$}
\put(10.2,15){$\scriptstyle d$}
\put(8.2,14.2){$\scriptstyle\text{(un)}$}
\put(1.2,11){$\scriptstyle p$}
\put(1.6,10){$\scriptstyle\text{(ram)}$}
\put(5.2,7){$\scriptstyle p$}
\put(6,6){$\scriptstyle\text{(ram)}$}
\put(9.2,11){$\scriptstyle p$}
\put(9.6,10){$\scriptstyle\text{(ram)}$}
\put(5.2,15){$\scriptstyle p$}
\put(5.8,14){$\scriptstyle\text{(ram)}$}
\end{picture}
\end{equation*}
in which $K/S$ is Galois, and $S$ and $L'$ are determined by
$(T,L,M,K)$. Thus if we let ${\mathcal W}$ denote the set of
all $(S,L',T,L,M,K)$ in this diagram such that $K/S$ is Galois, we
have $|{\mathcal Z}_2|=|{\mathcal W}|$. To compute $|{\mathcal W}|$,
we fix $S\in{\mathcal E}\bigl(F,\frac fd,e_0\bigr)$ and let $L'/S$ be unramified of degree $d$ and
$N/L'$ unramified of degree $p$. Write ${\mathfrak K}={\mathfrak K}(S,L';p^2)$
and ${\mathfrak K}'=\{K\in{\mathfrak K}:N\subset K\}={\mathfrak K}(S,N;p)$.
For each $K\in{\mathfrak K}$, let
\begin{equation}
\begin{split}
{\mathcal G}(K)&=\{(B,J):B\leq\text{Gal}(K/L'),\;|B|=p,\;
J\leq\text{Gal}(K/S), \\
&\hspace{2cm}|J|=d,\;\text{Gal}(K/S)=\text{Gal}(K/L')\cdot J\}.
\end{split}
\end{equation}
Note that if $(S,L',T,L,M,K)\in{\mathcal W}$, then $K\in{\mathfrak K}\setminus{\mathfrak K}'$.
Also note that for each $K\in{\mathfrak K}\setminus{\mathfrak K}'$
we have
\begin{equation} \label{GK}
|\{(T,L,M):(S,L',T,L,M,K)\in{\mathcal W}\}|=|{\mathcal G}(K)|.
\end{equation}
More precisely, $(T,L,M)\leftrightarrow (\text{Gal}(K/L),\text{Gal}(K/M))$ is a bijection
between the two sets in (\ref{GK}). Therefore
\begin{equation} \label{TLMKformula}
|\{(T,L,M,K):(S,L',T,L,M,K)\in{\mathcal W}\}|\:
=\sum_{K\in{\mathfrak K}}|{\mathcal G}(K)|-\sum_{K\in{\mathfrak K}'}|{\mathcal G}(K)|.
\end{equation}

To compute $\sum_{K\in{\mathfrak K}}|{\mathcal G}(K)|$,
let $H\mapsto K_H$ be the bijection between ${\mathcal H}=
{\mathcal H}\left(d,\frac 1de_0fn;p^2\right)$ and ${\mathfrak K}$
induced by class field theory, and
let ${\mathcal H}_1$ and ${\mathcal H}_2$ be the subsets of ${\mathcal H}$
corresponding to the two cases of Corollary~\ref{long}.
We have $\text{Gal}(K_H/L')\cong{\mathbb Z}/p^2{\mathbb Z}$ if
$H\in{\mathcal H}_1$ and $\text{Gal}(K_H/L')\cong({\mathbb Z}/p{\mathbb Z})^2$ if
$H\in{\mathcal H}_2$.  Therefore
\begin{equation}
|{\mathcal G}(K_H)|=
\begin{cases}
|{\mathcal J}(H)|&\text{if}\ H\in{\mathcal H}_1,\cr
(p+1)|{\mathcal J}(H)|&\text{if}\ H\in{\mathcal H}_2,\cr
\end{cases}
\end{equation}
where
\begin{equation}
{\mathcal J}(H)=\{J\leq\text{Gal}(K_H/S):|J|=d,\;\text{Gal}(K_H/S)=\text{Gal}(K_H/L')\cdot J\}.
\end{equation}
Thus
\begin{equation}
\sum_{K\in{\mathfrak K}}|{\mathcal G}(K)|
=\sum_{H\in{\mathcal H}}|{\mathcal G}(K_H)|
=\sum_{H\in{\mathcal H}_1}|{\mathcal J}(H)|+
(p+1)\!\sum_{H\in{\mathcal H}_2}|{\mathcal J}(H)|
\end{equation}
can be computed as before.

We now compute $\sum_{K\in{\mathfrak K}'}|{\mathcal G}(K)|$. We claim that
if $p\mid d$, then $|{\mathcal G}(K)|=0$ for all $K\in{\mathfrak K}'$.
Suppose to the contrary that there exists $(B,J)\in{\mathcal G}(K)$
for some $K\in{\mathfrak K}'$. Let $M$ denote the subfield of $K$ fixed by
$J$. Then we have the diagram
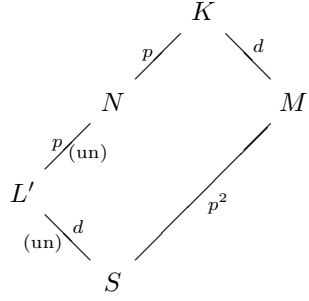
\setlength{\unitlength}{3mm}
\begin{equation*}
\begin{picture}(12,12)
\put(4,0){\makebox(0,0){$S$}}
\put(0,4){\makebox(0,0){$L'$}}
\put(4,8){\makebox(0,0){$N$}}
\put(12,8){\makebox(0,0){$M$}}
\put(8,12){\makebox(0,0){$K$}}
\put(1,3){\line(1,-1){2}}
\put(1,5){\line(1,1){2}}
\put(5,1){\line(1,1){6}}
\put(5,9){\line(1,1){2}}
\put(9,11){\line(1,-1){2}}
\put(0,1.5){$\scriptstyle\text{(un)}$}
\put(2,5.5){$\scriptstyle\text{(un)}$}
\put(2.2,2.2){$\scriptstyle d$}
\put(8.2,3.2){$\scriptstyle p^2$}
\put(1.3,6){$\scriptstyle p$}
\put(5.3,10){$\scriptstyle p$}
\put(10.2,10.2){$\scriptstyle d$}
\end{picture}
\end{equation*}
with $p\mid f(M/S)$. Thus $L'\cap M\ne S$, so $\text{Gal}(K/S)\ne
\text{Gal}(K/L')\cdot J$, contrary to the definition of
${\mathcal G}(K)$. Therefore we may
assume $p\nmid d$. Let $H\mapsto K_H$ be the bijection between
${\mathcal H}'={\mathcal H}\left(dp,\frac 1d e_0fn;p\right)$ and
${\mathfrak K}'$ induced by class field theory.
By Corollary~\ref{explicit}, every $H\in{\mathcal H}'$ is of
the form $H(\lambda,a)=a^\bot$ for some eigenvector $a$ of
$E\left(dp,\frac 1d e_0fn\right)^t$ with eigenvalue $\lambda$.
Furthermore, $\text{Gal}(K_{H(\lambda,a)}/S)$ is generated by
$\text{Gal}(K_{H(\lambda,a)}/N)=\langle\boldsymbol{\kappa}\rangle
\cong{\mathbb Z}/p{\mathbb Z}$
and an element $\boldsymbol{\theta}$ such that
$\boldsymbol{\theta}^{dp}=\boldsymbol{\kappa}^{c(a)}$ and
$\boldsymbol{\theta}\boldsymbol{\kappa}\boldsymbol{\theta}^{-1}=
\boldsymbol{\kappa}^{\lambda}$, where
$c(a)=\genfrac{[}{]}{0pt}{1}1\alpha^t\!a\in {\mathbb Z}/p{\mathbb Z}$
\,for some fixed $\alpha\in({\mathbb Z}/p{\mathbb Z})^{pe_0fn}$.
Moreover,
$\text{Gal}(K_{H(\lambda,a)}/L')$ is generated by
$\boldsymbol{\kappa}$ and $\boldsymbol{\theta}^d$
(cf.\ the proofs of Proposition~\ref{GalKHT} and Corollary~\ref{long}).
These explicit descriptions allow us to compute each $|{\mathcal G}
(K_{H(\lambda,a)})|$, and hence to compute $\sum_{H\in{\mathcal H}'}
|{\mathcal G}(K_H)|=\sum_{K\in{\mathfrak K}'}|{\mathcal G}(K)|$.
\vspace{.1cm}

Now $|\{(T,L,M,K):(S,L',T,L,M,K)\in{\mathcal W}\}|$ can be
computed using (\ref{TLMKformula}).
After multiplying the result by $\bigl|{\mathcal E}\bigl(F,\frac fd,e_0\bigr)\bigr|=e_0$,
we get the formula (\ref{Z2formula}) for $|{\mathcal W}|=|{\mathcal Z}_2|$.
\end{proof}

 From (\ref{XX2X3}), (\ref{X3Y1}), and (\ref{X2Z1Z2}), we have
\begin{equation}
|{\mathcal C}_d^0\setminus({\mathcal C}_d^1\cup{\mathcal C}_d^2)|
=|{\mathcal X}|-|{\mathcal Z}_1|+|{\mathcal Z}_2|-|{\mathcal Y}|+|{\mathcal Y}_2|,
\end{equation}
where $|{\mathcal X}|$, $|{\mathcal Y}|$, $|{\mathcal Y}_2|$, $|{\mathcal Z}_1|$, and $|{\mathcal Z}_2|$
are given in (\ref{XEF}), (\ref{Ycases}), (\ref{Y2cases}), (\ref{Z1formula}), and (\ref{Z2formula}).
Thus we obtain the main result of this section.

%%%%%%%%%%%%%%%%%%%%%%%%%%%%%%%%%%%%
%  Proposition 9.5
%%%%%%%%%%%%%%%%%%%%%%%%%%%%%%%%%%%%
\begin{prop} \label{Cd0Cd1Cd2}
We have
\begin{equation*}
\begin{split}
&|{\mathcal C}_d^0\setminus({\mathcal C}_d^1\cup{\mathcal C}_d^2)|=
\begin{cases}
0&\text{if}\ d\nmid f,\cr
\cr
\left. \begin{split}
&e_0p^2\Bigl[\frac 12(p-1,d)^2\frac{p(p^{\frac 1d e_0fn}-1)^2}{(p-1)^2} \\
&\;\;\;+\frac 12 (p^2-1,d)\frac{p(p^{\frac 1d 2e_0fn}-1)}{p^2-1}\cr
&\;\;\; -\frac{(p-1,d)}{p-1}
\Bigl(p^{\frac 1d e_0fn}(p^{\frac 1d e_0fn}-1) \\
&\;\;\;+(p^{1+\frac 1d e_0fn}-p+1)(p^{\frac 1d pe_0fn}-1)
\Bigr)\cr
&\;\;\; +p^{2+\frac 1d(p+1)e_0fn}-p^{2+\frac 1d pe_0fn} \\
&\;\;\;+p^{1+\frac 1d pe_0fn}-p+1\Bigr]
\end{split} \right\}&
\text{if}\ d\mid f,\ p\nmid d,\cr
\cr
\left. \begin{split}
&e_0p^2\Bigl[\frac 12(p-1,d)^2\frac{p(p^{\frac 1d e_0fn}-1)^2}{(p-1)^2} \\
&\;\;\;+\frac 12 (p^2-1,d)\frac{p(p^{\frac 1d 2e_0fn}-1)}{p^2-1}\cr
&\;\;\; -(p-1,d)\frac{(p^{\frac 1d pe_0fn}-1)(p^{1+\frac 1d e_0fn}-p+1)}{p-1}\cr
&\;\;\; +p^{2+\frac 1d(p+1)e_0fn}-p^{2+\frac 1d pe_0fn} \\
&\;\;\;+p^{1+\frac 1d pe_0fn}-p+1\Bigr]
\end{split} \right\} \hspace{-.3cm}
&\text{if}\ d\mid f,\ p\mid d.\cr
\end{cases} \\
\end{split}
\end{equation*}
\end{prop}
%%%%%%%%%%%%%%%%%%%%%%%%%%%%%%%%%%%%%%%%%%%%%%%%%%
%    Section 10
%%%%%%%%%%%%%%%%%%%%%%%%%%%%%%%%%%%%%%%%%%%%%%%%%%
\section{Conclusion of the Case $p^2\parallel e$} \label{conclude}

With $|{\mathcal C}_d^2|$, $|{\mathcal C}_d^1\setminus{\mathcal C}_d^2|$, and
$|{\mathcal C}_d^0\setminus({\mathcal C}_d^1\cup{\mathcal C}_d^2)|$ computed
in Sections \ref{d2} -- \ref{d012}, we are ready to state our
main result in the case $p^2\parallel e$.

%%%%%%%%%%%%%%%%%%%%%%%%
%   Theorem 10.1
%%%%%%%%%%%%%%%%%%%%%%%%
\begin{thm} \label{main2}
Let $F/{\mathbb Q}_p$ be a finite extension of degree $n$ with
residue degree $f_1$ and ramification index $e_1$. Let $f$, $e$,
and $e_0$ be positive integers such that $e=p^2e_0$,
$p\nmid e_0$, $(p^{f_1f}-1,e_0)=1$, and either $f(F(\zeta_p)/F)\nmid f$ or $e(F(\zeta_p)/F)>1$.  Write $f=p^it$ with $p\nmid t$. Then the number of
$F$-isomorphism classes of finite extensions $K/F$ with residue degree $f$ and
ramification index $e$ is given by
\begin{equation*}
\begin{split}
{\mathfrak I}(F,f,e)
=&\,\frac 1{p^it}\sum_{\tau\mid t}\phi(\tau)\Biggl[
\frac 12(p-1,\tau)^2\frac{(p^{\frac t\tau p^ie_0n}-1)^2}{p-1}
\\
&\;\;\;-(p-1,\tau)\frac{p^{\frac t\tau p^ie_0n}(p^{\frac t\tau p^ie_0n}-1)}{p-1}
+\frac 12 (p^2-1,\tau)\frac{(p^{\frac t\tau 2p^ie_0n}-1)}{p-1}
\\
&\;\;\;+(p+1)(p^{1+\frac t\tau (p+1) p^ie_0n}-p^{\frac t\tau p^ie_0n})
-(p^2-1)p^{\frac t\tau p^{i+1}e_0n}
+p^{\frac t\tau 2p^ie_0n}\cr
&\;\;\;+\sum_{j=1}^ip^{j-1}\biggl[
\frac{(p^{\frac t\tau p^{i-j}e_0n}-1)(p^{\frac t\tau p^{i-j}e_0n-1}-1)}{p+1} \\
&\hspace{1cm}+(p-1)(-p^{2+\frac t\tau p^{i-j+1}e_0n}+p^{\frac t\tau p^{i-j+1}e_0n}
+p^{\frac t\tau 2p^{i-j}e_0n})\cr
&\hspace{1cm}+(p^2-1)(p^{1+\frac t\tau (p+1)p^{i-j}e_0n}-p^{\frac t\tau p^{i-j}e_0n})
\\
&\hspace{1cm}+
\frac 12(p-1,\tau)^2(p^{\frac t\tau p^{i-j}e_0n}-1)^2
+\frac 12(p^2-1,\tau)(p^{\frac t\tau 2p^{i-j}e_0n}-1)
\biggr]\Biggr].
\end{split}
\end{equation*}
\end{thm}

\begin{proof}
By (\ref{Mobius}) and (\ref{Bd}) we have
\begin{equation}
{\mathfrak I}(F,f,e)=\frac 1{fe}\sum_{d>0}\phi(d)\bigl(n_0(d)+n_1(d)+n_2(d)\bigr),
\end{equation}
where
$n_0(d)=|{\mathcal C}_d^0\setminus({\mathcal C}_d^1\cup{\mathcal C}_d^2)|$,
$n_1(d)=|{\mathcal C}_d^1\setminus{\mathcal C}_d^2|$, and
$n_2(d)=|{\mathcal C}_d^2|$.
Since ${n_2(d)=n_2\bigl(\text{lcm}(p^2,d)\bigr)}$ and
$n_2(d)=0$ when $d\nmid p^2f$, it follows that
\begin{equation}
\begin{split}
\sum_{d>0}\phi(d)n_2(d)=&\sum_{d\mid p^{2+i}t}\phi(d)n_2(d)\cr
=\,&\sum_{\tau\mid t}\phi(\tau)\biggl(\phi(1)n_2(\tau)+\phi(p)n_2(p\tau) \\
&\hspace{2cm}
+\phi(p^2)n_2(p^2\tau)
+\sum_{j=3}^{2+i}\phi(p^j)n_2(p^j\tau)\biggr)\cr
=\,&\sum_{\tau\mid t}\phi(\tau)\biggl(p^2n_2(p^2\tau)+(p-1)
\sum_{j=3}^{2+i}p^{j-1}n_2(p^j\tau)\biggr).\cr
\end{split}
\end{equation}
By the same reasoning we have
\begin{align}
\sum_{d>0}\phi(d)n_0(d)&=\sum_{\tau\mid t}\phi(\tau)\biggl(n_0(\tau)+(p-1)
\sum_{j=1}^{i}p^{j-1}n_0(p^j\tau)\biggr), \\
\sum_{d>0}\phi(d)n_1(d)&=\sum_{\tau\mid t}\phi(\tau)\biggl(pn_1(p\tau)+(p-1)
\sum_{j=2}^{1+i}p^{j-1}n_1(p^j\tau)\biggr).
\end{align}
It follows that
\begin{equation}
\begin{split}
{\mathfrak I}(F,f,e)=&\frac 1{p^{2+i}te_0}\sum_{\tau\mid t}
\phi(\tau)\biggl(n_0(\tau)+pn_1(p\tau)+p^2n_2(p^2\tau)\cr
&\;\;+(p-1)\sum_{j=1}^ip^{j-1}\bigl(n_0(p^{j}\tau)+pn_1(p^{j+1}\tau)+
p^2n_2(p^{j+2}\tau)\bigr)\biggr).
\end{split}
\end{equation}
Using Propositions \ref{Cd2}, \ref{Cd1Cd2}, and
\ref{Cd0Cd1Cd2} to write out
$n_0(p^j\tau)+pn_1(p^{j+1}\tau)+p^2n_2(p^{j+2}\tau)$
explicitly for $0\le j \le i$, we obtain the final formula
for ${\mathfrak I}(F,f,e)$.
\end{proof}

%%%%%%%%%%%%%%%%%%%%%%%%%%%%%%%%%%%%%%%%%%%%%%%%%%
%    References
%%%%%%%%%%%%%%%%%%%%%%%%%%%%%%%%%%%%%%%%%%%%%%%%%%

\end{document}